\documentclass[10pt]{amsart}
\usepackage[T1]{fontenc}
\usepackage{geometry}
\usepackage[latin1] {inputenc}
\usepackage{amsmath}
\usepackage{amsfonts, amssymb, textcomp}
\usepackage[colorlinks=flase, linkcolor=red,urlcolor=green, citecolor=blue]{hyperref}
\usepackage{subeqnarray}

\usepackage{latexsym}
\usepackage{fancyhdr}
\usepackage{longtable}
\usepackage{amsmath, amssymb}
\usepackage{graphicx}
\usepackage{enumerate}

\numberwithin{equation}{section}

\theoremstyle{plain}
\newtheorem{proposition}{Proposition}[section]
\newtheorem{theorem}[proposition]{Theorem}
\newtheorem{lemma}[proposition]{Lemma}
\newtheorem{corollary}[proposition]{Corollary}
\newtheorem{definition}[proposition]{Definition}
\newtheorem{example}[proposition]{Example}

\newtheorem{remark}[proposition]{Remark}

\newcommand{\RR}{\mathbb{R}}

\newcommand{\NN}{\mathbb{N}}
\newcommand{\ZZ}{\mathbb{Z}}

\let\on=\operatorname

\newenvironment{demo}[1]{\par\smallskip\noindent{\bf #1.}}{\par\smallskip}

\makeatletter
\@namedef{subjclassname@2020}{%
  \textup{2020} Mathematics Subject Classification}
\makeatother

\newsavebox{\fmbox}
\newenvironment{fmpage}[1]
 {\begin{lrbox}{\fmbox}\begin{minipage}{#1}}
 {\end{minipage}\end{lrbox}\fbox{\usebox{\fmbox}}}

\title[On Orlicz classes defined in terms of associated weight functions]
{On Orlicz classes defined in terms of associated weight functions}

\author[G.~Schindl]{Gerhard Schindl}

\address{G.~Schindl: Fakult\"at f\"ur Mathematik, Universit\"at Wien, Oskar-Morgenstern-Platz~1, A-1090 Wien, Austria.}
\email{gerhard.schindl@univie.ac.at}

\begin{document}

\begin{abstract}
N-functions and their growth and regularity properties are crucial in order to introduce and study Orlicz classes and Orlicz spaces. We consider N-functions which are given in terms of so-called associated weight functions. These functions are frequently appearing in the theory of ultradifferentiable function classes and in this setting additional information is available since associated weight functions are defined in terms of a given weight sequence. We express and characterize several known properties for N-functions purely in terms of weight sequences which allows to construct (counter-)examples. Moreover, we study how for abstractly given N-functions this framework becomes meaningful and finally we establish a connection between the complementary N-function and the recently introduced notion of the so-called dual sequence.
\end{abstract}

\thanks{This research was funded in whole by the Austrian Science Fund (FWF) project 10.55776/P33417}
\date{\today}

\maketitle


\section{Introduction}\label{Introduction}
Let us start by recalling briefly the basic definitions of Orlicz classes and Orlicz spaces, we refer to \cite{orliczbook}, \cite{orliczbook1} and to the informative summary presented in \cite{orliczAlexopoulos04}. For this let $F$ be a so-called {\itshape N-function}, see Definition \ref{Nfunctiondef} below and \cite[Chapter 1, \S 1, 3, p. 6]{orliczbook}, \cite[Sect. 2.1, p. 13]{orliczbook1}, \cite[Def. 2.1]{orliczAlexopoulos04}. Moreover, let $\Omega$ be a bounded and closed set in $\RR^d$ and consider on $\Omega$ the usual Lebesgue measure. Then the {\itshape Orlicz class} is given by
$$\mathcal{L}_F(\Omega):=\{u:\Omega\rightarrow\RR,\;\text{measurable}:\;\;\;\int_{\Omega}F(u(x))dx<+\infty\},$$
and the {\itshape Orlicz space} by
$$\mathcal{L}^{*}_F(\Omega):=\{u:\Omega\rightarrow\RR,\;\text{measurable}: \int_{\Omega}u(x)v(x)dx<+\infty,\;\;\;\forall\;v\in\mathcal{L}_{F^c}(\Omega)\},$$
with $F^c$ denoting the so-called {\itshape complementary N-function} which is again an N-function, see Section \ref{complementarysubsection} and \cite[Chapter 1, \S 2, p. 11]{orliczbook}, \cite[Sect. 1.3, p. 6]{orliczbook1}, \cite[Def. 2.2]{orliczAlexopoulos04}. If we do not need to specify the set $\Omega$ we omit it and only write $\mathcal{L}_F$ resp. $\mathcal{L}^{*}_F$. In \cite[Chapter I, Thm. 2, Cor. 3]{orliczbook1}, by revisiting a result by de la Vall\'{e}e Poussin, it has been shown how such a growth restriction expressed in terms of certain convex functions $F$ is arising naturally. In the literature the aforementioned sets are occasionally also defined in an even more general measure theoretic setting and slightly different assumptions on $F$ are used; e.g. using {\itshape Young functions} in \cite{orliczbook1}, see Remark \ref{failurecounterrem} for more details.\vspace{6pt}

In order to study these classes several growth and regularity assumptions for $F$ and $F^c$ are considered frequently in the literature. Most prominent are the so-called $\Delta_2$, $\Delta_3$, $\Delta^2$ and $\Delta'$ condition for $F$, see e.g. \cite[Chapter I, \S 4-\S 6]{orliczbook}, \cite[Chapter II]{orliczbook1}, \cite[Sect. 2.2]{orliczAlexopoulos04} and Section \ref{conditionsection} in this work. If $F^c$ satisfies a ``$\Delta$-type'' property then by convention usually one writes that $F$ has the corresponding ``$\nabla$-type'' condition (and vice versa).\vspace{6pt}

The aim of this paper is to introduce and study N-functions $F_M$ which are given in terms of a given sequence $M\in\RR_{>0}^{\NN}$, see Definition \ref{assoNfctdef}, via the so-called {\itshape associated weight function} $\omega_M$ (see Section \ref{assofctsect}). Here the sequence is expressing the growth of $F_M$ and $M$ is assumed to satisfy mild standard and growth assumptions, see Section \ref{weightsequences}. Recall that $\omega_M$ is appearing frequently in the theory of classes of ultradifferentiable (and ultraholomorphic) functions defined in terms of weight sequences and it serves also as an example for an abstractly given weight function $\omega$ in the sense of Braun-Meise-Taylor, see \cite{BraunMeiseTaylor90}.

Consequently, $F_M$ contains additional information expressed in the underlying sequence $M$ and the idea is to exploit this fact, to ''combine'' the ultradifferentiable-type and the Orlicz-type setting and to treat the following questions/problems:

\begin{itemize}
\item[$(*)$] Study for $F_M$ the aforementioned known and important growth properties for abstractly given N-functions in terms of $M$. When given two N-functions $F_M$, $F_L$ expressed in terms of sequences $M$ and $L$, then study the crucial relation between N-functions (see \eqref{N-functionrelation}) in terms of a growth comparison between $M$ and $L$.

\item[$(*)$] Use this knowledge in order to construct (counter)-examples illustrating the relations and connections between the different growth conditions for N-functions.

\item[$(*)$] Compare the (partially) new growth properties for weight sequences with known conditions appearing in the ultradifferentiable setting.

\item[$(*)$] Check if these properties and conditions can be transferred from given $M$, $L$ to related constructed sequences, e.g. the point-wise product $M\cdot L$ and the convolution product $M\star L$ (see \eqref{convolvesequ}).

\item[$(*)$] Let $G$ be an abstractly given N-function. Is it then possible to associate with $G$ a weight sequence, say $M^G$, and to apply the derived results in order to get information for $G$ itself (via using $F_{M^G}$)?

\item[$(*)$] When given $M$ and $F_M$ study and establish the connection between the notions of the complementary N-function $F_M^c$ and the dual sequence $D$ (w.r.t. $M$) which has been introduced in \cite[Def. 2.1.40, p. 81]{dissertationjimenez}. This question has served as the main motivation for writing this article. The relevance of $D$ is given by the fact that the so-called orders and {\itshape Matuszewska indices} for $M$ and $D$ are ``reflected/inverted'' as it has been shown in the main result \cite[Thm. 2.1.43]{dissertationjimenez}. Concerning these notions we refer to \cite[Sect. 2.1.2]{dissertationjimenez}, \cite{index} and the citations there for more details and precise definitions.
\end{itemize}

However, it turns out that in the weight sequence setting we cannot expect that the relevant function $t\mapsto\varphi_{\omega_M}(t):=\omega_M(e^t)$ (see \eqref{intreprvarphi}) directly is an N-function, see Remark \ref{Nfunctionfailrem} for more explanations. One can overcome this technical problem by using the fact that $\varphi_{\omega_M}$ is the so-called {\itshape principal part} of an N-function $F_M$, see Definition \ref{principaldef} and Corollary \ref{princor}. On the other hand we mention that $\varphi_{\omega_M}$ also allows to compare different used notions for being a weight in the ``Orlicz-setting'', see Remark \ref{failurecounterrem} for more details.\vspace{6pt}

Note that the crucial conditions for $M$ in order to ensure the desired growth properties for $F_M$ are partially (slightly) different compared with the known ones used in the ultradifferentiable setting. This is mainly due to the fact that the relevant function under consideration is given by $\varphi_{\omega_M}$ and not by $\omega_M$ directly. For example, the prominent $\Delta_2$-property for N-functions (see Section \ref{Delta2section}) is also appearing as a known growth condition in the ultradifferentiable weight function setting (abbreviated by \hyperlink{om1}{$(\omega_1)$} in this work) but the crucial condition for $M$ is different (see Theorem \ref{delta2lemma} and the comments there).\vspace{6pt}

The paper is structured as follows: In Section \ref{weightsection} all relevant definitions concerning weight sequences and (associated) weight functions are given and in Section \ref{Nfunctionsect} we recall and introduce the notions of (associated) N-functions. In Section \ref{assoNcomparison} we focus on the study of the comparison between associated N-functions, see Theorems \ref{comparisoncharactthm}, \ref{comparisonprop0}, \ref{comparisonprop} and Corollary \ref{comparisonpropcor}, and give in Section \ref{sufficiencysection} several sufficient conditions on the sequences to ensure equivalence between the associated N-functions.

Section \ref{fromNtoassofct} is dedicated to the study of the meaning of the associated weight sequence $M^G$ when $G$ is an abstractly given N-function, see Theorems \ref{assoweightequivlemma} and \ref{assoseququotientcomparioncor}. In Section \ref{complementarysection} we introduce and study the complementary N-function $F^c_M$ (see Theorem \ref{princorcom1} and Corollary \ref{princorcom1cor}) and establish the connection between $F^c_M$ and the dual sequence $D$, see the main statement Theorem \ref{princorcom2}. Finally, in Section \ref{conditionsection} we provide a detailed study of growth and regularity conditions for N-functions in the weight sequence setting, see Theorems \ref{delta2lemma}, \ref{nabla2charact}, \ref{assoNfctDeltasquare}, \ref{Delta3equivforM} and Proposition \ref{Deltaprimefails}. Some (counter-)examples and their consequences are mentioned as well, see \eqref{counterexamplesequ} and Corollary \ref{nabla2notDelta3cor}.


\section{Weights and conditions}\label{weightsection}

\subsection{Weight sequences}\label{weightsequences}
We write $\NN=\{0,1,2,\dots\}$ and $\NN_{>0}:=\{1,2,\dots\}$ and $\RR_{>0}$ denotes the set $(0,+\infty)$ of positive real numbers.

Let a sequence $M=(M_j)_j$ of positive real numbers be given, i.e. $M\in\RR_{>0}^{\NN}$. We also use the corresponding sequence $\mu=(\mu_j)_j$ defined by $\mu_j:=\frac{M_j}{M_{j-1}}$, $\mu_0:=1$, and analogously for all other appearing sequences. $M$ is called {\itshape normalized} if $1=M_0\le M_1$ holds. For any $\ell>0$ we put $M^{\ell}:=(M^{\ell}_j)_{j\in\NN}$, i.e. the $\ell$-th power, and write $M\cdot L=(M_jL_j)_{j\in\NN}$. Finally, let us introduce the {\itshape convolved sequence} $M\star L$ by
\begin{equation}\label{convolvesequ}
M\star L_j:=\min_{0\le k\le j}M_kL_{j-k},\;\,\;j\in\NN,
\end{equation}
see \cite[$(3.15)$]{Komatsu73}.\vspace{6pt}

$M$ is called {\itshape log-convex} if
$$\forall\;j\in\NN_{>0}:\;M_j^2\le M_{j-1} M_{j+1},$$
equivalently if $(\mu_j)_j$ is non-decreasing. If $M$ is log-convex and normalized, then both $j\mapsto M_j$ and $j\mapsto(M_j)^{1/j}$ are non-decreasing and $(M_j)^{1/j}\le\mu_j$ for all $j\in\NN_{>0}$.\vspace{6pt}

$M$ (with $M_0=1$) has condition {\itshape moderate growth}, denoted by \hypertarget{mg}{$(\text{mg})$}, if
$$\exists\;C\ge 1\;\forall\;j,k\in\NN:\;M_{j+k}\le C^{j+k} M_j M_k.$$
In \cite{Komatsu73} this is denoted by $(M.2)$ and also known under the name {\itshape stability under ultradifferential operators.}\vspace{6pt}

For our purpose it is convenient to consider the following set of sequences
$$\hypertarget{LCset}{\mathcal{LC}}:=\{M\in\RR_{>0}^{\NN}:\;M\;\text{is normalized, log-convex},\;\lim_{j\rightarrow+\infty}(M_j)^{1/j}=+\infty\}.$$
We see that $M\in\hyperlink{LCset}{\mathcal{LC}}$ if and only if $1=\mu_0\le\mu_1\le\dots$, $\lim_{j\rightarrow+\infty}\mu_j=+\infty$ (see e.g. \cite[p. 104]{compositionpaper}) and there is a one-to-one correspondence between $M$ and $\mu=(\mu_j)_j$ by taking $M_j:=\prod_{k=0}^j\mu_k$. If $M,L\in\hyperlink{LCset}{\mathcal{LC}}$, then $M\cdot L, M\star L\in\hyperlink{LCset}{\mathcal{LC}}$ (for the convolution see \cite[Lemma 3.5]{Komatsu73}).\vspace{6pt}

Let $M,L\in\RR_{>0}^{\NN}$ be given, then write $M\le L$ if $M_j\le L_j$ for all $j\in\NN$ and $M\hypertarget{preceq}{\preceq}L$ if $\sup_{j\in\NN_{>0}}\left(\frac{M_j}{L_j}\right)^{1/j}<+\infty$. Sequences $M$ and $L$ are called {\itshape equivalent}, denoted by $M\hypertarget{approx}{\approx}L$, if $M\hyperlink{preceq}{\preceq}L$ and $L\hyperlink{preceq}{\preceq}M$.

\begin{example}
Frequently we will consider the following important examples belonging to the class \hyperlink{LCset}{$\mathcal{LC}$}:
\begin{itemize}
	\item[$(i)$] The Gevrey-sequences $G^s$, $s>0$, given by $G^s_j:=j!^s$.
	
	\item[$(ii)$] The sequences $M^{q,n}$, $q,n>1$, given by $M^{q,n}_j:=q^{j^n}$. If $n=2$, then $M^{q,2}$ is the so-called $q$-Gevrey-sequence.
\end{itemize}
\end{example}

\subsection{Associated weight function}\label{assofctsect}
Let $M\in\RR_{>0}^{\NN}$ (with $M_0=1$), then the {\itshape associated function} $\omega_M: \RR\rightarrow\RR\cup\{+\infty\}$ is defined by
\begin{equation*}\label{assofunc}
\omega_M(t):=\sup_{j\in\NN}\log\left(\frac{|t|^j}{M_j}\right)\;\;\;\text{for}\;t\in\RR,\;t\neq 0,\hspace{30pt}\omega_M(0):=0.
\end{equation*}
For an abstract introduction of the associated function we refer to \cite[Chapitre I]{mandelbrojtbook}, see also \cite[Definition 3.1]{Komatsu73}. Note that $\omega_M$ is here extended to whole $\RR$ in a symmetric (even) way.

If $\liminf_{j\rightarrow+\infty}(M_j)^{1/j}>0$, then $\omega_M(t)=0$ for sufficiently small $t$, since $\log\left(\frac{t^j}{M_j}\right)<0\Leftrightarrow t<(M_j)^{1/j}$ holds for all $j\in\NN_{>0}$. (In particular, if $M_j\ge 1$ for all $j\in\NN$, then $\omega_M$ is vanishing on $[0,1]$.) Moreover, under this assumption $t\mapsto\omega_M(t)$ is a continuous non-decreasing function, which is convex in the variable $\log(t)$ and tends faster to infinity than any $\log(t^j)$, $j\ge 1$, as $t\rightarrow+\infty$. $\lim_{j\rightarrow+\infty}(M_j)^{1/j}=+\infty$ implies that $\omega_M(t)<+\infty$ for each finite $t$ which shall be considered as a basic assumption for defining $\omega_M$.\vspace{6pt}

For given $M\in\hyperlink{LCset}{\mathcal{LC}}$ we define the {\itshape counting function} $\Sigma_M:[0,+\infty)\rightarrow\NN$ by
\begin{equation}\label{counting}
\Sigma_M(t):=|\{j\in\NN_{>0}:\;\;\;\mu_j\le t\}|,
\end{equation}
i.e. $\Sigma_M(t)$ is the maximal positive integer $j$ such that $\mu_j\le t$ (and $\Sigma_M(t)=0$ for $0\le t<\mu_1$). It is known that $\omega_M$ and $\Sigma_M$ are related by the following integral representation formula, see e.g. \cite[1.8. III]{mandelbrojtbook} and \cite[$(3.11)$]{Komatsu73}:
\begin{equation}\label{intrepr}
\omega_M(t)=\int_0^t\frac{\Sigma_M(u)}{u}du=\int_{\mu_1}^t\frac{\Sigma_M(u)}{u}du.
\end{equation}
Consequently, $\omega_M$ vanishes on $[0,\mu_1]$, in particular on the unit interval.

By definition of $\omega_M$ the following formula is immediate:
\begin{equation}\label{powersub}
\forall\;\ell>0\;\forall\;t\ge 0:\;\;\;\ell\omega_M(t^{1/\ell})=\omega_{M^{\ell}}(t).
\end{equation}
In \cite[Lemma 3.5]{Komatsu73} for given $M,L\in\hyperlink{LCset}{\mathcal{LC}}$ it is shown that
$$\forall\;t\ge 0:\;\;\;\Sigma_{M\star L}(t)=\Sigma_M(t)+\Sigma_L(t),$$
which implies by \eqref{intrepr}
$$\forall\;t\ge 0:\;\;\;\omega_{M\star L}(t)=\omega_M(t)+\omega_L(t).$$

Finally, if $M\in\hyperlink{LCset}{\mathcal{LC}}$, then we can compute $M$ by involving $\omega_M$ as follows, see \cite[Chapitre I, 1.4, 1.8]{mandelbrojtbook} and also \cite[Prop. 3.2]{Komatsu73}:
\begin{equation}\label{Prop32Komatsu}
M_j=\sup_{t\ge 0}\frac{t^j}{\exp(\omega_{M}(t))},\;\;\;j\in\NN.
\end{equation}

\begin{remark}\label{strictincrem}
Let $M\in\hyperlink{LCset}{\mathcal{LC}}$ be given, we comment on the surjectivity of $\Sigma_M$.
\begin{itemize}
\item[$(*)$] Obviously $\Sigma_M(t)\in\NN$ for all $t\ge 0$ and $\Sigma_M$ is surjective if and only if $\mu_j<\mu_{j+1}$ for all $j\in\NN_{>0}$, i.e. if $j\mapsto\mu_j$ is strictly increasing: In this case we have $\Sigma_M(t)=j$ for all $\mu_j\le t<\mu_{j+1}$, $j\in\NN_{>0}$, and $\Sigma_M(t)=0$ for $t\in[0,\mu_1)$.

\item[$(*)$] Note that $\mu_j<\mu_{j+1}$ for all $j$ does not hold automatically for all sequences belonging to the set \hyperlink{LCset}{$\mathcal{LC}$}.

However, when given $M\in\hyperlink{LCset}{\mathcal{LC}}$, then we can always find $\widetilde{M}\in\hyperlink{LCset}{\mathcal{LC}}$ such that $M$ and $\widetilde{M}$ are equivalent and such that the corresponding sequence of quotients $(\widetilde{\mu}_j)_{j\ge 1}$ is strictly increasing, see \cite[Lemma 3.18]{optimalflat}. This formal switch allows to avoid technical complications resp. to simplify arguments.

More precisely, in \cite[Lemma 3.18]{optimalflat} it has been shown that even
\begin{equation}\label{congrelation}
0<\inf_{j\in\NN}\frac{\mu_j}{\widetilde{\mu}_j}\le\sup_{j\in\NN}\frac{\mu_j}{\widetilde{\mu}_j}<+\infty,
\end{equation}
which clearly implies $M\hyperlink{approx}{\approx}\widetilde{M}$. We write $M\hypertarget{cong}{\cong}N$ if \eqref{congrelation} holds for the corresponding sequences of quotients $\mu$, $\nu$.
\end{itemize}
\end{remark}

\subsection{Growth properties for abstractly given weight functions}
Let $\omega:[0,+\infty)\rightarrow[0,+\infty)$, we introduce the following growth and regularity conditions
$$\hypertarget{om1}{(\omega_1)}:\;\;\;\omega(2t)=O(\omega(t))\hspace{15pt}t\rightarrow+\infty,$$

$$\hypertarget{om3}{(\omega_3)}:\;\;\;\log(t)=o(\omega(t))\hspace{15pt}t\rightarrow+\infty,$$

$$\hypertarget{om4}{(\omega_4)}:\;\;\;\varphi_{\omega}:t\mapsto\omega(e^t)\;\text{is a convex function on}\;\RR.$$
These conditions are named after \cite{dissertation}. \hyperlink{om1}{$(\omega_1)$}, \hyperlink{om3}{$(\omega_3)$} and \hyperlink{om4}{$(\omega_4)$} are standard assumptions in the theory of ultradifferentiable functions defined by so-called {\itshape Braun-Meise-Taylor weight functions} $\omega$, see \cite{BraunMeiseTaylor90}.

We write that $\omega$ has $\hypertarget{om0}{(\omega_0)}$ if $\omega$ is continuous, non-decreasing, $\omega(t)=0$ for all $t\in[0,1]$ (normalization) and $\lim_{t\rightarrow+\infty}\omega(t)=+\infty$. Finally let us put
$$\hypertarget{omset0}{\mathcal{W}_0}:=\{\omega:[0,\infty)\rightarrow[0,\infty): \omega\;\text{has}\;\hyperlink{om0}{(\omega_0)},\hyperlink{om3}{(\omega_3)},\hyperlink{om4}{(\omega_4)}\}.$$

If $M\in\hyperlink{LCset}{\mathcal{LC}}$ then $\omega_M\in\hyperlink{omset0}{\mathcal{W}_0}$, see e.g. \cite[Lemma 3.1]{sectorialextensions1} and the citations there.

\section{N-functions in the weight sequence setting}\label{Nfunctionsect}
\subsection{Basic definitions and abstractly given N-functions}\label{Nfctdefsect}
We revisit the basic definitions from \cite[Chapter I, \S 1, 3, p. 6]{orliczbook} and \cite[Def. 2.1]{orliczAlexopoulos04}. Consider $f:[0,+\infty)\rightarrow[0,+\infty)$ with the following properties:
\begin{itemize}
\item[$(I)$] $f$ is right-continuous and non-decreasing;

\item[$(II)$] $f(0)=0$ and $f(t)>0$ for all $t>0$;

\item[$(III)$] $\lim_{t\rightarrow+\infty}f(t)=+\infty$.
\end{itemize}

Then we give the following definition, see e.g. \cite[Chapter I, \S 1, 3, p. 6]{orliczbook}.

\begin{definition}\label{Nfunctiondef}
Let $f:[0,+\infty)\rightarrow[0,+\infty)$ be having $(I)$, $(II)$ and $(III)$. The function $F:\RR\rightarrow[0,+\infty)$ defined by
\begin{equation}\label{Ndef}
F(x):=\int_0^{|x|}f(t)dt,
\end{equation}
is called an {\itshape N-function.}
\end{definition}

\begin{remark}\label{RaoRendefrem}
In \cite[Sect. 1.3, Thm. 1, Cor. 2; Sect. 2.1]{orliczbook1} an analogous integral representation for N- resp. even Young functions has been obtained with the integrand $f$ (``density'') being non-decreasing and {\itshape left-continuous}. Thus the authors are working with the left-derivative of $F$. Since we are focusing on the weight sequence case, see \eqref{intreprvarphi} in Section \ref{fromassotoN}, we have to involve the counting function $\Sigma_M$ from \eqref{counting} which is right-continuous by the very definition and so we prefer to work within the above setting.
\end{remark}

Every N-function $F$ satisfies the following properties, see \cite[Chapter I, \S 1, 4, p. 7]{orliczbook}:
\begin{itemize}
\item[$(*)$] $F(0)=0$ (normalization) and $F(x)>0$ for all $x\neq 0$,

\item[$(*)$] $F$ is even, non-decreasing, continuous, and convex.

\item[$(*)$] The convexity and $F(0)=0$ imply that
\begin{equation}\label{equ114}
\forall\;0\le t\le 1\;\forall\;u\ge 0:\;\;\;F(tu)\le tF(u),
\end{equation}
see \cite[$(1.14)$]{orliczbook}. This holds since by convexity we have $F(tx+(1-t)y)\le tF(x)+(1-t)F(y)$ for all $0\le t\le 1$ and $x,y\ge 0$ and then set $y=0$.

\item[$(*)$] Finally, let us recall
\begin{equation}\label{equ115116}
\lim_{t\rightarrow 0}\frac{F(t)}{t}=0,\hspace{20pt}\lim_{t\rightarrow+\infty}\frac{F(t)}{t}=+\infty,
\end{equation}
see \cite[$(1.15), (1.16)$]{orliczbook}, and which follows from $(II)$ resp. $(III)$ for $f$.
\end{itemize}

When given two N-functions (or even arbitrary functions) $F_1,F_2: [0,+\infty)\rightarrow[0,+\infty)$, then write $F_1\hypertarget{ompreceq}{\preceq_{\mathfrak{c}}}F_2$ if
\begin{equation}\label{N-functionrelation}
\exists\;K>0\;\exists\;t_0>0\;\forall\;t\ge t_0:\;\;\;F_1(t)\le F_2(Kt).
\end{equation}
If either $F_1$ or $F_2$ is non-decreasing then w.l.o.g. we can restrict to $K\in\NN_{>0}$ and this relation is clearly reflexive and transitive. It has been introduced in \cite[Chapter I, \S 3]{orliczbook} and \cite[Sect. 2.2, Def. 1]{orliczbook1}. In \cite{orliczbook}, $F_1$ and $F_2$ are called {\itshape comparable,} if either $F_1\hyperlink{ompreceq}{\preceq_{\mathfrak{c}}}F_2$ or $F_2\hyperlink{ompreceq}{\preceq_{\mathfrak{c}}}F_1$. In \cite{orliczbook1}, when $F_1$, $F_2$ are related by \eqref{N-functionrelation}, then $F_2$ has been called {\itshape stronger} than $F_1$.

For this relation we are gathering several equivalent reformulations.

\begin{lemma}\label{relationrem}
Let $F_1,F_2:[0,+\infty)\rightarrow[0,+\infty)$ be non-decreasing. Assume that either $F_1$ or $F_2$ is normalized, convex and tending to infinity as $t\rightarrow+\infty$. Then the following are equivalent:

\begin{itemize}
\item[$(i)$] $F_1\hyperlink{ompreceq}{\preceq_{\mathfrak{c}}}F_2$ holds.

\item[$(ii)$] We have that
\begin{equation}\label{relationconstoutside}
\exists\;K,K_1\ge 1\;\exists\;t_0>0\;\forall\;t\ge t_0:\;\;\;F_1(t)\le K_1F_2(Kt).
\end{equation}

\item[$(iii)$] We have that
\begin{equation}\label{implwithC}
\exists\;C>0\;\exists\;K\ge 1\;\forall\;t\ge 0:\;\;\;F_1(t)\le F_2(Kt)+C.
\end{equation}

\item[$(iv)$] We have that
\begin{equation}\label{implwithCextra}
\exists\;K,K_1\ge 1\;\exists\;C>0\;\forall\;t\ge 0:\;\;\;F_1(t)\le K_1F_2(Kt)+C.
\end{equation}
\end{itemize}
\end{lemma}
In particular, the above characterization applies if both $F_1$ and $F_2$ are N-functions.

\demo{Proof}
$(i)\Rightarrow(ii)$ is trivial and $(ii)\Rightarrow(i)$ follows by \eqref{equ114}: If $F_2$ is normalized and convex, then when given $K_1>1$ we put $t:=K_1^{-1}$ in \eqref{equ114} and hence $K_1F_2(u)\le F_2(K_1u)$ for all $u\ge 0$. Similarly, if $F_1$ is normalized and convex, then the assumption gives $K_1^{-1}F_1(t)\le F_2(Kt)$ and so $F_1(t)\le K_1^{-1}F_1(K_1t)\le F_2(KK_1t)$ for all $t\ge t_0$ holds which shows $F_1\hyperlink{ompreceq}{\preceq_{\mathfrak{c}}}F_2$ when choosing $K_2:=KK_1$.

$(i)\Rightarrow(iii)$ is clear, since $F_2(t)\ge 0$ and $F_1$ is non-decreasing take e.g. $C:=F_1(t_0)$.

$(iii)\Rightarrow(ii)$ When given $C\ge 1$ then we have $F_2(Kt)+C\le 2F_2(Kt)$ for all sufficiently large $t$ if $\lim_{t\rightarrow+\infty}F_2(t)=+\infty$. Thus \eqref{relationconstoutside} is verified with $K_1:=2$ and the same $K$. If $\lim_{t\rightarrow+\infty}F_1(t)=+\infty$, then by \eqref{implwithC} also $\lim_{t\rightarrow+\infty}F_2(t)=+\infty$ and the rest follows as before.

$(iii)\Rightarrow(iv)$ is trivial and $(iv)\Rightarrow(iii)$ holds as $(ii)\Rightarrow(i)$.
\qed\enddemo

This motivates the following definition, see \cite[Chapter I, \S 3]{orliczbook}.

\begin{definition}
We call two functions $F_1$ and $F_2$ {\itshape equivalent,} written $F_1\hypertarget{sim}{\sim_{\mathfrak{c}}}F_2$, if $F_1\hyperlink{ompreceq}{\preceq_{\mathfrak{c}}}F_2$ and $F_2\hyperlink{ompreceq}{\preceq_{\mathfrak{c}}}F_1$.
\end{definition}

In particular, for any N-function $F$ we have that all $F_k:t\mapsto F(kt)$, $k>0$, are equivalent.

In \cite[Thm. 13.2]{orliczbook} it has been shown that $F_1\hyperlink{sim}{\sim_{\mathfrak{c}}}F_2$ if and only if $\mathcal{L}^{*}_{F_1}=\mathcal{L}^{*}_{F_2}$.

\begin{remark}\label{convexityrem}
We comment on relation \hyperlink{sim}{$\sim_{\mathfrak{c}}$} for given N-functions $F_1, F_2$ and their corresponding right-derivatives $f_1,f_2$ appearing in \eqref{Ndef}:
\begin{itemize}
\item[$(i)$] On \cite[p. 15]{orliczbook} it is mentioned that if
\begin{equation}\label{equ34}
\exists\;b\in(0,+\infty):\;\;\;\lim_{t\rightarrow+\infty}\frac{F_1(t)}{F_2(t)}=b,
\end{equation}
then $F_1\hyperlink{sim}{\sim_{\mathfrak{c}}}F_2$ is valid. Indeed, this implication holds for any non-decreasing functions $F_1,F_2:[0,+\infty)\rightarrow[0,+\infty)$ such that either $F_1$ or $F_2$ is assumed to be convex and normalized:

For any $0<a\le 1$ we clearly have $aF_2(u)\le F_2(u)$ and if $a>1$, then as in the proof of Lemma \ref{relationrem} the estimate \eqref{equ114} applied to $t:=a^{-1}$ gives $aF_2(u)\le F_2(au)$ for all $u\ge 0$. The proof for $F_1$ is analogous.

In particular, \eqref{equ34} holds (with $b=1$) if $F_1(t)=F_2(t)$ for all $t$ large.

\item[$(ii)$] Moreover, if $\lim_{t\rightarrow+\infty}F_i(t)=+\infty$, $i=1,2$, then \eqref{equ34} holds with $b=1$ if
\begin{equation}\label{equivwithC}
\exists\;C,D\ge 1\;\forall\;t\ge 0:\;\;\;F_1(t)-C\le F_2(t)\le F_1(t)+D.
\end{equation}

\item[$(iii)$] In \cite[Lemma 3.1]{orliczbook} it has been shown that $f_1\hyperlink{ompreceq}{\preceq_{\mathfrak{c}}}f_2$ implies $F_1\hyperlink{ompreceq}{\preceq_{\mathfrak{c}}}F_2$ and in \cite[Lemma 3.2]{orliczbook} that
\begin{equation}\label{equ38}
\exists\;b\in(0,+\infty):\;\;\;\lim_{t\rightarrow+\infty}\frac{f_1(t)}{f_2(t)}=b
\end{equation}
implies $F_1\hyperlink{sim}{\sim_{\mathfrak{c}}}F_2$.

\item[$(iv)$] Finally, let us characterize this relation for N-functions $F_1$, $F_2$ in terms of the right-derivatives $f_1$, $f_2$; see also \cite[Sect. 2.2, Thm. 2]{orliczbook1} and \cite[Lemma 3.1, p. 17-18]{orliczbook}: We have $F_1\hyperlink{ompreceq}{\preceq_{\mathfrak{c}}}F_2$ if and only if
    \begin{equation}\label{rightderrelation}
    \exists\;k>1\;\exists\;t_1>0\;\forall\;t\ge t_1:\;\;\;f_1(t)\le kf_2(kt).
    \end{equation}
    Note: $f_1\hyperlink{ompreceq}{\preceq_{\mathfrak{c}}}f_2$ implies the above relation and since $f_2$ is non-decreasing, equivalently we can use in \eqref{rightderrelation} the control $k_1f_2(k_2t)$ for some $k_1,k_2>1$.

    On the one hand, for all $x\ge t_1$ with $t_1$ denoting the value in \eqref{rightderrelation}
    \begin{align*}
    F_1(x)&=\int_0^xf_1(s)ds=\int_0^{t_1}f_1(s)ds+\int_{t_1}^xf_1(s)ds\le F_1(t_1)+\int_{t_1}^xkf_2(ks)ds\le F_1(t_1)+k\int_0^xf_2(ks)ds
    \\&
    =F_1(t_1)+k\int_0^{kx}f_2(u)\frac{du}{k}=F_1(t_1)+F_2(kx),
    \end{align*}
    hence \eqref{implwithC} is verified with $C:=2F_1(t_1)$ and $K:=k$. Conversely, let $x>0$ with $x\ge t_0/2$, $t_0$ from \eqref{N-functionrelation}, and estimate as follows:
    \begin{align*}
    xf_1(x)&\le\int_x^{2x}f_1(s)ds\le\int_0^{2x}f_1(s)ds=F_1(2x)\le F_2(2Kx)=\int_0^{2Kx}f_2(s)ds\le 2Kxf_2(2Kx).
    \end{align*}
    So \eqref{rightderrelation} is verified with $k:=2K$ and $t_1:=t_0/2$.
\end{itemize}
\end{remark}

\begin{remark}\label{preceqrem}
In the theory of Orlicz classes also the following slightly different relation between functions is considered:

\begin{itemize}
\item[$(*)$] Let $F_1,F_2$ be N-functions. In \cite[Thm. 8.1]{orliczbook}, see also \cite[Thm. 3.3]{orliczAlexopoulos04}, it has been shown that $\mathcal{L}_{F_1}(\Omega)\subseteq\mathcal{L}_{F_2}(\Omega)$ (as sets) if and only if
\begin{equation}\label{bigOrelation}
\exists\;t_0>0\;\exists\;K\ge 1\;\forall\;t\ge t_0:\;\;\;F_2(t)\le KF_1(t),
\end{equation}
i.e. $F_2(t)=O(F_1(t))$ as $t\rightarrow+\infty$. The sufficiency of \eqref{bigOrelation} for having the inclusion $\mathcal{L}_{F_1}(\Omega)\subseteq\mathcal{L}_{F_2}(\Omega)$ is clear.

\item[$(*)$] Given $F_1,F_2:[0,+\infty)\rightarrow[0,+\infty)$ we write $F_1\preceq F_2$ if \eqref{bigOrelation} holds and $F_1\sim F_2$ if $F_1\preceq F_2$ and $F_2\preceq F_1$.

Note that this relation $\sim$ is precisely \cite[$(8.6)$]{orliczbook} and it is also the crucial one for the characterization of inclusions (resp. equalities) of classes in the ultradifferentiable weight function setting, see \cite[Sect. 5]{compositionpaper}.

Moreover, let us write $F_1\vartriangleleft F_2$ if $F_2(t)=o(F_1(t))$ as $t\rightarrow+\infty$.

\item[$(*)$] In view of \eqref{equ114}, for given N-functions $F_1$, $F_2$ we have that $F_1\preceq F_2$ implies $F_2\hyperlink{ompreceq}{\preceq_{\mathfrak{c}}}F_1$. For this implication one only requires that either $F_1$ or $F_2$ is normalized and convex, see the proof of $(ii)\Rightarrow(i)$ in Lemma \ref{relationrem}. Consequently, if N-functions $F_1$ and $F_2$ are related by $F_1\sim F_2$, then they are also equivalent.
\end{itemize}
\end{remark}

For the sake of completeness let us summarize more relations between (N-)functions mentioned in \cite[Sect. 2.2, Def. 1]{orliczbook1}:

\begin{itemize}
\item[$(*)$] $F_2$ is {\itshape essentially stronger} than $F_1$, if
$$\forall\;K>0\;\exists\;t_0>0\;\forall\;t\ge t_0:\;\;\;F_1(t)\le F_2(Kt),$$
and write $F_1\preceq\preceq_{\mathfrak{c}}F_2$ for this relation which is, of course, stronger than $F_1\hyperlink{ompreceq}{\preceq_{\mathfrak{c}}}F_2$. When taking $K:=1$ it also implies $F_2\preceq F_1$. Lemma \ref{relationrem} transfers to this relation and using this analogously to $(iv)$ in Remark \ref{convexityrem} one can prove that $F_1\preceq\preceq_{\mathfrak{c}}F_2$ if and only if
\begin{equation}\label{rightderrelationstrong}
\forall\;k>1\;\exists\;t_1>0\;\forall\;t\ge t_1:\;\;\;f_1(t)\le kf_2(kt).
\end{equation}

\item[$(*)$] $F_2$ is {\itshape completely stronger} than $F_1$, if
$$\forall\;\epsilon>0\;\;\exists\;K>0\;\exists\;t_0>0\;\forall\;t\ge t_0:\;\;\;F_1(t)\le K F_2(\epsilon t).$$
The choice $\epsilon:=1$ implies again $F_2\preceq F_1$.

\item[$(*)$] $F_2$ is {\itshape increasing more rapidly} than $F_1$, if
$$\forall\;\epsilon>0\;\exists\;\delta>0\;\exists\;t_0>0\;\forall\;t\ge t_0:\;\;\;\frac{1}{\delta}F_1(\delta t)\le\epsilon F_2(t).$$
If $F_2\vartriangleleft F_1$ is valid, then $F_2$ is increasing more rapidly than $F_1$ by the uniform choice $\delta:=1$ for all $\epsilon>0$.
\end{itemize}

\subsection{From weight sequences to associated N-functions}\label{fromassotoN}
Let $M\in\hyperlink{LCset}{\mathcal{LC}}$ be given. When rewriting \eqref{intrepr} we obtain for all $t\ge 0$ (note that $\mu_1\ge 1$):
\begin{equation}\label{intreprvarphi}
\varphi_{\omega_M}(t)=\omega_M(e^t)=\int_{\mu_1}^{e^t}\frac{\Sigma_M(s)}{s}ds=\int_{\log(\mu_1)}^t\frac{\Sigma_M(e^u)}{s}sdu=\int_{\log(\mu_1)}^t\Sigma_M(e^u)du=\int_0^t\Sigma_M(e^u)du,
\end{equation}
since $\Sigma_M(e^u)=0$ for $0\le u<\log(\mu_1)$. This formula should be compared with \cite[Thm. 1.1, $(1.10)$]{orliczbook}. $\Sigma_M\circ\exp$ is right-continuous, non-decreasing and clearly $\lim_{t\rightarrow+\infty}\Sigma_M(e^t)=+\infty$.

Recall that $\omega_M\in\hyperlink{omset0}{\mathcal{W}_0}$ and so $\varphi_{\omega_M}$ is convex, $t\mapsto\frac{\varphi_{\omega_M}(t)}{t}$ is non-decreasing with $\lim_{t\rightarrow+\infty}\frac{\varphi_{\omega_M}(t)}{t}=+\infty$ and finally $\varphi_{\omega_M}(0)=0$ is valid, see e.g. \cite[Rem. 1.3, Lemma 1.5]{BraunMeiseTaylor90} and also \cite[p. 7]{orliczbook}.

\begin{remark}\label{Nfunctionfailrem}
However, requirement $(II)$ cannot be achieved for $\Sigma_M\circ\exp$ for any $M\in\hyperlink{LCset}{\mathcal{LC}}$: If $M_1=M_0(=1)$, and so $\mu_1=1$, then \eqref{counting} yields $\Sigma_M(e^0)=\Sigma_M(\mu_1)\ge 1\neq 0$. If $M_1>M_0\Leftrightarrow\mu_1>1$, then $\Sigma_M(e^t)=0$ for all $0\le t<\log(\mu_1)$.

Finally remark that, if $M$ is log-convex with $\lim_{j\rightarrow+\infty}(M_j)^{1/j}=+\infty$ but such that normalization fails, then $0<\mu_1<1$ and so $\Sigma_M(e^t)\ge 1$ for any $t\ge 0$. Thus also in this case the first requirement in $(II)$ is violated.

This failure is related to the fact that the first property in \eqref{equ115116} and $\frac{\varphi_{\omega_M}(t)}{t}>0$ for all $t>0$ are not satisfied automatically for $\varphi_{\omega_M}$, see the proofs and arguments in \cite[Chapter I, \S 1, 5, p. 8-9]{orliczbook}. Thus $\varphi_{\omega_M}$ is formally not an N-function according to Definition \ref{Nfunctiondef}.
\end{remark}

In order to overcome this technical problem we recall the following notion, see \cite[Chapter I, \S 3, 3, p. 16]{orliczbook}:

\begin{definition}\label{principaldef}
A convex function $Q$ is called the principal part of an N-function $F$ if
$$\exists\;t_0>0\;\forall\;t\ge t_0:\;\;\;Q(t)=F(t).$$
\end{definition}

We have the following result, see \cite[Thm. 3.3]{orliczbook} and the proof there:

\begin{theorem}\label{princthm}
Let $Q:[0,+\infty)\rightarrow[0,+\infty)$ be a convex function such that $\lim_{t\rightarrow+\infty}\frac{Q(t)}{t}=+\infty$. Then there exists an N-function $F$ such that $Q$ is the principal part of $F$.

More precisely, we even get that
\begin{equation}\label{princthmequ}
\exists\;t_0>0\;\forall\;t\ge t_0:\;\;\;f(t)=q(t),
\end{equation}
with $f$ denoting the function appearing in \eqref{Ndef} of $F$ and $q$ denoting the non-decreasing and right-continuous function appearing in the representation
\begin{equation}\label{Qrepr}
Q(t)=\int_a^tq(s)ds,
\end{equation}
see \cite[$(1.10)$]{orliczbook}. Here $a\ge 0$ is such that $Q(a)=0$ and we have $t_0>a$.
\end{theorem}

\begin{proposition}\label{princprop}
Let $Q$ be a convex function such that $\lim_{t\rightarrow+\infty}\frac{Q(t)}{t}=+\infty$ and let $F$ be the N-function according to Theorem \ref{princthm}. Then we get
\begin{equation}\label{princpropequ}
\exists\;C,D\ge 1\;\forall\;t\ge 0:\;\;\;Q(t)-C\le F(t)\le Q(t)+D,
\end{equation}
cf. \eqref{equivwithC}. This relation implies $\lim_{t\rightarrow+\infty}\frac{F(t)}{Q(t)}=1$ and so both $F\hyperlink{sim}{\sim_{\mathfrak{c}}}Q$ and $F\sim Q$ holds.
\end{proposition}

\demo{Proof}
More generally, when for given functions $F_1, F_2:[0,+\infty)\rightarrow[0,+\infty)$ we get $F_1(t)=F_2(t)$ for all $t$ large then we have $F_1(t)\le F_2(t)+C$, $F_2(t)\le F_1(t)+D$ for all $t\ge 0$ with $C:=\max\{F_1(t): 0\le t\le t_0\}$ and $D:=\max\{F_2(t): 0\le t\le t_0\}$.

Hence \eqref{princpropequ} follows by Theorem \ref{princthm} (recall Definition \ref{principaldef}).

Note that $Q$ is convex but normalization for $Q$ (i.e. $a=0$ in \eqref{Qrepr}) is not guaranteed in general.
\qed\enddemo

\begin{remark}\label{convexdivrem}
Conversely, each convex function $Q:[0,+\infty)\rightarrow[0,+\infty)$ admitting the representation \eqref{Qrepr} for a non-decreasing and right-continuous function $q$ with $q(t)\rightarrow+\infty$ satisfies also $\lim_{t\rightarrow+\infty}\frac{Q(t)}{t}=+\infty$. This holds since for all $t\ge a$ (see the proof of $(1.16)$ in \cite[Chapter I, \S 1, p. 7]{orliczbook}):
$$Q(2t)=\int_a^{2t}q(s)ds\ge\int_t^{2t}q(s)ds\ge q(t)t.$$
\end{remark}

In particular, when applying these results to $Q=\varphi_{\omega_M}$ we get the following consequence:

\begin{corollary}\label{princor}
Let $M\in\hyperlink{LCset}{\mathcal{LC}}$ be given. Then there exists an N-function $F_M$ such that $\varphi_{\omega_M}$ is the principal part of $F_M$ and so
\begin{equation}\label{equivwithCforM}
\exists\;C,D\ge 1\;\forall\;t\ge 0:\;\;\;\varphi_{\omega_M}(t)-C\le F_M(t)\le\varphi_{\omega_M}(t)+D.
\end{equation}
This implies $\varphi_{\omega_M}\sim F_M$ and hence also $\varphi_{\omega_M}\hyperlink{sim}{\sim_{\mathfrak{c}}}F_M$.

Moreover, if $f_M$ denotes the function appearing in the representation \eqref{Ndef} of $F_M$, then we even get
\begin{equation}\label{princorequ}
\exists\;t_0>0\;\forall\;t\ge t_0:\;\;\;f_M(t)=\Sigma_M(e^t).
\end{equation}
In view of this equality we call $\Sigma_M\circ\exp$ the principal part of $f_M$ (see \cite[p. 18]{orliczbook}).
\end{corollary}

\demo{Proof}
We can apply Theorem \ref{princthm} to $Q\equiv\varphi_{\omega_M}$ because $\lim_{t\rightarrow+\infty}\frac{\varphi_{\omega_M}(t)}{t}=\lim_{s\rightarrow+\infty}\frac{\omega_M(s)}{\log(s)}=+\infty$ by \hyperlink{om3}{$(\omega_3)$} (recall \cite[Lemma 3.1]{sectorialextensions1} and the citations there). In fact \hyperlink{om3}{$(\omega_3)$} for $\omega_M$ is precisely \cite[(3.6)]{orliczbook} for $Q=\varphi_{\omega_M}$.

\eqref{equivwithCforM} follows by Proposition \ref{princprop}, and \eqref{princorequ} holds by taking into account the representation \eqref{intreprvarphi}.

Note that by normalization of $M$ we get $\varphi_{\omega_M}(0)=\omega_M(1)=0$ and so in \eqref{Qrepr} we have $a=0$ and $q=\Sigma_M\circ\exp$.
\qed\enddemo

The following is an immediate consequence of \eqref{equivwithCforM}:

\begin{corollary}\label{Bigorelationcor}
Let $M,L\in\hyperlink{LCset}{\mathcal{LC}}$ be given. Then $F_M\hyperlink{ompreceq}{\preceq_{\mathfrak{c}}}F_L$ if and only if $\varphi_{\omega_M}\hyperlink{ompreceq}{\preceq_{\mathfrak{c}}}\varphi_{\omega_L}$ and $F_M\preceq F_L$ if and only if $\varphi_{\omega_M}\preceq\varphi_{\omega_L}$.

Moreover, also all further relations from \cite[Sect. 2.2, Def. 1]{orliczbook1} hold between $F_M$, $F_L$  if and only if they are valid between $\varphi_{\omega_M}$, $\varphi_{\omega_L}$.
\end{corollary}

\begin{definition}\label{assoNfctdef}
Let $M\in\hyperlink{LCset}{\mathcal{LC}}$ be given. Then the N-function $F_M$ from Corollary \ref{princor} is called the associated N-function.
\end{definition}

We close this section by commenting on the relation between $\varphi_{\omega_M}$ and other notions of defining functions in the Orlicz setting.

\begin{remark}\label{failurecounterrem}
As seen above, for any given $M\in\hyperlink{LCset}{\mathcal{LC}}$ we cannot expect that $\varphi_{\omega_M}$ is formally an N-function. On the other hand $\varphi_{\omega_M}$ can be used to illustrate the differences between appearing definitions for Orlicz classes in the literature. In \cite{Osan15} an exhaustive study is provided and the different notions and conditions for the defining functions are compared, see also the literature citations there.

\begin{itemize}
\item[$(*)$] $\varphi_{\omega_M}$ coincides with an N-function (with $F_M$) for sufficiently large values.

\item[$(*)$] For any $M\in\hyperlink{LCset}{\mathcal{LC}}$ the function $\varphi_{\omega_M}$ is always a {\itshape Young function}, see \cite[Def. 1.4]{Osan15} and \cite[Sect. 1.3, p. 6]{orliczbook1}: $\varphi_{\omega_M}$ is convex, satisfies $\varphi_{\omega_M}(0)=\omega_M(1)=0$ by normalization and $\varphi_{\omega_M}(t)\rightarrow+\infty$ as $t\rightarrow+\infty$ (and it can be extend to $\RR$ in an even way). Note that formally for Young functions $F$ it is allowed that $F(a)=+\infty$ for some $a\in\RR$.

\item[$(*)$] $\varphi_{\omega_M}$ is a {\itshape strong Young function} (see \cite[Def. 1.7]{Osan15}) if and only if $\mu_1=1$: Continuity is clear and $\varphi_{\omega_M}(t)>0$ for all $t>0$ follows if and only if $\Sigma_M(e^t)>0$ for all $t\ge 0$, see \eqref{intreprvarphi}. This is clearly equivalent to $\mu_1=1$.

\item[$(*)$] Finally, $\varphi_{\omega_M}$ is always an {\itshape Orlicz function} (see \cite[Def. 1.9]{Osan15}), since $\varphi_{\omega_M}$ is never identically zero or infinity (follows again by \eqref{intreprvarphi}).

\item[$(*)$] Consequently, $\varphi_{\omega_M}$ provides (counter-)examples for the first two strict implications in \cite[Cor. 2.7]{Osan15}, see \cite[Cor. 2.8]{Osan15}: Each N-function is a strong Young function and each strong Young function is an Orlicz function but each implication cannot be reversed in general; one can take $M\in\hyperlink{LCset}{\mathcal{LC}}$ with $\mu_1=1$ for the first and $M\in\hyperlink{LCset}{\mathcal{LC}}$ with $\mu_1>1$ for the second part.
\end{itemize}
\end{remark}

\section{Comparison between associated N-functions}\label{assoNcomparison}
The goal of this Section is to give a connection resp. comparison between the growth relation $\preceq$ for weight sequences, crucially appearing in the theory of ultradifferentiable and ultraholomorphic functions, and the previously defined relations \hyperlink{ompreceq}{$\preceq_{\mathfrak{c}}$} and $\preceq$ for (associated) N-functions.

\subsection{Main statements}
The first main result establishes a characterization of relation $F_M\hyperlink{ompreceq}{\preceq_{\mathfrak{c}}}F_L$ in terms of a growth comparison between $M$ and $L$. However, the characterization is not given via $\preceq$ but expressed in terms of the corresponding sequences of quotients $\mu=(\mu_j)_{j\in\NN}$ and $\lambda=(\lambda_j)_{j\in\NN}$. In explicit applications and for constructing weight sequences $M$ it is often convenient to start with $\mu$. Note that when involving $\mu$ we get automatically growth conditions for the counting function $\Sigma_M$ as well.

\begin{theorem}\label{comparisoncharactthm}
Let $M,L\in\hyperlink{LCset}{\mathcal{LC}}$ be given. Then the following are equivalent:
\begin{itemize}
\item[$(a)$] We have that
\begin{equation}\label{comparisonlemmaequ}
\exists\;A\ge 1\;\exists\;k>0\;\exists\;t_0>0\;\forall\;t\ge t_0:\;\;\;\Sigma_M(t)\le A\Sigma_L(t^k).
\end{equation}
\item[$(b)$] We have that
\begin{equation}\label{comparisonlemmaequ1}
\exists\;B\ge 1\;\exists\;k>0\;\exists\;j_0\in\NN_{>0}\;\forall\;j\ge j_0:\;\;\;\lambda_{\lceil j/B\rceil}\le\mu_j^k.
\end{equation}

\item[$(c)$] We have that $F_M\hyperlink{ompreceq}{\preceq_{\mathfrak{c}}}F_L$ (equivalently $\varphi_{\omega_M}\hyperlink{ompreceq}{\preceq_{\mathfrak{c}}}\varphi_{\omega_L}$) holds.

\item[$(d)$] We have that \eqref{rightderrelation} holds between $f_M$ and $f_L$, i.e.
$$\exists\;k>1\;\exists\;t_1>0\;\forall\;t\ge t_1:\;\;\;f_M(t)\le kf_L(kt).$$
\end{itemize}
The proof shows that we can take $A=B=k$ in $(a)$ and $(b)$ and the result becomes trivial for $M=L$ (set $A=B=k=1$).
\end{theorem}

The analogous characterization for $F_M\preceq\preceq_{\mathfrak{c}}F_L$ is obtained when replacing \eqref{comparisonlemmaequ} by \eqref{rightderrelationstrong} and so \eqref{comparisonlemmaequ1} by the condition
$$\forall\;k>0\;\exists\;j_0\in\NN_{>0}\;\forall\;j\ge j_0:\;\;\;\lambda_{\lceil j/k\rceil}\le\mu_j^k.$$

\demo{Proof}
$(a)\Rightarrow(b)$ Let $j_0\in\NN_{>0}$ be minimal to ensure $\mu_{j_0}\ge t_0$ (note that $\lim_{j\rightarrow+\infty}\mu_j=+\infty$). For any $t\ge\mu_{j_0}$ we have $\mu_j\le t<\mu_{j+1}$ for some $j\in\NN_{>0}$, $j\ge j_0$. Then by assumption $j=\Sigma_M(t)\le A\Sigma_L(t^k)\Leftrightarrow\Sigma_L(t^k)\ge\frac{j}{A}$ and so $t^k\ge\lambda_{\lceil j/A\rceil}$ (note that $\Sigma_L(t^k)\in\NN$). When taking $t:=\mu_j$ we get \eqref{comparisonlemmaequ1} with $B:=A$ and the same $k>0$ for all $j\ge j_0$ with $\mu_j<\mu_{j+1}$.

If $j\ge j_0$ with $\mu_j=\mu_{j+1}$, then $\mu_j=\mu_{j+1}=\dots=\mu_{j+d}<\mu_{j+d+1}$ for some $d\in\NN_{>0}$ and thus $j+d=\Sigma_M(t)$ for $\mu_{j+d}\le t<\mu_{j+d+1}$. This yields $\Sigma_L(t^k)\ge\frac{j+d}{A}$, i.e. $t^k\ge\lambda_{\lceil (j+d)/A\rceil}\ge\lambda_{\lceil (j+i)/A\rceil}$ for all $0\le i\le d$. Put $t:=\mu_{j+d}(=\dots=\mu_j)$ in order to verify \eqref{comparisonlemmaequ1} with $B:=A$ and the same $k>0$ for all $j\ge j_0$.\vspace{6pt}

$(b)\Rightarrow(a)$ Let $t\ge 0$ be such that $\mu_j\le t<\mu_{j+1}$ for some $j\ge j_0$. Then $\mu_j^k\le t^k<\mu^k_{j+1}$ and so $t^k\ge\mu_j^k\ge\lambda_{\lceil j/B\rceil}$ which gives $\Sigma_M(t)=j$ and $\lceil j/B\rceil\le\Sigma_L(t^k)$. Thus \eqref{comparisonlemmaequ} follows when $j\le A j/B\le A\lceil j/B\rceil$ is ensured. So we can take $A:=B$, the same $k>0$ and $t_0:=\mu_{j_0}$.\vspace{6pt}

$(a)\Leftrightarrow(c)$ First, \eqref{comparisonlemmaequ} precisely means that $\Sigma_M(e^s)\le A\Sigma_L(e^{sk})$ for some $A\ge 1$, $k>0$ and all sufficiently large $s$. Thus the desired equivalence holds by taking into account the representation \eqref{intreprvarphi}, following the estimates in $(iv)$ in Remark \ref{convexityrem} with $f_1$ being replaced by $\Sigma_M\circ\exp$ and $f_2$ by $\Sigma_L\circ\exp$ (\eqref{comparisonlemmaequ} is precisely \eqref{rightderrelation} for these choices), and finally using Corollary \ref{princor}.\vspace{6pt}

$(a)\Leftrightarrow(d)$ This holds by \eqref{princorequ}.
\qed\enddemo

We continue with the following result providing a complete characterization for the relation $\preceq$ between associated N-functions.

\begin{theorem}\label{comparisonprop0}
Let $M,L\in\hyperlink{LCset}{\mathcal{LC}}$ be given. Then the following are equivalent:
\begin{itemize}
\item[$(i)$] The associated N-functions satisfy $F_M\preceq F_L$.

\item[$(ii)$] The functions $\varphi_{\omega_M}$ and $\varphi_{\omega_L}$ satisfy $\varphi_{\omega_M}\preceq\varphi_{\omega_L}$.

\item[$(iii)$] The sequences $M$ and $L$ are related by
\begin{equation}\label{comparisonprop0equ}
\exists\;A\ge 1\;\exists\;c\in\NN_{>0}\;\forall\;j\in\NN:\;\;\;M_j\le A(L_{cj})^{1/c}.
\end{equation}
\end{itemize}
\end{theorem}

\demo{Proof}
$(i)\Leftrightarrow(ii)$ holds by Corollary \ref{Bigorelationcor}.

$(ii)\Leftrightarrow(iii)$ First note that $\varphi_{\omega_M}\preceq\varphi_{\omega_L}$ precisely means
$$\exists\;K,D\ge 1\;\forall\;t\ge 0:\;\;\;\omega_L(e^t)=\varphi_{\omega_L}(t)\le K\varphi_{\omega_M}(t)+D=K\omega_M(e^t)+D.$$
Since $\omega_M(t)=\omega_L(t)=0$ for all $0\le t\le 1$ by normalization of $M$ and $L$ we get $\omega_L(t)\le K\omega_M(t)+D$ for all $t\ge 0$, i.e. $\omega_L(t)=O(\omega_M(t))$ and so $\omega_M\preceq\omega_L$. Similarly, $\omega_M\preceq\omega_L$ implies $\varphi_{\omega_M}\preceq\varphi_{\omega_L}$ as well.

Consequently, the desired equivalence $(ii)\Leftrightarrow(iii)$ follows by the first part in \cite[Lemma 6.5]{PTTvsmatrix}.
\qed\enddemo

Next let us gather some more immediate consequences concerning relation \hyperlink{ompreceq}{$\preceq_{\mathfrak{c}}$}.

\begin{remark}\label{relationconnectionrem}
Let $M,L\in\hyperlink{LCset}{\mathcal{LC}}$ be given.

\begin{itemize}
\item[$(i)$] If $L\le M$, then by definition $\omega_M(t)\le\omega_L(t)$ for all $t\ge 0$ and so $\varphi_{\omega_M}(t)\le\varphi_{\omega_L}(t)$ for all $t\ge 0$. In view of Corollary \ref{princor} we get $F_M\hyperlink{ompreceq}{\preceq_{\mathfrak{c}}}F_L$.

\item[$(ii)$] More generally, if $L\hyperlink{preceq}{\preceq}M$, then by definition and since $M_0=N_0=1$ we have $\omega_M(t)\le\omega_L(ht)$ for some $h>1$ and all $t\ge 0$. Hence
    $$\exists\;s_h>0\;\forall\;s\ge s_h:\;\;\;\varphi_{\omega_M}(s)=\omega_M(e^s)\le\omega_L(e^{s+\log(h)})=\varphi_{\omega_L}(s+\log(h))\le\varphi_{\omega_L}(sh),$$ since $\varphi_{\omega_L}$ is non-decreasing and $s+\log(h)\le sh\Leftrightarrow\log(h)\le s(h-1)$ for all $s$ large enough.

    This verifies $\varphi_{\omega_M}\hyperlink{ompreceq}{\preceq_{\mathfrak{c}}}\varphi_{\omega_L}$ and Corollary \ref{princor} implies again $F_M\hyperlink{ompreceq}{\preceq_{\mathfrak{c}}}F_L$.

\item[$(iii)$] Consequently, equivalent weight sequences yield equivalent associated N-functions and, in particular, this applies to the situation described in Remark \ref{strictincrem}.

\item[$(iv)$] However, the converse implication in $(iii)$ is not true in general. For this recall that by \eqref{powersub} we get
\begin{equation*}\label{powersubvarphi}
\forall\;\ell>0\;\forall\;t\ge 0:\;\;\;\varphi_{\omega_{M^{\ell}}}(t)=\omega_{M^{\ell}}(e^t)=\ell\omega_M(e^{t/\ell})=\ell\varphi_{\omega_M}(t/\ell).
\end{equation*}
Then, by following the arguments in $(i)$ in Remark \ref{convexityrem}, we see that
\begin{equation}\label{powersubvarphiequiv}
\forall\;\ell>0:\;\;\;\varphi_{\omega_M}\hyperlink{sim}{\sim_{\mathfrak{c}}}\varphi_{\omega_{M^{\ell}}},
\end{equation}
hence by Corollary \ref{princor} also $F_M\hyperlink{sim}{\sim_{\mathfrak{c}}}F_{M^{\ell}}$ for any $\ell>0$. However, for $\ell\neq 1$ the sequences $M$ and $M^{\ell}$ are not equivalent since $\lim_{j\rightarrow+\infty}(M_j)^{1/j}=+\infty$.

\item[$(v)$] In particular, $(iv)$ applies to the Gevrey-sequence $M\equiv G^s$, $s>0$. It is known that $\omega_{G^s}\sim t\mapsto t^{1/s}$, i.e. $\varphi_{\omega_{G^s}}\sim t\mapsto e^{t/s}$ (see the proof of Theorem \ref{comparisonprop0}). So \eqref{powersubvarphiequiv} is verified but clearly $G^s$ is not equivalent to $G^{s'}$ if $s\neq s'$.
\end{itemize}
\end{remark}

The next result provides a second characterization for $F_M\hyperlink{ompreceq}{\preceq_{\mathfrak{c}}}F_L$ in terms of a growth relation between $M$ and $L$ directly.

\begin{theorem}\label{comparisonprop}
Let $M,L\in\hyperlink{LCset}{\mathcal{LC}}$ be given. Consider the following assertions:
\begin{itemize}
\item[$(i)$] $L\hyperlink{preceq}{\preceq}M$ is valid.

\item[$(ii)$] The associated N-functions $F_M$ and $F_L$ (see Corollary \ref{princor}) satisfy
$$F_M\hyperlink{ompreceq}{\preceq_{\mathfrak{c}}}F_L,$$
equivalently $\varphi_{\omega_M}\hyperlink{ompreceq}{\preceq_{\mathfrak{c}}}\varphi_{\omega_L}$ is valid.

\item[$(iii)$] The sequences $M$ and $L$ satisfy
\begin{equation}\label{comparisonpropequ}
\exists\;c\in\NN_{>0}\;\exists\;A\ge 1\;\forall\;j\in\NN:\;\;\;L_j\le AM_{cj}.
\end{equation}
\end{itemize}
Then $(i)\Rightarrow(ii)\Rightarrow(iii)$ holds. If either $M$ or $L$ has in addition \hyperlink{mg}{$(\on{mg})$}, then also $(iii)\Rightarrow(ii)$.
\end{theorem}

By $(iv)$ in Remark \ref{relationconnectionrem} the implication $(i)\Rightarrow(ii)$ cannot be reversed in general.

\demo{Proof}
$(i)\Rightarrow(ii)$ This is shown in $(i),(ii)$ in Remark \ref{relationconnectionrem}.\vspace{6pt}

$(ii)\Rightarrow(iii)$ By Lemma \ref{relationrem} relation $F_M\hyperlink{ompreceq}{\preceq_{\mathfrak{c}}}F_L$ is equivalent to
$$\exists\;K\ge 1\;\exists\;C\ge 1\;\forall\;t\ge 0:\;\;\;F_M(t)\le F_L(Kt)+C,$$
and so by Corollary \ref{princor} (take w.l.o.g. $K\in\NN_{>0}$)
$$\exists\;K\in\NN_{>0}\;\exists\;D\ge 1\;\forall\;t\ge 0:\;\;\;\omega_M(e^t)=\varphi_{\omega_M}(t)\le\varphi_{\omega_L}(Kt)+D=\omega_L(e^{tK})+D.$$
We set $s:=e^t$ and hence this is equivalent to
$$\exists\;K\in\NN_{>0}\;\exists\;D\ge 1\;\forall\;s\ge 1:\;\;\;\omega_M(s)\le\omega_L(s^{K})+D.$$
Recall that $M,L\in\hyperlink{LCset}{\mathcal{LC}}$ implies (by normalization) $\omega_M(s)=\omega_L(s)=0$ for all $s\in[0,1]$ and so the previous estimate holds for any $s\ge 0$ (with the same constants). Thus \eqref{Prop32Komatsu} yields for all $j\in\NN$:
\begin{align*}
M_{Kj}&=\sup_{t\ge 0}\frac{t^{Kj}}{\exp(\omega_M(t))}\ge\frac{1}{e^{D}}\sup_{t\ge 0}\frac{t^{Kj}}{\exp(\omega_L(t^{K}))}=\frac{1}{e^{D}}\sup_{s\ge 0}\frac{s^j}{\exp(\omega_L(s))}=\frac{1}{e^{D}}L_j,
\end{align*}
and we are done when taking $c:=K$ and $A:=e^{D}$.\vspace{6pt}

$(iii)\Rightarrow(ii)$ Assume that \hyperlink{mg}{$(\on{mg})$} holds for $M$. By assumption \eqref{comparisonpropequ} and an iterated application of \hyperlink{mg}{$(\on{mg})$} we get
\begin{equation}\label{comparisonpropauxil}
\exists\;c\in\NN_{>0}\;\exists\;A,B\ge 1\;\forall\;k\in\NN:\;\;\;L_k\le AM_{ck}\le AB^kM_k^c,
\end{equation}
thus by definition of associated weight functions $\omega_{M^c}(t)\le\omega_L(Bt)+\log(A)$ for all $t\ge 0$. Consequently, since $\varphi_{\omega_L}(t)\rightarrow+\infty$ and by \eqref{equ114} we get for all $t$ large enough $$\varphi_{\omega_{M^c}}(t)\le\varphi_{\omega_L}(t+\log(B))+\log(A)\le\varphi_{\omega_L}(2t)+\log(A)\le 2\varphi_{\omega_L}(2t)\le\varphi_{\omega_L}(4t),$$
hence $\varphi_{\omega_{M^c}}\hyperlink{ompreceq}{\preceq_{\mathfrak{c}}}\varphi_{\omega_L}$. By \eqref{powersubvarphiequiv} we have that $\varphi_{\omega_{M^c}}\hyperlink{sim}{\sim_{\mathfrak{c}}}\varphi_{\omega_M}$ and so $\varphi_{\omega_M}\hyperlink{ompreceq}{\preceq_{\mathfrak{c}}}\varphi_{\omega_L}$ is verified. Corollary \ref{princor} yields the conclusion.

If $L$ has in addition \hyperlink{mg}{$(\on{mg})$}, then by \eqref{comparisonpropequ} and iterating \hyperlink{mg}{$(\on{mg})$} first we get
$$\exists\;c\in\NN_{>0}\;\exists\;A,B_1\ge 1\;\forall\;j\in\NN:\;\;\;L_{cj}\le B_1^{cj}L_j^c\le A^cB_1^{cj}M^c_{cj}.$$
Thus \eqref{comparisonpropauxil} is verified for all $k\in\NN$ with $k=cj$, $j\in\NN$ arbitrary. For the remaining cases let $k$ with $cj<k<cj+c$ for some $j\in\NN$ and then, since both $M$ and $L$ are also non-decreasing, we get for some $C\ge 1$
\begin{align*}
L_k&\le L_{cj+c}\le C^{c(j+1)}L_cL_{cj}\le C^{c(j+1)}L_c A^cB_1^{cj}M^c_{cj}\le(AC)^cL_c(B_1C)^kM^c_k.
\end{align*}
Summarizing, \eqref{comparisonpropauxil} is verified for all $k\in\NN$ and the rest follows as above.
\qed\enddemo

Thus we have the following characterization:

\begin{corollary}\label{comparisonpropcor}
Let $M,L\in\hyperlink{LCset}{\mathcal{LC}}$ be given. Assume that either $M$ or $L$ has in addition \hyperlink{mg}{$(\on{mg})$}, then the following are equivalent:
\begin{itemize}
\item[$(i)$] The associated N-functions $F_M$ and $F_L$ are equivalent.

\item[$(ii)$] The functions $\varphi_{\omega_M}$ and $\varphi_{\omega_L}$ are equivalent.

\item[$(iii)$] The sequences $M$ and $L$ satisfy
$$\exists\;c,d\in\NN_{>0}\;\exists\;A,B\ge 1\;\forall\;j\in\NN:\;\;\;M_j\le AL_{cj},\;\;\;\;\;L_j\le BM_{dj}.$$
\end{itemize}
\end{corollary}

The proof of $(ii)\Rightarrow(iii)$ in Theorem \ref{comparisonprop} transfers to relation $\preceq\preceq_{\mathfrak{c}}$ immediately: If $F_M\preceq\preceq_{\mathfrak{c}}F_L,$ equivalently if $\varphi_{\omega_M}\preceq\preceq_{\mathfrak{c}}\varphi_{\omega_L}$, then the sequences $M$ and $L$ satisfy
\begin{equation*}
\forall\;c\in\NN_{>0}\;\exists\;A\ge 1\;\forall\;j\in\NN:\;\;\;L_j\le AM_{cj}.
\end{equation*}

We finish by comparing the characterizing conditions for $M$ and $L$ in the previous results.

\begin{remark}\label{diffrelrem}
Let $M,L\in\hyperlink{LCset}{\mathcal{LC}}$ be given.
\begin{itemize}
\item[$(*)$] We have $L_j\ge 1$ for all $j\in\NN$ and so \eqref{comparisonprop0equ} implies \eqref{comparisonpropequ} (with $M$, $L$ interchanged).

\item[$(*)$] On the other hand, recall that by \eqref{equ114} we get that $F_M\preceq F_L$ implies $F_L\hyperlink{ompreceq}{\preceq_{\mathfrak{c}}}F_M$. Summarizing,
\begin{equation}\label{diffrelremequ}
\eqref{comparisonprop0equ}\Longleftrightarrow F_M\preceq F_L\Longrightarrow F_L\hyperlink{ompreceq}{\preceq_{\mathfrak{c}}}F_M\Longrightarrow\eqref{comparisonpropequ},
\end{equation}
and the last implication can be reversed provided that either $M$ or $L$ has in addition \hyperlink{mg}{$(\on{mg})$}.
\end{itemize}
\end{remark}

\subsection{On sufficiency conditions}\label{sufficiencysection}
We want to find some sufficient conditions for given sequences $M$, $L$ in order to ensure $F_M\hyperlink{sim}{\sim_{\mathfrak{c}}}F_L$. More precisely, the aim is to ensure an asymptotic behavior of the counting functions near infinity (see \eqref{comparisonlemma1equ}) which might be important for applications in the ultradifferentiable setting as well.

\begin{lemma}\label{comparisonlemma1}
Let $M,L\in\hyperlink{LCset}{\mathcal{LC}}$ be given. Assume that
\begin{equation}\label{comparisonlemma1equ}
\exists\;b\in(0,+\infty):\;\;\;\lim_{t\rightarrow+\infty}\frac{\Sigma_M(t)}{\Sigma_L(t)}=b.
\end{equation}
Then we get $F_M\hyperlink{sim}{\sim_{\mathfrak{c}}}F_L$ (equivalently $\varphi_{\omega_M}\hyperlink{sim}{\sim_{\mathfrak{c}}}\varphi_{\omega_L}$).
\end{lemma}

{\itshape Note:} For $M=L$ property \eqref{comparisonlemma1equ} holds trivially with $b=1$.

\demo{Proof}
By \eqref{princorequ} and \eqref{comparisonlemma1equ} we have $\lim_{t\rightarrow+\infty}\frac{f_M(t)}{f_L(t)}=b>0$ and so \cite[Lemma 3.2]{orliczbook} yields $F_M\hyperlink{sim}{\sim_{\mathfrak{c}}}F_L$.
\qed\enddemo

Let us now study condition \eqref{comparisonlemma1equ} in more detail.

\begin{lemma}\label{comparisonlemma2}
Let $M,L\in\hyperlink{LCset}{\mathcal{LC}}$ be given. Assume that
\begin{equation}\label{comparisonlemma2equ}
\exists\;c,j_0\in\NN_{>0}\;\forall\;j\ge j_0,\;\text{s.th.}\;\mu_j<\mu_{j+1},\;\exists\;d_j\in\NN_{>0}:\;\;\;\lambda_{cj}\le\mu_j<\mu_{j+1}\le\lambda_{cj+d_j},
\end{equation}
with $d_j\in\NN_{>0}$ such that $\lim_{j\rightarrow+\infty}\frac{d_j}{j}=0$. Then \eqref{comparisonlemma1equ} holds (with $b=\frac{1}{c}$).
\end{lemma}

\demo{Proof}
For all $t\ge\mu_{j_0}$ we find $j\ge j_0$ such that $\mu_j\le t<\mu_{j+1}$. Then $\Sigma_M(t)=j$ and by \eqref{comparisonlemma2equ} we get $cj\le\Sigma_L(t)<cj+d_j$. Thus
$$\frac{j}{cj+d_j}<\frac{\Sigma_M(t)}{\Sigma_L(t)}\le\frac{j}{cj},\;\;\;\forall\;\mu_j\le t<\mu_{j+1},\;j\ge j_0,$$
and since $\lim_{j\rightarrow+\infty}\frac{d_j}{j}=0$ we get \eqref{comparisonlemma1equ} with $b:=\frac{1}{c}$.
\qed\enddemo

Note that it is enough to require the existence of a sequence $d_j$ for $j$ such that $\mu_j<\mu_{j+1}$: If $j\ge j_0$ and $\mu_j=\mu_{j+1}=\dots=\mu_{j+\ell}<\mu_{j+\ell+1}$ for some $\ell\in\NN_{>0}$, then by \eqref{comparisonlemma2equ} we get
$$\lambda_{cj}\le\lambda_{cj+c\ell}\le\mu_j=\mu_{j+1}=\dots=\mu_{j+\ell}<\mu_{j+\ell+1}\le\lambda_{cj+c\ell+d_{j+\ell}},$$
and so for $t$ with $\mu_{j+\ell}\le t<\mu_{j+\ell+1}$ we have $\Sigma_M(t)=j+\ell$ and $cj+c\ell=c(j+\ell)\le\Sigma_L(t)<cj+c\ell+d_{j+\ell}=c(j+\ell)+d_{j+\ell}$ yielding the same estimate as above. (On the other hand, by switching to an equivalent sequence, we can assume $\mu_j<\mu_{j+1}$ for all $j\in\NN_{>0}$, see Remark \ref{strictincrem}.)\vspace{6pt}

We prove now the following characterization for \eqref{comparisonlemma1equ}.

\begin{lemma}\label{comparisonlemma3}
Let $M,L\in\hyperlink{LCset}{\mathcal{LC}}$ be given such that $\mu_j<\mu_{j+1}$ for all $j\in\NN_{>0}$. Then the following are equivalent:
\begin{itemize}
\item[$(i)$] We have that
\begin{equation}\label{comparisonlemma3equ1}
\exists\;b\in(0,+\infty):\;\;\;\lim_{t\rightarrow+\infty}\frac{\Sigma_M(t)}{\Sigma_L(t)}=b,
\end{equation}
i.e. \eqref{comparisonlemma1equ}.

\item[$(ii)$] We have that
\begin{equation}\label{comparisonlemma3equ}
\exists\;b\in(0,+\infty)\;\exists\;d\in\NN\;\forall\;0<c_1<b<c_2\;\exists\;j_0\in\NN_{>0}\;\forall\;j\ge j_0:\;\;\;\lambda_{\lceil j/c_2\rceil-d}\le\mu_j<\lambda_{\lceil j/c_1\rceil+d}.
\end{equation}

\item[$(iii)$] We have that
\begin{equation}\label{comparisonlemma3equ0}
\exists\;b\in(0,+\infty)\;\forall\;0<c_1<b<c_2\;\exists\;j_0\in\NN_{>0}\;\forall\;j\ge j_0\;\exists\;d_j\in\NN_{>0}:\;\;\;\lambda_{\lceil j/c_2\rceil-d_j}\le\mu_j<\lambda_{\lceil j/c_1\rceil+d_j},
\end{equation}
with $d_j\in\NN_{>0}$ such that $\lim_{j\rightarrow+\infty}\frac{d_j}{j}=0$.
\end{itemize}
The proof shows that in all assertions the value $b$ is the same.
\end{lemma}

\demo{Proof}
$(i)\Rightarrow(ii)$ By assumption,
$$\forall\;\epsilon>0\;\exists\;t_{\epsilon}>0\;\forall\;t\ge t_{\epsilon}:\;\;\;\Sigma_L(t)(b-\epsilon)<\Sigma_M(t)<(b+\epsilon)\Sigma_L(t),$$
thus we find $j_{\epsilon}\in\NN_{>0}$ large such that $\mu_{j_{\epsilon}}\ge t_{\epsilon}$ so that for all $t$ with $\mu_j\le t<\mu_{j+1}$ for some $j\ge j_{\epsilon}$ we get
\begin{equation}\label{comparisonlemma3proofequ}
\Sigma_L(t)(b-\epsilon)<j=\Sigma_M(t)<(b+\epsilon)\Sigma_L(t).
\end{equation}
The first estimate in \eqref{comparisonlemma3proofequ} yields $\Sigma_L(t)<\frac{j}{b-\epsilon}\le\lceil\frac{j}{b-\epsilon}\rceil$ and since $\Sigma_L(t)\in\NN$ we have $\Sigma_L(t)\le\lceil\frac{j}{b-\epsilon}\rceil-1$. Let now $c_1<b$ be arbitrary but fixed and take $\epsilon_1$ small enough to ensure $\frac{1}{b-\epsilon_1}\le\frac{1}{c_1}\Leftrightarrow c_1\le b-\epsilon_1$. Note that the choice of $\epsilon_1$ is only depending on chosen $c_1$ but not on $j$.

Thus $\Sigma_L(t)<\lceil\frac{j}{c_1}\rceil$ and so $t<\lambda_{\lceil\frac{j}{c_1}\rceil}$. We take $t:=\mu_j$ and get $\mu_j<\lambda_{\lceil\frac{j}{c_1}\rceil}$, i.e. the second half of \eqref{comparisonlemma3equ} for all $j\ge j_{\epsilon_1}$ (since $\mu_j<\mu_{j+1}$ for all $j$) with $d:=0$.\vspace{6pt}

Similarly, the second estimate in \eqref{comparisonlemma3proofequ} gives $\lfloor\frac{j}{b+\epsilon}\rfloor\le\frac{j}{b+\epsilon}<\Sigma_L(t)$ and so $t\ge\lambda_{\lfloor\frac{j}{b+\epsilon}\rfloor}$. We let $c_2>b$ be arbitrary but fixed and choose $\epsilon_2>0$ small enough to ensure $\frac{1}{b+\epsilon_2}\ge\frac{1}{c_2}\Leftrightarrow c_2\ge b+\epsilon_2$ and so $\lfloor\frac{j}{b+\epsilon_2}\rfloor\ge\lceil\frac{j}{b+\epsilon_2}\rceil-1\ge\lceil\frac{j}{c_2}\rceil-1$. Take $t:=\mu_j$ to get the first half of \eqref{comparisonlemma3equ} for all $j\ge j_{\epsilon_2}$ with $d:=1$. Summarizing, we are done when taking $j_0:=\max\{j_{\epsilon_1}, j_{\epsilon_2}\}$ and note that $d$ can be taken uniformly and not depending on chosen $c_1,c_2$.\vspace{6pt}

$(ii)\Rightarrow(iii)$ This is clear.\vspace{6pt}

$(iii)\Rightarrow(i)$ Let $b\in(0,+\infty)$ be given and take arbitrary (close) $0<c_1<b<c_2$ but from now on fixed. Let $t\ge\mu_{j_0}$ and so $\mu_j\le t<\mu_{j+1}$ for some $j\ge j_0$. Then \eqref{comparisonlemma3equ0} gives
$$\lambda_{\lceil j/c_2\rceil-d_j}\le\mu_j\le t<\mu_{j+1}<\lambda_{\lceil (j+1)/c_1\rceil+d_{j+1}}.$$
So $\Sigma_M(t)=j$ and the first estimate implies $\Sigma_L(t)\ge\lceil\frac{j}{c_2}\rceil-d_j\ge\frac{j}{c_2}-d_j$, whereas the last estimate yields $\Sigma_L(t)\le\lceil\frac{j+1}{c_2}\rceil+d_{j+1}-1\le\frac{j+1}{c_1}+d_{j+1}$. Summarizing, for all such $t$ we obtain
$$\frac{j}{j/c_1+1/c_1+d_{j+1}}\le\frac{\Sigma_M(t)}{\Sigma_L(t)}\le\frac{j}{j/c_2-d_j}.$$
Then note that $\frac{d_{j+1}}{j}=\frac{d_{j+1}}{j+1}\frac{j+1}{j}\rightarrow 0$ as $j\rightarrow+\infty\Leftrightarrow t\rightarrow+\infty$. Hence, as $c_1,c_2\rightarrow b$ we get that $\frac{\Sigma_M(t)}{\Sigma_L(t)}\rightarrow b$ as $t\rightarrow+\infty$.
\qed\enddemo

\section{From N-functions to associated weight sequences}\label{fromNtoassofct}
The aim is to reverse the construction from Section \ref{fromassotoN}; i.e. we start with an abstractly given N-function $G$, associate to it a sequence $M^G$ and study then the relation between $G$ and the associated N-function $F_{M^G}$.

Let $F$ be an N-function and first we introduce
\begin{equation}\label{Flogweight}
\omega_F(t):=F(\log(t)),\;\;\;t\ge 1,\hspace{15pt}\omega_F(t):=0,\;\;\;0\le t<1\;\;\;\text{(normalization)}.
\end{equation}
By the properties of $F$ we have that $\omega_F$ belongs to the class \hyperlink{omset0}{$\mathcal{W}_0$}: Concerning \hyperlink{om4}{$(\omega_4)$} note that $\varphi_{\omega_F}(t)=\omega_F(e^t)=F(t)$ for all $t\ge 0$, concerning \hyperlink{om3}{$(\omega_3)$} we remark that $\frac{\omega_F(t)}{\log(t)}=\frac{\omega_F(e^s)}{s}=\frac{F(s)}{s}$ for all $t>1\Leftrightarrow s>0$ and the last quotient tends to $+\infty$ as $s\rightarrow+\infty$, see \eqref{equ115116}.

{\itshape Note:} When given $\omega\in\hyperlink{omset0}{\mathcal{W}_0}$, then one can put
\begin{equation}\label{Flogweightinverse}
F_{\omega}(t):=\varphi_{\omega}(|t|)=\omega(e^{|t|}),\;\;\;t\in\RR.
\end{equation}
$F_{\omega}$ satisfies all properties to be formally an N-function (see Definition \ref{Nfunctiondef}) except necessarily the first part in \eqref{equ115116} and also $F_{\omega}(t)>0$ for all $t\neq 0$ is not clear. Thus the set of all $N$-functions does not coincide with the class \hyperlink{omset0}{$\mathcal{W}_0$}. However, one can overcome this technical problem for $F_{\omega}$ by passing to an equivalent (associated) N-function when taking into account Theorem \ref{princthm} analogously as it has been done before with $\varphi_{\omega_M}$.\vspace{6pt}

We consider the so-called {\itshape Legendre-Fenchel-Young-conjugate} of $\varphi_{\omega_F}$ which is given by
\begin{equation}\label{LFYconjugate}
\varphi^{*}_{\omega_F}(s):=\sup\{|s|t-\varphi_{\omega_F}(t): t\ge 0\}=\sup\{|s|t-F(t): t\ge 0\},\;\;\;s\in\RR.
\end{equation}
$\varphi_{\omega_F}$ is non-decreasing, convex by assumption, $\varphi_{\omega_F}(0)=\omega_F(e^0)=F(0)=0$ and $\lim_{t\rightarrow+\infty}\frac{\varphi_{\omega_F}(t)}{t}=\lim_{s\rightarrow+\infty}\frac{\omega_F(s)}{\log(s)}=+\infty$ as seen before. Thus we get, see e.g. \cite[Rem. 1.3]{BraunMeiseTaylor90}:

$\varphi^{*}_{\omega_F}$ is convex, $\varphi^{*}_{\omega_F}(0)=0$, $s\mapsto\frac{\varphi^{*}_{\omega_F}(s)}{s}$ is non-decreasing and tending to $+\infty$ as $s\rightarrow+\infty$. Finally $\varphi^{**}_{\omega_F}(s)=\varphi_{\omega_F}(s)$ holds for all $s\ge 0$.

Next let us introduce the associated sequence $M^F=(M^F_j)_j$ by
\begin{equation}\label{assosequ}
M^F_j:=\sup_{t>0}\frac{t^j}{\exp(\omega_F(t))}=\sup_{t\ge 1}\frac{t^j}{\exp(\omega_F(t))}=\sup_{t\ge 1}\frac{t^j}{\exp(F(\log(t)))},\;\;\;j\in\NN.
\end{equation}
The second equality holds by normalization and this formula should be compared with \eqref{Prop32Komatsu}. Hence for any $j\in\NN$ we get
\begin{align*}
M^F_j&=\sup_{t>0}\frac{t^j}{\exp(\omega_F(t))}=\exp\sup_{t>0}\left(j\log(t)-\omega_F(t)\right)=\exp\sup_{s\in\RR}\left(js-\omega_F(e^s)\right)=\exp\sup_{s\ge 0}\left(js-\omega_F(e^s)\right)
\\&
=\exp\sup_{s\ge 0}\left(js-\varphi_{\omega_F}(s)\right)=\exp(\varphi^{*}_{\omega_F}(j)),
\end{align*}
which implies $M^F\in\hyperlink{LCset}{\mathcal{LC}}$ by the properties of the conjugate. Note that by definition (normalization) one has $\omega_F(e^s)=0$ for all $-\infty<s\le 0$.

Finally, by \cite[Thm. 4.0.3]{dissertation} applied to $\omega=\omega_F$, see also \cite[Lemma 5.7]{compositionpaper} and \cite[Lemma 2.5]{sectorialextensions} with weight matrix parameter $x=1$ there, we obtain
\begin{equation}\label{assoweightequiv}
\exists\;C\ge 1\;\forall\;t\ge 0:\;\;\;\omega_{M^F}(t)\le\omega_F(t)\le 2\omega_{M^F}(t)+C,
\end{equation}
hence by \eqref{Flogweight}
\begin{equation}\label{assoweightequiv1}
\exists\;C\ge 1\;\forall\;t\ge 0:\;\;\;\varphi_{\omega_{M^F}}(t)\le\omega_F(e^t)=F(t)\le 2\varphi_{\omega_{M^F}}(t)+C.
\end{equation}

In order to avoid confusion let from now in this context $G$ be the given N-function, $M^G$ the sequence defined in \eqref{assosequ} and $F_{M^G}$ the associated N-function from Corollary \ref{princor} (applied to the sequence $M^G$).

\begin{definition}
Let $G$ be an N-function. Then the sequence $M^G$ is called the associated weight sequence.
\end{definition}

\begin{theorem}\label{assoweightequivlemma}
Let $G$ be an N-functions and let $F_{M^G}$ be the N-function associated with the sequence $M^G$. Then we get
\begin{equation}\label{assoweightequiv2}
\exists\;A,B\ge 1\;\forall\;t\ge 0:\;\;\;F_{M^G}(t)\le G(t)+A\le 2F_{M^G}(t)+B\le F_{M^G}(2t)+B,
\end{equation}
and this implies both $G\sim F_{M^G}\sim\varphi_{\omega_{M^G}}$ and $G\hyperlink{sim}{\sim_{\mathfrak{c}}}F_{M^G}\hyperlink{sim}{\sim_{\mathfrak{c}}}\varphi_{\omega_{M^G}}$.
\end{theorem}

\demo{Proof}
Corollary \ref{princor} applied to $M^G$ and \eqref{assoweightequiv1} yield
$$\exists\;C,C_1\ge 1\;\forall\;t\ge 0:\;\;\;G(t)\le 2\varphi_{\omega_{M^G}}(t)+C\le 2F_{M^G}(t)+C+C_1,$$
and
$$\exists\;C_2\ge 1\;\forall\;t\ge 0:\;\;\;F_{M^G}(t)\le\varphi_{\omega_{M^G}}(t)+C_2\le G(t)+C_2.$$
These estimates prove the first two parts in \eqref{assoweightequiv2} and the last one there follows by \eqref{equ114}, so by convexity and normalization of $F_{M^G}$, see also the estimate in $(i)$ in Remark \ref{convexityrem} applied to $a:=2$.

The desired relations follow now by \eqref{assoweightequiv2}, Lemma \ref{relationrem} and Corollary \ref{princor}.
\qed\enddemo

The next result gives the (expected) equivalence when starting with a weight sequence.

\begin{proposition}
Let $L\in\hyperlink{LCset}{\mathcal{LC}}$ be given and let $F_L$ be the associated N-function. Then we get
$$\exists\;A,B\ge 1\;\forall\;j\in\NN:\;\;\;\frac{1}{A}L_j\le M^{F_L}_j\le BL_j.$$
This estimate implies $M^{F_L}\hyperlink{approx}{\approx}L$, $F_L\hyperlink{sim}{\sim_{\mathfrak{c}}}F_{M^{F_L}}$ and $F_L\sim F_{M^{F_L}}$.
\end{proposition}

\demo{Proof}
We apply the previous constructions to the associated N-function $F_L$. For all $t\ge 1$ we have $\omega_{F_L}(t)=F_L(\log(t))$ and so by \eqref{equivwithCforM}
\begin{equation}\label{startingwithsequence}
\exists\;C,D\ge 1\;\forall\;t\ge 1:\;\;\;\omega_L(t)-C=\varphi_{\omega_L}(\log(t))-C\le\omega_{F_L}(t)\le\varphi_{\omega_L}(\log(t))+D=\omega_L(t)+D.
\end{equation}
Because $L\in\hyperlink{LCset}{\mathcal{LC}}$ we have $\omega_L(t)=0$ for all $t\in[0,1]$ (normalization) and $\omega_{F_L}(t)=0$ holds for all $t\in[0,1]$ by \eqref{Flogweight}. Thus \eqref{startingwithsequence} is valid for any $t\ge 0$. By combining \eqref{startingwithsequence} with \eqref{assosequ} and \eqref{Prop32Komatsu} we arrive at
$$\exists\;C,D\ge 1\;\forall\;j\in\NN:\;\;\;e^{-D}L_j\le M^{F_L}_j\le e^CL_j,$$
hence the conclusion.

This relation clearly implies $M^{F_L}\hyperlink{approx}{\approx}L$ and so, by Theorem \ref{comparisonprop} (recall also Remark \ref{relationconnectionrem}) the equivalence $F_L\hyperlink{sim}{\sim_{\mathfrak{c}}}F_{M^{F_L}}$. Moreover, \eqref{comparisonprop0equ} is verified with $c=1$ (pair-wise) and hence Theorem \ref{comparisonprop0} gives that $F_L\sim F_{M^{F_L}}$ as well.
\qed\enddemo

We continue with the following observations:

\begin{itemize}
\item[$(*)$] Theorem \ref{assoweightequivlemma} suggests that for abstractly given N-functions $G$ it is important to have information about $F_{M^G}$ resp. $\varphi_{\omega_{M^G}}$ and to study the associated weight sequence $M^G$. In particular, in view of \eqref{intreprvarphi} and \eqref{princorequ} the knowledge about the counting function $\Sigma_{M^G}$ is useful and this amounts to study the sequence of quotients $\mu^G$.

\item[$(*)$] On the other hand, when desired growth behaviours are expressed in terms of $\mu^G$, it is an advantage not to compute first $M^G$ via given $G$ and then the corresponding quotient sequence $\mu^G$ but to come up with a property for $G$ directly.

\item[$(*)$] This can be achieved by relating $\mu^G$ directly to $G$ as follows: Put
\begin{equation}\label{vfromF}
v_G(t):=\exp(-G(\log(t))),\;\;\;t\ge 1,\hspace{15pt}v_G(t):=1,\;\;\;0\le t<1,
\end{equation}
and so \eqref{assosequ} takes alternatively the following form:
$$M^G_j=\sup_{t>0}t^jv_G(t),\;\;\;j\in\NN.$$
This should be compared with \cite[Sect. 6]{solidassociatedweight}, \cite[Sect. 2.6, 2.7]{weightedentireinclusion}. Indeed, $M^G$ coincides with the crucial associated sequence $M^{v_G}$ in \cite{solidassociatedweight}, \cite{weightedentireinclusion}. By the assumptions on $G$ we get that $v_G$ is a (normalized and convex) weight in the notion of \cite{solidassociatedweight}, \cite{weightedentireinclusion}, i.e. in the setting of weighted spaces of entire functions, see also the literature citations in these papers.

\item[$(*)$] Let now $(t_j)_{j\in\NN}$ be the sequence such that $t_j\ge 0$ is denoting a/the global maximum point of the mapping $t\mapsto t^jv_G(t)$.

\item[$(*)$] The advantage when considering $v_G$ is that in \cite[$(6.5)$]{solidassociatedweight} we have derived the following relation between the sequence of quotients $\mu^G$ (set $\mu^G_0:=1$) and $(t_j)_{j\in\NN}$:
\begin{equation}\label{rmucomparison}
\forall\;j\in\NN:\;\;\;t_j\le\mu^G_{j+1}\le t_{j+1}.
\end{equation}
{\itshape Note:} For $j=0$ we have put $t_0:=1$ and so equality with $\mu^G_0$. In fact for $j=0$ we can choose any $t\in[0,1]$ as $t_0$ since $v_G$ is non-increasing and normalized, and by the convention $0^0:=1$. So $t_j\ge 1$ for any $j\in\NN$ because $j\mapsto t_j$ is non-decreasing and since $\lim_{j\rightarrow+\infty}\mu^G_j=+\infty$ we also get $t_j\rightarrow+\infty$. For concrete given $G$ (belonging to $\mathcal{C}^1$) the concrete computation for the values of $t_j$ might be not too difficult.
\end{itemize}

Then put
\begin{equation}\label{Fcounting}
\Sigma^G(t):=|\{j\in\NN_{>0}:\;\;\;t_j\le t\}|,\;\;\;t\ge 0,
\end{equation}
and $\Sigma^G:[0,+\infty)\rightarrow\NN$ is a right-continuous non-decreasing function with $\Sigma^G(t)=0$ for all $0\le t<1$ and $\Sigma^G(t)\rightarrow+\infty$ as $t\rightarrow+\infty$. By taking into account \eqref{rmucomparison} and the definition of $\Sigma_{M^G}$ and $\Sigma^G$ we get:

\begin{proposition}\label{assoseququotientcomparion}
Let $G$ be an N-function and let $M^G$ be the associated weight sequence. Then the counting functions $\Sigma^G$ and $\Sigma_{M^G}$ are related by
\begin{equation}\label{assoseququotientcomparionequ}
\forall\;t\ge 0:\;\;\;\Sigma_{M^G}(t)\le\Sigma^G(t)+1\le\Sigma_{M^G}(t)+1,
\end{equation}
which yields
$$\lim_{t\rightarrow+\infty}\frac{\Sigma_{M^G}(t)}{\Sigma^G(t)}=1.$$
\end{proposition}

Using $\Sigma^G$ we introduce
\begin{equation}\label{Gvequation}
\varphi^G(x):=\int_0^{|x|}\Sigma^G(e^t)dt,\;\;\;x\in\RR.
\end{equation}
Note that, analogously to the explanations for $\varphi_{\omega_M}$ given in Remark \ref{Nfunctionfailrem} in Section \ref{fromassotoN} we have that $\varphi^G$ is formally not an N-function since $t_j\ge 1$ for all $j\ge 1$ and so by definition $\Sigma^G(e^0)\ge 1\neq 0$ if $t_1=1$ or $\Sigma^G(e^t)=0$ for all $0\le t<\log(t_1)$ if $t_1>1$. We summarize the whole information in the final statement of this section:

\begin{theorem}\label{assoseququotientcomparioncor}
Let $G$ be an N-function, $M^G$ the associated weight sequence, $F_{M^G}$ the associated N-function and $\varphi^G$ given by \eqref{Gvequation}. Then
$$G\hyperlink{sim}{\sim_{\mathfrak{c}}}F_{M^G}\hyperlink{sim}{\sim_{\mathfrak{c}}}\varphi_{\omega_{M^G}}\hyperlink{sim}{\sim_{\mathfrak{c}}}\varphi^G,\hspace{15pt}G\sim F_{M^G}\sim\varphi_{\omega_{M^G}}\sim\varphi^G.$$
\end{theorem}

\demo{Proof}
The first two parts are shown in Theorem \ref{assoweightequivlemma}. The last one follows by Proposition \ref{assoseququotientcomparion}: We use \eqref{assoseququotientcomparionequ} and the representations \eqref{Gvequation} and \eqref{intreprvarphi} in order to get:
\begin{equation}\label{assoseququotientcomparioncorequ}
\exists\;C\ge 1\;\forall\;t\ge 0:\;\;\;\varphi^G(t)\le\varphi_{\omega_{M^G}}(t)\le\varphi^G(t)+t\le 2\varphi^G(t)+C\le\varphi^G(2t)+C.
\end{equation}
The second last estimate in \eqref{assoseququotientcomparioncorequ} holds since $\lim_{t\rightarrow+\infty}\frac{\varphi^G(t)}{t}=+\infty$ and the last one by \eqref{equ114} applied to $\varphi^G$. Both requirements on $\varphi^G$ are valid by the representation \eqref{Gvequation}, see Remark \ref{convexdivrem}. Then \eqref{assoseququotientcomparioncorequ} implies both $\varphi_{\omega_{M^G}}\hyperlink{sim}{\sim_{\mathfrak{c}}}\varphi^G$ and $\varphi_{\omega_{M^G}}\sim\varphi^G$.
\qed\enddemo

Theorem \ref{assoseququotientcomparioncor} shows that, up to equivalence, the whole information concerning growth and regularity properties of $G$ is already expressed by involving a certain associated weight sequence $M^G$ and its related/associated N-function $F_{M^G}$.

\section{On complementary N-functions and dual sequences}\label{complementarysection}
We give a connection between the so-called complementary N-functions $F^c$ and the dual sequences $D$ in the weight sequence setting. Note that in the theory of N-functions one naturally has the pair $(F,F^c$) and thus also $(F_M,F_M^c)$. Similarly for each $M\in\hyperlink{LCset}{\mathcal{LC}}$ one can naturally assign the dual sequence $D\in\hyperlink{LCset}{\mathcal{LC}}$ and we show that $F_M^c$ is closely related to $D$ via an integral representation formula using the counting-function $\log\circ\Sigma_D$.

\subsection{Complementary N-functions in the weight sequence setting}\label{complementarysubsection}
Let $F$ be an N-function given by the representation \eqref{Ndef} with right-derivative $f$. First, introduce
\begin{equation}\label{fcdef}
f^c(s):=\sup\{t\ge 0: f(t)\le s\},\;\;\;s\ge 0,
\end{equation}
thus $f^c$ is the right-inverse of $f$, i.e. it is the right-inverse of the right-derivative of $F$, see \cite[$(2.1)$]{orliczbook}. It is known that $f^c$ satisfies $(I)-(III)$ in Section \ref{Nfctdefsect}. \eqref{fcdef} should be compared with \cite[$(13)$, p. 10]{orliczbook1}; the difference appears since in \cite{orliczbook1} the authors work with the integral representation involving the left-continuous integrand/density (see Remark \ref{RaoRendefrem}).

\begin{definition}
The {\itshape complementary N-function} $F^c$, see \cite[Chapter I, \S 2, p. 11]{orliczbook}, is defined by
\begin{equation}\label{Ncomdef}
F^c(x):=\int_0^{|x|}f^c(t)dt,\;\;\;x\in\RR.
\end{equation}
\end{definition}

In \cite[Chapter I, $(2.9)$, p. 13]{orliczbook} it is mentioned that
\begin{equation}\label{Youngtrafo}
F^c(s)=\max_{t\ge 0}\{|s|t-F(t)\},\;\;\;s\in\RR,
\end{equation}
and this formula can be considered as an equivalent definition for $F^c$.

\begin{remark}\label{compequivrem}
We comment on the comparison of growth relations \hyperlink{ompreceq}{$\preceq_{\mathfrak{c}}$} and $\preceq$ between (associated) N-functions $F$, $G$ and their complementary N-functions $F^c$, $G^c$.
\begin{itemize}
\item[$(i)$] In \cite[Thm. 3.1, Thm. 3.2]{orliczbook} it has been shown that $F\hyperlink{ompreceq}{\preceq_{\mathfrak{c}}}G$ implies $G^c\hyperlink{ompreceq}{\preceq_{\mathfrak{c}}}F^c$ and by \cite[Sect. 2.2, Thm. 2 $(a)$]{orliczbook1} in fact this is an equivalence. By \cite[Sect. 2.2, Thm. 2 $(b)$]{orliczbook1} the corresponding statements holds w.r.t. relation $\preceq\preceq_{\mathfrak{c}}$ as well. In particular, two N-functions are equivalent if and only if their complementary N-functions so are.

\item[$(ii)$] If $F\preceq G$, then $G(t)\le CF(t)+D$ for some $C,D\ge 1$ and all $t\ge 0$ and so \eqref{Youngtrafo} gives for any $s\in\RR$
 \begin{align*}
 G^c(s)&=\max_{t\ge 0}\{|s|t-G(t)\}\ge\max_{t\ge 0}\{|s|t-CF(t)\}-D=C\max_{t\ge 0}\{C^{-1}|s|t-F(t)\}-D
 \\&
 =CF^c(s/C)-D,
 \end{align*}
see also the proof of \cite[Lemma 5.16]{compositionpaper}. This relation implies \eqref{implwithCextra} and so Lemma \ref{relationrem} yields that $F^c\hyperlink{ompreceq}{\preceq_{\mathfrak{c}}}G^c$. When one can choose $C=1$, then also $G^c\preceq F^c$ follows but in general this implication is not clear; see also Lemma \ref{delta2growthcomp}.
\end{itemize}
\end{remark}

By combining \eqref{LFYconjugate} and \eqref{Youngtrafo} we get
$$\varphi^{*}_{\omega_F}=F^c,$$
and since $M^F_j=\exp(\varphi^{*}_{\omega_F}(j))$, for this see \eqref{assosequ} and the computations below this equation, it follows that
$$\forall\;j\in\NN:\;\;\;M^F_j=\exp(F^c(j)).$$

Let now $M\in\hyperlink{LCset}{\mathcal{LC}}$ be given, then write $F^c_M$ and $f^c_M$ for the functions considered before w.r.t. to the associated N-function $F_M$ and the corresponding right-derivative $f_M$. Moreover, in view of \eqref{fcdef}, let us introduce
\begin{equation}\label{complementary}
\Gamma_M(s):=\sup\{t\ge 0: \Sigma_M(e^t)\le s\},\;\;\;s\ge 0.
\end{equation}
Finally we set
\begin{equation}\label{complementary1}
\varphi^c_{\omega_M}(x):=\int_0^{|x|}\Gamma_M(s)ds,\;\;\;x\in\RR.
\end{equation}

By \eqref{princorequ} it follows that
\begin{equation}\label{princorequcom}
\exists\;s_0>0\;\forall\;s\ge s_0:\;\;\;f^c_M(s)=\Gamma_M(s),
\end{equation}
because both functions are non-decreasing and $f_M(t)=\Sigma_M(e^t)$ for all $t\ge t_0$, with $t_0$ the value appearing in \eqref{princorequ}. So \eqref{princorequcom} is valid for all $s\ge s_0:=f_M(t_0)$. Using this identity we can prove the analogous result of Corollary \ref{princor} for the functions $F^c_M$ and $\varphi^c_{\omega_M}$.

\begin{proposition}\label{princorcom}
Let $M\in\hyperlink{LCset}{\mathcal{LC}}$ be given. Then we get
$$\exists\;C,D\ge 1\;\forall\;t\ge 0:\;\;\;\varphi^c_{\omega_M}(t)-C\le F^c_M(t)\le\varphi^c_{\omega_M}(t)+D.$$
This implies that both $\varphi^c_{\omega_M}\hyperlink{sim}{\sim_{\mathfrak{c}}}F^c_M$ and $\varphi^c_{\omega_M}\sim F^c_M$ hold.
\end{proposition}

\demo{Proof}
We use \eqref{princorequcom} and the representations \eqref{Ncomdef} and \eqref{complementary1} (recall also the arguments in the proof of Proposition \ref{princprop}). For all $s\ge s_0$ (with $s_0$ from \eqref{princorequcom}) we get
\begin{align*}
\varphi^c_{\omega_M}(s)&=\int_0^s\Gamma_M(t)dt=\int_0^{s_0}\Gamma_M(t)dt+\int_{s_0}^s\Gamma_M(t)dt=\int_0^{s_0}\Gamma_M(t)dt+\int_{s_0}^sf^c_M(t)dt
\\&
\le\int_0^{s_0}\Gamma_M(t)dt+\int_0^sf_M(t)dt=\varphi^c_{\omega_M}(s_0)+F_M^c(s),
\end{align*}
hence $\varphi^c_{\omega_M}(s)\le F_M^c(s)+C$ for all $s\ge 0$ when choosing $C:=\varphi^c_{\omega_M}(s_0)$. Similarly, for all $s\ge s_0$
\begin{align*}
F^c_M(s)&=\int_0^sf_M^c(t)dt=\int_0^{s_0}f_M^c(t)dt+\int_{s_0}^sf_M^c(t)dt=\int_0^{s_0}f_M^c(t)dt+\int_{s_0}^s\Gamma_M(t)dt
\\&
\le\int_0^{s_0}f_M^c(t)dt+\int_0^s\Gamma_M(t)dt=F^c_M(s_0)+\varphi^c_{\omega_M}(s),
\end{align*}
hence $F^c_M(s)\le\varphi^c_{\omega_M}(s)+D$ for all $s\ge 0$ when choosing $D:=F^c_M(s_0)$.
\qed\enddemo

On the other hand, formula \eqref{Youngtrafo} applied to $\varphi_{\omega_M}$ yields
\begin{equation}\label{Youngtrafoasusual}
\sup_{t\ge 0}\{|s|t-\varphi_{\omega_M}(t)\}=\varphi^{*}_{\omega_M}(s),\;\;\;s\in\RR,
\end{equation}
the {\itshape Legendre-Fenchel-Young-conjugate} of $\varphi_{\omega_M}$. By combining \eqref{equivwithCforM}, \eqref{Youngtrafo} applied to $F_M$ and \eqref{Youngtrafoasusual} we get
\begin{equation}\label{princorcomequ}
\exists\;C,D\ge 1\;\forall\;s\in\RR:\;\;\;-C+F^c_M(s)\le\varphi^{*}_{\omega_M}(s)\le F^c_M(s)+D,
\end{equation}
see the estimate in $(ii)$ in Remark \ref{compequivrem}. Hence $F^c_M\hyperlink{sim}{\sim_{\mathfrak{c}}}\varphi^{*}_{\omega_M}$ and $F^c_M\sim\varphi^{*}_{\omega_M}$ hold. When combining this with Proposition \ref{princorcom} we arrive at the following result:

\begin{theorem}\label{princorcom1}
Let $M\in\hyperlink{LCset}{\mathcal{LC}}$ be given. Then $$\varphi^c_{\omega_M}\hyperlink{sim}{\sim_{\mathfrak{c}}}F^c_M\hyperlink{sim}{\sim_{\mathfrak{c}}}\varphi^{*}_{\omega_M},\hspace{15pt}\varphi^c_{\omega_M}\sim F^c_M\sim\varphi^{*}_{\omega_M}.$$
\end{theorem}

Finally, we apply Theorem \ref{princorcom1} to the associated weight sequence $M^G$.

\begin{corollary}\label{princorcom1cor}
Let $G$ be an N-function, let $M^G\in\hyperlink{LCset}{\mathcal{LC}}$ be the associated sequence defined via \eqref{assosequ} and $\varphi^G$ be given by \eqref{Gvequation}. Then
$$\varphi^{*}_{\omega_{M^G}}\hyperlink{sim}{\sim_{\mathfrak{c}}}\varphi^c_{\omega_{M^G}}\hyperlink{sim}{\sim_{\mathfrak{c}}}F^c_{M^G}\hyperlink{sim}{\sim_{\mathfrak{c}}}G^c\hyperlink{sim}{\sim_{\mathfrak{c}}}(\varphi^G)^c,\hspace{15pt}\varphi^{*}_{\omega_{M^G}}\sim\varphi^c_{\omega_{M^G}}\sim F^c_{M^G}.$$
\end{corollary}

\demo{Proof}
Concerning \hyperlink{sim}{$\sim_{\mathfrak{c}}$}, the first and the second equivalence hold by Theorem \ref{princorcom1}, the third and the fourth one by taking into account \eqref{assoweightequiv2}, Theorem \ref{assoseququotientcomparioncor} (see \eqref{assoseququotientcomparioncorequ}) and $(i)$ in Remark \ref{compequivrem}. Here $(\varphi^G)^c$ is given in terms of \eqref{Youngtrafo}.

For $\sim$ we use the same results, however $F^c_{M^G}\sim G^c$ and $G^c\sim(\varphi^G)^c$ are not clear in general: In order to conclude we want to apply \eqref{Youngtrafo} and so in the relations only an additive constant should appear, see $(ii)$ in Remark \ref{compequivrem}. However, in both \eqref{assoweightequiv2} and \eqref{assoseququotientcomparioncorequ} we also have the multiplicative constant $2$.
\qed\enddemo

\subsection{Complementary associated N-functions versus dual sequences}
In this section the aim is to see that $\Gamma_M$ in \eqref{complementary} and crucially appearing in the representation \eqref{complementary1} is closely connected to the counting function $\Sigma_D$, with $D$ denoting the so-called {\itshape dual sequence} of $M$. Moreover, \eqref{complementary} should be compared with the formula for the so-called {\itshape bidual sequence} of $M$ in \cite[Def. 2.1.41, Thm. 2.1.42]{dissertationjimenez}.

For given $M\in\hyperlink{LCset}{\mathcal{LC}}$ we introduce the {\itshape dual sequence} $D=(D_j)_j$ defined in terms of the corresponding quotient sequence $\delta$ as follows, see \cite[Def. 2.1.40, p. 81]{dissertationjimenez}:
\begin{equation}\label{dualdef}
\forall\;j\ge\mu_1(\ge 1):\;\;\;\delta_{j+1}:=\Sigma_M(j),\hspace{20pt}\delta_{j+1}:=1\;\;\;\forall\;j\in\ZZ,\;-1\le j<\mu_1.
\end{equation}
Then we set
$$D_j:=\prod_{k=0}^j\delta_k,\;\;\;j\in\NN,$$
hence $D\in\hyperlink{LCset}{\mathcal{LC}}$ with $1=D_0=D_1$ follows. Recall that by \cite[Def. 2.1.27]{dissertationjimenez} the function $\nu_{\mathbf{m}}$ in \cite{dissertationjimenez} precisely denotes the counting function $\Sigma_M$ (see \eqref{counting}) and note that concerning the sequence of quotients there exists an index-shift: more precisely we have $m_j\equiv\mu_{j+1}$ for all $j\in\NN$ with $\mathbf{m}=(m_j)_j$ used in \cite{dissertationjimenez} and \cite{index}. Moreover, a weight sequence in \cite{dissertationjimenez} means a sequence having all requirements from class \hyperlink{LCset}{$\mathcal{LC}$} except $M_0\le M_1$; see \cite[Sect. 1.1.1, Def. 1.1.8]{dissertationjimenez}. Finally, we mention that this notion should be compared with the definition of the counting function $\varphi$ on \cite[Sect. 1.2, p. 3]{orliczbook1}.\vspace{6pt}

Let $M\in\hyperlink{LCset}{\mathcal{LC}}$ be given, we analyze now $\Gamma_M$:\vspace{6pt}

\begin{itemize}
\item[$(*)$] Obviously, $\lim_{s\rightarrow+\infty}\Gamma_M(s)=+\infty$ and since $\lim_{j\rightarrow+\infty}\mu_j=+\infty$ we have $\mu_j<\mu_{j+1}$ for infinitely many indices $j$.

\item[$(*)$] Then recall that $\Sigma_M(e^t)\in\NN$ and that $\Sigma_M(e^t)=0$ for all $0\le t<\log(\mu_1)$. By normalization we get $\log(\mu_j)\ge\log(\mu_1)\ge\log(1)=0$ for all $j\in\NN_{>0}$.

\item[$(*)$] {\itshape Case I:} Assume that $1=\mu_1=\dots=\mu_d<\mu_{d+1}$ for some $d\in\NN_{>0}$. Then $e^t\ge\mu_d=1$ and so $\Sigma_M(e^t)\ge d$ for all $t\ge 0$. Hence the set of values $t\ge 0$ with $\Sigma_M(e^t)\le s$ in the definition of $\Gamma_M(s)$ is empty for all values $s$ with $0\le s<d$. In this case we put $\Gamma_M(s):=0$. Note that $d$ is finite because $\lim_{j\rightarrow+\infty}\mu_j=+\infty$.

    Let now $s$ be such that $j\le s<j+1$ for some $j\in\NN_{>0}$ with $j\ge d$. In view of Remark \ref{strictincrem} in general it is not clear that we can find $t\ge 0$ such that $\Sigma_M(e^t)=j$. In fact this holds if and only if $\mu_j<\mu_{j+1}$ because for all $t\ge 0$ with $\log(\mu_j)\le t<\log(\mu_{j+1})$ we get $\Sigma_M(e^t)=j$. We have $\Sigma(e^{\log(\mu_{j+1})})\ge j+1>s$ and consequently for all such indices $j$, in particular for $j=d$, we obtain $\Gamma_M(s)=\log(\mu_{j+1})$.

If $j\ge d$ is such that $\mu_j=\mu_{j+1}$, then $j\ge d+1$. Thus there exist $\ell,c\in\NN_{>0}$ with $d\le\ell\le j-1$ and such that $\mu_{\ell}<\mu_{\ell+1}=\mu_{\ell+2}=\dots=\mu_{\ell+c}=\mu_j=\mu_{j+1}$. For all $t$ with $\log(\mu_{\ell})\le t<\log(\mu_{\ell+1})$ we get $\Sigma_M(e^t)=\ell<j$ and since $\Sigma_M(e^{\log(\mu_{\ell+1})})=\Sigma_M(e^{\log(\mu_{j+1})})\ge j+1>s$ we have again $\Gamma_M(s)=\log(\mu_{j+1})$.

\item[$(*)$] {\itshape Case II:} Assume that $\mu_1>1$. First we have $\Gamma_M(s)=\log(\mu_1)(>0)$ for all $0\le s<1$ since $\Sigma_M(e^t)=0$ for all $0\le t<\log(\mu_1)$ and $\Sigma_M(e^{\log(\mu_1)})\ge 1$.

Let now $j\le s<j+1$ for some $j\in\NN_{>0}$. Similarly as above, if $\mu_j<\mu_{j+1}$ then we get $\Gamma_M(s)=\log(\mu_{j+1})$.

If $\mu_j=\mu_{j+1}$, then we distinguish: Either there exist $\ell,c\in\NN_{>0}$ with $1\le\ell\le j-1$ such that $\mu_{\ell}<\mu_{\ell+1}=\mu_{\ell+2}=\dots=\mu_{\ell+c}=\mu_j=\mu_{j+1}$. Then again $\Gamma_M(s)=\log(\mu_{j+1})$ by the same reasons as in Case I before.

Finally, if this choice is not possible, then this precisely means $\mu_1=\dots=\mu_j=\mu_{j+1}$. Since $\mu_1>1$ we have $\Sigma_M(e^t)=0$ for all $0\le t<\log(\mu_1)$ and $\Sigma_M(e^t)\ge j+1>s$ for all $t\ge\log(\mu_{j+1})=\log(\mu_1)$. Thus $\Gamma_M(s)=\log(\mu_{j+1})=\log(\mu_1)$ holds.
\end{itemize}
\vspace{6pt}

Summarizing, we have shown the following statement:

\begin{proposition}\label{GammaMprop}
Let $M\in\hyperlink{LCset}{\mathcal{LC}}$ be given with quotient sequence $(\mu_j)_j$.
\begin{itemize}
\item[$(a)$] If $1=\mu_1=\dots=\mu_d<\mu_{d+1}$ for some $d\in\NN_{>0}$, then
$$\Gamma_M(s)=0,\;\;\;\forall\;0\le s<d,\hspace{25pt}\Gamma_M(s)=\log(\mu_{j+1}),\;\;\;\;\forall\;j\le s<j+1,\;\forall\;j\ge d.$$

\item[$(b)$] If $\mu_1>1$, then
$$\Gamma_M(s)=\log(\mu_{j+1}),\;\;\;\forall\;j\le s<j+1,\;\forall\;j\in\NN.$$
\end{itemize}
\end{proposition}

Next let us study $\Sigma_D$ in detail, see also the proof of \cite[Thm. 2.1.42]{dissertationjimenez}:

\begin{itemize}
\item[$(*)$] Let $t\ge 0$, then the value $\Sigma_D(t)$ is equal to the maximal integer $j\in\NN_{>0}$ such that $\delta_j\le t$, if such $j$ exists, and otherwise equal to $0$. By definition in \eqref{dualdef} we have $\delta_j\in\NN_{>0}$ for all $j\in\NN$ and so $\Sigma_D(t)=0$ for $0\le t<1$.

\item[$(*)$] In fact $\delta_j=1$ for all $j\in\NN_{>0}$ such that $j<\mu_1+1$. The largest of those integers $j$ coincides with $\lfloor\mu_1+1\rfloor\ge 2$ if $\mu_1\notin\NN_{>0}$ (in this case $\mu_1>1$) and with $\mu_1$ if $\mu_1\in\NN_{>0}$. Moreover, $\delta_j=1$ for all $j$ such that $\mu_1+1\le j<\mu_2+1$ if such an integer $j$ exists (case I). In particular, this holds if $\mu_2\ge\mu_1+1$.

    If not, then $\delta_j=d\in\NN_{\ge 2}$ for the minimal integer $j$ with $j\ge\mu_1+1$ (case II). This occurs if for this minimal integer already $j\ge\mu_d+1$ is valid. Note that $d$ is finite since $\lim_{j\rightarrow+\infty}\mu_j=+\infty$ and let $d$ be such that $\mu_d+1\le j<\mu_{d+1}+1$ for $j$ minimal such that $j\ge\mu_1+1$.

\item[$(*)$] Consider $1\le t<2$ and distinguish: In the first case $\Sigma_D(t)$ is the largest integer $j$ with $\mu_1+1\le j<\mu_2+1$ and in the second case $\Sigma_D(t)$ coincides with the largest integer $j$ with $j<\mu_1+1$. Since $\mu_1+1\ge 2$ in both cases the existence of such an integer $j$ is ensured.

\item[$(*)$] More generally, in case I if $n\le t<n+1$ with $n\in\NN_{>0}$, then $\Sigma_D(t)$ is the largest integer $j$ such that $j<\mu_{n+1}+1$.

\item[$(*)$] In case II, if $n\le t<n+1$ for some $n\in\NN_{>0}$ with $n+1\le d$, then $\Sigma_D(t)$ is (still) the largest integer $j$ such that $j<\mu_1+1$ since for the minimal integer $j\ge\mu_1+1$ we have $\delta_{j}=\Sigma_M(j-1)\ge\Sigma_M(\mu_d)=d>t$. The last equality holds since $\mu_d<\mu_{d+1}$ by the choice of $d$.

    If $n\le t<n+1$ for some $n\in\NN_{>0}$ with $n\ge d$, then $\delta_j$ coincides with largest integer $j$ such that $j<\mu_{n+1}+1$. For $n=d$ recall that at least one integer $j$ exists with $\mu_d+1\le j<\mu_{d+1}+1$ and so $\delta_j=\Sigma_M(j-1)=d$.
\end{itemize}

Summarizing, we have shown the following statement (see also \cite[Thm. 2.1.42]{dissertationjimenez}):

\begin{proposition}\label{DeltaDprop}
Let $M\in\hyperlink{LCset}{\mathcal{LC}}$ be given and $D$ its dual sequence. Then
$$\Sigma_D(t)=0,\;\;\;\forall\;0\le t<1,$$
and:
\begin{itemize}
\item[$(a)$] If for the minimal integer $j$ such that $j\ge\mu_1+1$ one has $\mu_d+1\le j<\mu_{d+1}+1$ for some $d\in\NN_{>0}$, $d\ge 2$, then
\begin{equation*}\label{SigmaDequ}
\Sigma_D(t)=\max\{j\in\NN_{>0}:\;j<\mu_1+1\},\;\;\;\forall\;n\in\NN_{>0},\;n<d,\;\;\;\forall\;n\le t<n+1,
\end{equation*}
\begin{equation*}\label{SigmaDequ1}
\Sigma_D(t)=\max\{j\in\NN_{>0}:\;j<\mu_{n+1}+1\},\;\;\;\forall\;n\in\NN_{>0},\;n\ge d,\;\;\;\forall\;n\le t<n+1.
\end{equation*}

\item[$(b)$] If for the minimal integer $j$ such that $j\ge\mu_1+1$ one has $\mu_1+1\le j<\mu_2+1$, then
\begin{equation*}\label{SigmaDequ2}
\Sigma_D(t)=\max\{j\in\NN_{>0}:\;j<\mu_{n+1}+1\},\;\;\;\forall\;n\in\NN_{>0},\;\;\;\forall\;n\le t<n+1.
\end{equation*}
\end{itemize}
\end{proposition}

When combining Propositions \ref{GammaMprop} and \ref{DeltaDprop} we are able to establish a connection between the counting functions $\Gamma_M$ and $\Sigma_D$.

\begin{theorem}\label{GammaMvsDprop}
Let $M\in\hyperlink{LCset}{\mathcal{LC}}$ be given and $D$ its dual sequence. Then
\begin{equation}\label{GammaMvsDpropequ}
\exists\;t_0>0\;\forall\;t\ge t_0:\;\;\;\Gamma_M(t)\le\log(\Sigma_D(t))\le\Gamma_M(t)+1,
\end{equation}
and this estimate implies
$$\lim_{t\rightarrow+\infty}\frac{\Gamma_M(t)}{\log(\Sigma_D(t))}=1.$$
In \eqref{GammaMvsDpropequ} we can take $t_0:=d\in\NN_{>0}$ such that for the minimal integer $j\ge\mu_1+1$ we get $\mu_d+1\le j<\mu_{d+1}+1$.
\end{theorem}

\demo{Proof}
First, by Proposition \ref{DeltaDprop} we get that $\Sigma_D(t)$ coincides with the maximal integer $j<\mu_{n+1}+1$ for all $t\ge 0$ having $n\le t<n+1$ and all $n\ge d$ with $d\in\NN_{>0}$ such that for the minimal integer $j\ge\mu_1+1$ we get $\mu_d+1\le j<\mu_{d+1}+1$. In particular, for all such $n$ we get $\mu_{n+1}\ge\mu_{d+1}>\mu_1\ge 1$.

Then recall that by Proposition \ref{GammaMprop} we get $\Gamma_M(t)=\log(\mu_{n+1})$ for all $t$ such that $n\le t<n+1$ with $n\in\NN$ such that $\mu_{n+1}>1$. In particular, as mentioned before, this holds for all $n\ge d$.

So let $n\in\NN_{>0}$ be such that $n\ge d$. Take $t$ with $n\le t<n+1$, then one has $\Sigma_D(t)<\mu_{n+1}+1$ and
$$\log(\Sigma_D(t))<\log(\mu_{n+1}+1)\le\log(e\mu_{n+1})=\log(\mu_{n+1})+1=\Gamma_M(t)+1,$$
showing the second part of \eqref{GammaMvsDpropequ}.

On the other hand for all such $t$ we get $\Sigma_D(t)\ge\mu_{n+1}=\exp(\Gamma_M(t))$, in fact even $\Sigma_D(t)\ge\lceil\mu_{n+1}\rceil$ holds, and which proves the first estimate in \eqref{GammaMvsDpropequ}.
\qed\enddemo

Using the counting function $\Sigma_D$ we set
\begin{equation}\label{complementary2}
F_{\widetilde{\Gamma}_D}(x):=\int_0^{|x|}\widetilde{\Gamma}_D(s)ds,\;\;\;x\in\RR,
\end{equation}
with
$$\widetilde{\Gamma}_D(s):=\log(\Sigma_D(s)),\;\;\;s\ge 1,\hspace{15pt}\widetilde{\Gamma}_D(s):=0,\;\;\;0\le s<1.$$
$\widetilde{\Gamma}_D$ is non-decreasing and right-continuous and tending to infinity. Recall that $\Sigma_D(s)\in\NN_{>0}$ for all $s\ge 1$ and this technical modification is unavoidable: In \eqref{complementary2} we cannot consider $\log(\Sigma_D(s))$ directly for all $s\ge 0$ because $\Sigma_D(s)=0$ for $0\le s<1$ by definition. Then in view of \eqref{GammaMvsDpropequ} we get
\begin{equation}\label{GammaMvsDpropequ1}
\exists\;C,D\ge 1\;\forall\;t\ge 0:\;\;\;-C+\Gamma_M(t)\le\widetilde{\Gamma}_D(t)\le\Gamma_M(t)+D,
\end{equation}
which implies
$$\lim_{t\rightarrow+\infty}\frac{\Gamma_M(t)}{\widetilde{\Gamma}_D(t)}=1.$$
Note that $F_{\widetilde{\Gamma}_D}$ is formally not an N-function by analogous reasons as in Remark \ref{Nfunctionfailrem}: We have $\widetilde{\Gamma}_D(s)=0$ for (at least) all $0\le s<1$, see the proof of Proposition \ref{DeltaDprop}.

Now we are able to prove the main statement of this section.

\begin{theorem}\label{princorcom2}
We have the following equivalences:

\begin{itemize}
\item[$(i)$] Let $M\in\hyperlink{LCset}{\mathcal{LC}}$ be given and $D$ its dual sequence, then $$F_{\widetilde{\Gamma}_D}\hyperlink{sim}{\sim_{\mathfrak{c}}}\varphi^c_{\omega_M}\hyperlink{sim}{\sim_{\mathfrak{c}}}F^c_M\hyperlink{sim}{\sim_{\mathfrak{c}}}\varphi^{*}_{\omega_M},\hspace{15pt}F_{\widetilde{\Gamma}_D}\sim\varphi^c_{\omega_M}\sim F^c_M\sim\varphi^{*}_{\omega_M}.$$

\item[$(ii)$] Let $G$ be an N-function, $M^G\in\hyperlink{LCset}{\mathcal{LC}}$ the associated sequence defined via \eqref{assosequ} and $D^G$ the corresponding dual sequence. Finally, let $\varphi^G$ be given by \eqref{Gvequation}. Then
$$F_{\widetilde{\Gamma}_{D^G}}\hyperlink{sim}{\sim_{\mathfrak{c}}}\varphi^{*}_{\omega_{M^G}}\hyperlink{sim}{\sim_{\mathfrak{c}}}\varphi^c_{\omega_{M^G}}\hyperlink{sim}{\sim_{\mathfrak{c}}}F^c_{M^G}\hyperlink{sim}{\sim_{\mathfrak{c}}} G^c\hyperlink{sim}{\sim_{\mathfrak{c}}}(\varphi^G)^c,\hspace{15pt}F_{\widetilde{\Gamma}_{D^G}}\sim\varphi^{*}_{\omega_{M^G}}\sim\varphi^c_{\omega_{M^G}}\sim F^c_{M^G}.$$
\end{itemize}
\end{theorem}

\demo{Proof}
$(i)$ In view of Theorem \ref{princorcom1} only $F_{\widetilde{\Gamma}_D}\hyperlink{sim}{\sim_{\mathfrak{c}}}\varphi^c_{\omega_M}$ resp. $F_{\widetilde{\Gamma}_D}\sim\varphi^c_{\omega_M}$ has to be verified. This follows by using \eqref{GammaMvsDpropequ1} and the representations \eqref{complementary2} and \eqref{complementary1}. Note that these representations imply $\lim_{t\rightarrow+\infty}\frac{F_{\widetilde{\Gamma}_D}(t)}{t}=\lim_{t\rightarrow+\infty}\frac{\varphi^c_{\omega_M}(t)}{t}=+\infty$, see Remark \ref{convexdivrem}.\vspace{6pt}

$(ii)$ By Corollary \ref{princorcom1cor} it suffices to verify $F_{\widetilde{\Gamma}_{D^G}}\hyperlink{sim}{\sim_{\mathfrak{c}}}\varphi^c_{\omega_{M^G}}$ and $F_{\widetilde{\Gamma}_{D^G}}\sim\varphi^c_{\omega_{M^G}}$. Both relations follow by applying the first part to $M^G$ and $D^G$. Concerning $\sim$ note that here both functions are given directly by the representations \eqref{complementary2} resp. \eqref{complementary1}, so we are not involving formula \eqref{Youngtrafo} and since in \eqref{GammaMvsDpropequ1} only additive constants appear the problem described in the proof of Corollary \ref{princorcom1cor} (see $(ii)$ in Remark \ref{compequivrem}) does not occur.
\qed\enddemo

We remark that $F_{\widetilde{\Gamma}_D}$ does not coincide with $F_D$ directly but nevertheless the dual sequence $D$ can be used to get an alternative (equivalent) representation and description of the complementary N-function $F^c_M$.

\section{Growth and regularity conditions for associated N-functions}\label{conditionsection}
In the theory of Orlicz classes and Orlicz spaces several conditions for (abstractly given) N-functions appear frequently. The aim of this section is to study these known growth and regularity assumptions in the weight sequence setting in terms of given $M$.\vspace{6pt}

We remark that all appearing conditions are naturally preserved under \hyperlink{sim}{$\sim_{\mathfrak{c}}$} for (associated) N-functions, see the given citations in the forthcoming sections resp. Remark \ref{Delta2complementaryrem}. However, by inspecting the proofs one can see that this fact also holds for a wider class of functions, e.g. when having normalization, convexity, being non-decreasing and tending to infinity (in particular for $\varphi_{\omega_M}$).

Therefore recall that the technical failure of $\varphi_{\omega_M}$ to be formally an N-function occurs at the point $0$ (see Remark \ref{Nfunctionfailrem}) whereas for all crucial conditions under consideration large values $t\ge t_0>0$ are relevant.\vspace{6pt}

Of course, it makes also sense to consider the conditions in this section for arbitrary functions $F:[0,+\infty)\rightarrow[0,+\infty)$.

\subsection{The $\Delta_2$-condition}\label{Delta2section}
The most prominent property is the so-called $\Delta_2$-condition, see e.g. \cite[Chapter I, \S 4, p. 23]{orliczbook} and \cite[Sect. 2.3, Def. 1, p. 22]{orliczbook1}, which reads as follows:
\begin{equation}\label{delta2}
\limsup_{t\rightarrow+\infty}\frac{F(2t)}{F(t)}<+\infty.
\end{equation}
This growth-condition is precisely \hyperlink{om1}{$(\omega_1)$} for $F$ and thus also frequently appears for so-called weight functions $\omega$ in the sense of Braun-Meise-Taylor, see \cite{BraunMeiseTaylor90}. It is straight-forward that $\Delta_2$ is preserved under relation $\sim$ and also under \hyperlink{sim}{$\sim_{\mathfrak{c}}$}, see \cite[p. 23]{orliczbook}.

It is known, see e.g. \cite[Thm. 8.2]{orliczbook} and \cite[Thm. 3.4]{orliczAlexopoulos04}, that $\mathcal{L}_F$ is a linear space if and only if $F$ satisfies the $\Delta_2$-condition. Moreover, we also get:

\begin{lemma}\label{delta2growthcomp}
Let $F_1$ and $F_2$ be two N-functions such that either $F_1$ or $F_2$ has $\Delta_2$. Then $F_1\hyperlink{ompreceq}{\preceq_{\mathfrak{c}}}F_2$ if and only if $F_2\preceq F_1$.

In fact, in order to conclude, we only require that either $F_1$ or $F_2$ is normalized and convex and that either $F_1$ or $F_2$ has $\Delta_2$.
\end{lemma}

\demo{Proof}
By \eqref{equ114} we have that $F_1\preceq F_2$ implies $F_2\hyperlink{ompreceq}{\preceq_{\mathfrak{c}}}F_1$ (see Remark \ref{preceqrem}). The converse implication holds by an iterated application of $\Delta_2$ for either $F_1$ or $F_2$.
\qed\enddemo

In the weight sequence setting in view of Corollary \ref{princor} we require \eqref{delta2} (i.e. \hyperlink{om1}{$(\omega_1)$}) not for $\omega_M$ directly but for $\varphi_{\omega_M}$. Note that in \cite[Thm. 3.1]{subaddlike} we have already given a characterization of \hyperlink{om1}{$(\omega_1)$} for $\omega_M$ in terms of $M$ but $\Delta_2$ for $\varphi_{\omega_M}$ (resp. equivalently for $F_M$) does precisely mean
\begin{equation}\label{delta2foromega}
\limsup_{t\rightarrow+\infty}\frac{\omega_M(t^2)}{\omega_M(t)}<+\infty,
\end{equation}
which is obviously stronger than \hyperlink{om1}{$(\omega_1)$} for $\omega_M$. Concerning this requirement we formulate the following result.

\begin{theorem}\label{delta2lemma}
Let $M\in\hyperlink{LCset}{\mathcal{LC}}$ be given. Then the following are equivalent:
\begin{itemize}
\item[$(i)$] The associated N-function $F_M$ (see Corollary \ref{princor}) satisfies the $\Delta_2$-condition.

\item[$(ii)$] $\varphi_{\omega_M}$ satisfies the $\Delta_2$-condition.

\item[$(iii)$] $\omega_M$ satisfies \eqref{delta2foromega}.

\item[$(iv)$] $\omega_M$ satisfies
\begin{equation}\label{om7equ}
\exists\;C>0\;\exists\;H>0\;\forall\;t\ge 0:\;\;\;\omega_M(t^2)\le C\omega_M(Ht)+C.
\end{equation}
\item[$(v)$] $M$ satisfies
\begin{equation}\label{om7equforM}
\exists\;k\in\NN_{>0}\;\exists\;A,B\ge 1\;\forall\;j\in\NN:\;\;\;(M_j)^{2k}\le AB^jM_{kj}.
\end{equation}
\end{itemize}
\end{theorem}

\eqref{om7equ} has already appeared (crucially) in different contexts; it is denoted by $(\omega_7)$ in \cite{sectorialextensions1} and in \cite{dissertation}; by $(\omega_8)$ in \cite{compositionpaper} and by $\Xi$ in \cite{scales}.

\demo{Proof}
$(i)\Leftrightarrow(ii)$ This is clear by Corollary \ref{princor} since $F_M\hyperlink{sim}{\sim_{\mathfrak{c}}}\varphi_{\omega_M}$.\vspace{6pt}

$(ii)\Leftrightarrow(iii)$ This is immediate.\vspace{6pt}

$(iii)\Rightarrow(iv)$ This is clear since $\omega_M$ is non-decreasing, $\lim_{t\rightarrow+\infty}\omega_M(t)=+\infty$ and w.l.o.g. $H\ge 1$.\vspace{6pt}

$(iv)\Rightarrow(iii)$ As mentioned at the beginning of \cite[Appendix A]{sectorialextensions1} each non-decreasing function with property \eqref{om7equ} has already \hyperlink{om1}{$(\omega_1)$}. By iterating this property we get $\omega_M(t^2)\le C_1\omega_M(t)+C_1$ for some $C_1\ge 1$ and all $t\ge 0$, i.e. \eqref{delta2foromega}.\vspace{6pt}

$(iv)\Leftrightarrow(v)$ This has been shown in \cite[Lemma 6.3, Cor. 6.4]{scales}.
\qed\enddemo

In \cite[Sect. 2.3, Thm. 3 $(a)$]{orliczbook1} a characterization of $\Delta_2$ in terms of the left-continuous density/left-derivative of $F$ (see Remark \ref{RaoRendefrem}) has been obtained, too.

\begin{example}\label{Delta2rem}
We list some examples and their consequences:
\begin{itemize}
\item[$(*)$] According to \cite[Lemma 6.5]{scales} we know that any $M\in\hyperlink{LCset}{\mathcal{LC}}$ cannot satisfy \hyperlink{mg}{$(\on{mg})$} and \eqref{om7equforM} simultaneously. In particular, the Gevrey sequences $G^s:=(j!^s)_{j\in\NN}$, $s>0$, are violating \eqref{om7equforM}.

\item[$(*)$] Consider the sequences $M^{q,n}:=(q^{j^n})_{j\in\NN}$ with $q,n>1$. Then \eqref{om7equforM} is valid (with $A=B=1$ and $k$ such that $k\ge 2^{1/(n-1)}$) as it is shown in \cite[Example 6.6 $(1)$]{scales}.

\item[$(*)$] Combining Lemma \ref{delta2growthcomp}, Theorem \ref{delta2lemma} and Remark \ref{diffrelrem} yields the following: Let $M,L\in\hyperlink{LCset}{\mathcal{LC}}$ be given and assume that either $M$ or $L$ has \hyperlink{mg}{$(\on{mg})$} and that either $M$ or $L$ has \eqref{om7equforM}. Then
$$\eqref{comparisonprop0equ}\Longleftrightarrow F_M\preceq F_L\Longleftrightarrow F_L\hyperlink{ompreceq}{\preceq_{\mathfrak{c}}}F_M\Longleftrightarrow\eqref{comparisonpropequ},$$
i.e. in \eqref{diffrelremequ} the second implication can also be reversed.
\end{itemize}
\end{example}

\begin{corollary}\label{delta2lemmacor0}
Let $M,L\in\hyperlink{LCset}{\mathcal{LC}}$ be given. Assume that both $F_M$ and $F_L$ satisfy the $\Delta_2$-condition. Then $F_{M\cdot L}$ and $F_{M\star L}$ satisfy the $\Delta_2$-condition, too.
\end{corollary}

\demo{Proof}
By assumption $M$ and $L$ both have \eqref{om7equforM} and so it is immediate that $M\cdot L$ has this property as well.

Concerning the convolution note that we get $\omega_{M\star L}=\omega_M+\omega_L$ and so $\omega_{M\star L}$ has \eqref{om7equ} because both $\omega_M$ and $\omega_L$ have this property.
\qed\enddemo

Finally we apply Theorem \ref{delta2lemma} to $M=M^G$ and get the following.

\begin{corollary}\label{delta2lemmacor}
Let $G$ be an N-function. Let $\omega_G$ be the function from \eqref{Flogweight}, $M^G\in\hyperlink{LCset}{\mathcal{LC}}$ the associated sequence defined via \eqref{assosequ} and $F_{M^G}$ the associated N-function. Then the following are equivalent:
\begin{itemize}
\item[$(i)$] $G$ satisfies the $\Delta_2$-condition.

\item[$(ii)$] The associated N-function $F_{M^G}$ satisfies the $\Delta_2$-condition.

\item[$(iii)$] $\varphi_{\omega_{M^G}}$ satisfies the $\Delta_2$-condition.

\item[$(iv)$] The function $\varphi^G$ (see \eqref{Gvequation}) satisfies the $\Delta_2$-condition.

\item[$(v)$] $\omega_G$ (equivalently $\omega_{M^G}$) satisfies \eqref{delta2foromega}.

\item[$(vi)$] $\omega_G$ (equivalently $\omega_{M^G}$) satisfies \eqref{om7equ}.

\item[$(vii)$] $M^G$ satisfies \eqref{om7equforM}.
\end{itemize}
\end{corollary}

\demo{Proof}
Everything is immediate from Theorem \ref{delta2lemma} applied to $M^G$: $(i)\Leftrightarrow(ii)\Leftrightarrow(iii)\Leftrightarrow(iv)$ holds by Theorem \ref{assoseququotientcomparioncor} and since $\Delta_2$ is preserved under equivalence. Concerning $(v)$ and $(vi)$ note that condition \eqref{delta2foromega} resp. \eqref{om7equ} holds equivalently for $\omega_G$ and $\omega_{M^G}$ (by \eqref{assoweightequiv} and since these conditions are preserved under $\sim$).
\qed\enddemo

\subsection{The $\nabla_2$-condition}\label{nabla2section}
Let $F$ be an N-function and $F^c$ its complementary function. In \cite[Thm. 4.2]{orliczbook} and \cite[Sect. 2.3, Thm. 3 $(a)$]{orliczbook1} it has been shown that $F^c$ satisfies the $\Delta_2$-condition if and only if $F$ has
\begin{equation}\label{Delta2complementary}
\exists\;\ell>1\;\exists\;t_0>0\;\forall\;t\ge t_0:\;\;\;2\ell F(t)\le F(\ell t),
\end{equation}
also known under the name $\nabla_2$ for $F$ (see also \cite[Sect. 2.3, Def. 2, p. 22]{orliczbook1}). A direct characterization for $\nabla^2$ for $F$ in terms of the left-derivative of $F$ has been obtained in \cite[Sect. 2.3, Thm. 3 $(b)$]{orliczbook1}.

\begin{remark}\label{Delta2complementaryrem}
We summarize some properties for $\nabla_2$. For $(i)$ and $(ii)$ it is sufficient to require that $F,G:[0,+\infty)\rightarrow[0,+\infty)$ are non-decreasing, normalized and convex.
\begin{itemize}
\item[$(i)$] $\nabla_2$ is preserved under relation \hyperlink{sim}{$\sim_{\mathfrak{c}}$}; this follows either from $(i)$ in Remark \ref{compequivrem} and since $\Delta_2$ is preserved under \hyperlink{sim}{$\sim_{\mathfrak{c}}$}, or it can be seen directly as follows:

    Assume that $F$ has $\nabla_2$ and let $G$ be another N-function such that $F\hyperlink{sim}{\sim_{\mathfrak{c}}}G$. So
    $$\exists\;k>1\;\exists\;t_0>0\;\forall\;t\ge t_0:\;\;\;G(tk^{-1})\le F(t)\le G(tk).$$
    When iterating $\nabla_2$ we find (since $\ell>1$)
    $$\exists\;t_1>0\;\forall\;n\in\NN_{>0}\;\forall\;t\ge t_1:\;\;\;F(t)\le\frac{1}{(2\ell)^n}F(\ell^nt).$$
    We choose $d\in\NN_{>0}$ such that $\ell^d\ge k^2$ and $n\in\NN_{>0}$ such that $2^{n-1}\ge\ell^d$. Then we estimate for all $t\ge\max\{t_0,t_1\}$ as follows:
\begin{align*}
G(t)\le F(kt)\le\frac{1}{(2\ell)^n}F(\ell^n kt)\le\frac{1}{(2\ell)^n}F(\ell^{n+d}k^{-1}t)\le\frac{1}{(2\ell)^n}G(\ell^{n+d}t)\le\frac{1}{2\ell^{n+d}}G(\ell^{n+d}t),
\end{align*}
i.e. $\nabla_2$ for $G$ holds with $\ell^{n+d}$ and for all $t\ge t_2:=\max\{t_0,t_1\}$.

\item[$(ii)$] Similarly, we show that $\nabla_2$ is preserved under relation $\sim$: Assume that $F$ has $\nabla_2$ and let $G$ be another N-function such that $F\sim G$. So
    $$\exists\;k>1\;\exists\;t_0>0\;\forall\;t\ge t_0:\;\;\;k^{-1}G(t)\le F(t)\le kG(t),$$
    and then we take $n\in\NN_{>0}$ such that $k^2\le 2^{n-1}$. Iterating $n$-times property $\nabla_2$ gives for all $t$ sufficiently large that
    $$\frac{1}{k}(2\ell)^nG(t)\le(2\ell)^nF(t)\le F(\ell^nt)\le kG(\ell^nt).$$
    Since $G$ is an N-function we get $k^2G(\ell^nt)\le G(k^2\ell^nt)$ (recall \eqref{equ114}) and $2k^2\ell^n\le(2\ell)^n$ by the choice of $n$. Consequently, $2k^2\ell^nG(t)\le G(k^2\ell^nt)$ is verified for all sufficiently large $t$, i.e. $\nabla_2$ for $G$ with $k^2\ell^n$.

\item[$(iii)$] Let us prove that for any non-decreasing $F:[0,+\infty)\rightarrow[0,+\infty)$ with $\lim_{t\rightarrow+\infty}F(t)=+\infty$ we have that \eqref{Delta2complementary} is equivalent to
\begin{equation}\label{Delta2complementarywithC}
\exists\;C\ge 1\;\exists\;\ell>1\;\forall\;t\ge 0:\;\;\;2\ell F(t)\le F(\ell t)+2\ell C.
\end{equation}
\eqref{Delta2complementary} implies \eqref{Delta2complementarywithC} with the same $\ell$, e.g. take $C:=F(t_0)$, because $F$ is non-decreasing. For the converse first we iterate \eqref{Delta2complementarywithC} and get $4\ell^2F(t)\le 2\ell F(\ell t)+4\ell^2C\le F(\ell^2t)+2\ell C+4\ell^2C$. Then, since $F(t)\rightarrow+\infty$ as $t\rightarrow+\infty$ we can find $t_0>0$ such that $F(\ell^2t)\ge 2\ell C+4\ell^2C$ for all $t\ge t_0$. Thus $4\ell^2F(t)\le 2F(\ell^2t)$ for all $t\ge t_0$ is verified, i.e. \eqref{Delta2complementary} with the choice $\ell^2$.
\end{itemize}
\end{remark}

In the weight sequence setting we are interested in having \eqref{Delta2complementary} for $F_M$ resp. for $\varphi_{\omega_M}$. Since $\nabla_2$ is preserved under equivalence, via Corollary \ref{princor} this condition transfers into
\begin{equation}\label{Delta2complementary1}
\exists\;\ell>1\;\exists\;s_0>1\;\forall\;s\ge s_0:\;\;\;2\ell\omega_M(s)\le\omega_M(s^{\ell}).
\end{equation}

The aim is to characterize now \eqref{Delta2complementary1} in terms of $M$.

\begin{theorem}\label{nabla2charact}
Let $M\in\hyperlink{LCset}{\mathcal{LC}}$ be given. Then the following are equivalent:
\begin{itemize}
\item[$(i)$] The associated N-function $F_M$ satisfies the $\nabla_2$-condition.

\item[$(ii)$] $\varphi_{\omega_M}$ satisfies the $\nabla_2$-condition.

\item[$(iii)$] $\omega_M$ satisfies \eqref{Delta2complementary1}.

\item[$(iv)$] $\omega_M$ satisfies
\begin{equation}\label{Delta2complementary2}
\exists\;C\ge 1\;\exists\;\ell>1\;\forall\;s\ge 0:\;\;\;2\ell\omega_M(s)\le\omega_M(s^{\ell})+2\ell C.
\end{equation}

\item[$(v)$] The sequence $M$ satisfies
\begin{equation}\label{Delta2complementary3}
\exists\;A\ge 1\;\exists\;\ell>1\;\forall\;j\in\NN:\;\;\;M_{2j}\le AM_j^{2\ell}.
\end{equation}
\end{itemize}
\end{theorem}

The proof shows that in $(iv)$ and $(v)$ we can take the same choice for $\ell$ and the correspondence between $C$ and $A$ is given by $A=e^{2\ell C}$. Consequently, if \eqref{Delta2complementary2} holds for $\ell>1$ then also for all $\ell'\ge\ell$ (even with the same choice for $C$).

\demo{Proof}
$(i)\Leftrightarrow(ii)$ follows from Corollary \ref{princor} and $(i)$ in Remark \ref{Delta2complementaryrem}. $(ii)\Leftrightarrow(iii)$ is clear and $(iii)\Leftrightarrow(iv)$ follows as in resp. by $(iii)$ in Remark \ref{Delta2complementaryrem}.\vspace{6pt}

$(iv)\Rightarrow(v)$ By using \eqref{Prop32Komatsu} we get for all $j\in\NN$:
\begin{align*}
M_{2j}&=\sup_{t\ge 0}\frac{t^{2j}}{\exp(\omega_M(t))}=\sup_{t\ge 0}\frac{t^{2j\ell}}{\exp(\omega_M(t^{\ell}))}\le e^{2\ell C}\sup_{t\ge 0}\frac{t^{2j\ell}}{\exp(2\ell\omega_M(t))}
\\&
=e^{2\ell C}\left(\sup_{t\ge 0}\frac{t^j}{\exp(\omega_M(t))}\right)^{2\ell}=e^{2\ell C}M_j^{2\ell},
\end{align*}
so \eqref{Delta2complementary3} is verified with $A:=e^{2\ell C}$ and the same $\ell$.\vspace{6pt}

$(v)\Rightarrow(iv)$ We have $\frac{t^j}{M_j^{\ell}}\le\sqrt{A}\frac{t^j}{(M_{2j})^{1/2}}$ for all $j\in\NN$ and $t\ge 0$. This yields by definition of associated weight functions
\begin{equation}\label{nabla2charactequ}
\forall\;t\ge 0:\;\;\;\omega_{M^{\ell}}(t)\le\omega_{\widetilde{M}^2}(t)+A_1,
\end{equation}
with $A_1:=\frac{\log(A)}{2}$, $M^{\ell}:=(M_j^{\ell})_{j\in\NN}$ and the auxiliary sequence $\widetilde{M}^2:=(M_{2j}^{1/2})_{j\in\NN}$, see \cite[Sect.3]{modgrowthstrange}, \cite[$(2.5), (2.6)$]{subaddlike} and also \cite[Lemma 6.5]{PTTvsmatrix}. When taking $c=2$ in \cite[Lemma 6.5 $(6.7)$]{PTTvsmatrix}, then
\begin{equation}\label{goodequivalenceclassic}
\exists\;D\ge 1\;\forall\;t\ge 0:\;\;\;\omega_{\widetilde{M}^2}(t)\le 2^{-1}\omega_M(t)\le 2\omega_{\widetilde{M}^2}(t)+D.
\end{equation}
Using \eqref{powersub} and the first half of \eqref{goodequivalenceclassic} we continue \eqref{nabla2charactequ} and get
$$\forall\;t\ge 0:\;\;\;\ell\omega_M(t^{1/\ell})=\omega_{M^{\ell}}(t)\le\omega_{\widetilde{M}^2}(t)+A_1\le 2^{-1}\omega_M(t)+A_1,$$
hence \eqref{Delta2complementary2} is verified with the same $\ell$ and $C:=\frac{A_1}{\ell}=\frac{\log(A)}{2\ell}$.
\qed\enddemo

\begin{example}\label{mgstrongerrem}
We provide some examples of sequences such that \eqref{Delta2complementary3} is valid.
\begin{itemize}
\item[$(i)$] Any $M\in\hyperlink{LCset}{\mathcal{LC}}$ with \hyperlink{mg}{$(\on{mg})$} does also have \eqref{Delta2complementary3}: By \cite[Thm. 9.5.1]{dissertation} or \cite[Thm. 9.5.3]{dissertation} applied to the matrix $\mathcal{M}=\{M\}$ (see also \cite[Thm. 1]{matsumoto}) we know that \hyperlink{mg}{$(\on{mg})$} is equivalent to
\begin{equation}\label{diagmg}
\exists\;B\ge 1\;\forall\;j\in\NN:\;\;\;M_{2j}\le B^jM_j^2.
\end{equation}
Then note that $B^jM_j^2\le AM_j^{2\ell}\Leftrightarrow B\le A^{1/j}M_j^{2(\ell-1)/j}$ holds for all $j\in\NN_{>0}$ if $A\ge 1$ is chosen sufficiently large because $\lim_{j\rightarrow+\infty}(M_j)^{1/j}=+\infty$ and $\ell>1$.

Thus, in particular, the Gevrey sequences $G^s$, $s>0$, have \eqref{Delta2complementary3}.

\item[$(ii)$] However, the converse implication is not valid in general: Consider again the sequences $M^{q,n}:=(q^{j^n})_{j\in\NN}$ with $q,n>1$. Condition \eqref{Delta2complementary3} is valid because for given $n>1$ we choose $\ell\ge 2^{n-1}$ and get for all $j\in\NN$ (and $q>1$) that
    $$q^{(2j)^n}=M^{q,n}_{2j}\le(M^{q,n}_j)^{2\ell}=(q^{j^n})^{2\ell}\Leftrightarrow 1\le q^{2j^n(\ell-2^{n-1})}.$$
    But \eqref{diagmg} is violated: This requirement means $q^{(2j)^n}\le B^jq^{2j^n}$ and since $n>1$ this is impossible for any choice of $B$ as $j\rightarrow+\infty$.
\end{itemize}
\end{example}

\begin{corollary}\label{nabla2charactcor0}
Let $M,L\in\hyperlink{LCset}{\mathcal{LC}}$ be given and assume that both $F_M$ and $F_L$ satisfy the $\nabla_2$-condition. Then $F_{M\cdot L}$ and $F_{M\star L}$ have the $\nabla_2$-condition as well.
\end{corollary}

\demo{Proof}
By assumption both sequences satisfy \eqref{Delta2complementary3} and so it is immediate that $M\cdot L$ has this property, too.

Concerning the convolution we have $\omega_{M\star L}=\omega_M+\omega_L$ and so $\omega_{M\star L}$ has \eqref{Delta2complementary2} because both $\omega_M$ and $\omega_L$ have this property. For this recall that if \eqref{Delta2complementary2} holds for some $\ell>1$ then for all $\ell'\ge\ell$ as well.
\qed\enddemo

Finally, let us apply Theorem \ref{nabla2charact} to $M=M^G$ from Section \ref{fromNtoassofct}.

\begin{corollary}\label{nabla2charactcor}
Let $G$ be an given N-function. Let $\omega_G$ be the associated weight function from \eqref{Flogweight}, $M^G$ the associated weight sequence (see \eqref{assosequ}) and finally $F_{M^G}$ the N-function associated with $M^G$. Then the following are equivalent:
\begin{itemize}
\item[$(i)$] $G$ satisfies the $\nabla_2$-condition (equivalently the complementary N-function $G^c$ satisfies the $\Delta_2$-condition).

\item[$(ii)$] $F_{M^G}$ satisfies the $\nabla_2$-condition (which is equivalent to the fact that the complementary N-function $F^c_{M^G}$ satisfies the $\Delta_2$-condition).

\item[$(iii)$] $\varphi_{\omega_{M^G}}$ satisfies the $\nabla_2$-condition (equivalently the complementary function $\varphi^c_{\omega_{M^G}}$ satisfies the $\Delta_2$-condition).

\item[$(iv)$] The function $\varphi^G$ (see \eqref{Gvequation}) satisfies the $\nabla_2$-condition.

\item[$(v)$] $\omega_G$ (equivalently $\omega_{M^G}$) satisfies \eqref{Delta2complementary1}.

\item[$(vi)$] $\omega_G$ (equivalently $\omega_{M^G}$) satisfies \eqref{Delta2complementary2}.

\item[$(vii)$] $M^G$ satisfies \eqref{Delta2complementary3}.
\end{itemize}
\end{corollary}

\demo{Proof}
We apply Theorem \ref{nabla2charact} to $M^G$. Again $(i)\Leftrightarrow(ii)\Leftrightarrow(iii)\Leftrightarrow(iv)$ follows by Theorem \ref{assoseququotientcomparioncor} and the fact that both the $\Delta_2$ and the $\nabla_2$-condition are preserved under equivalence, recall $(i)$ in Remark \ref{Delta2complementaryrem} for the latter one.

Concerning $(v)$ and $(vi)$ note that by $(ii)$ in Remark \ref{Delta2complementaryrem} both \eqref{Delta2complementary1} and \eqref{Delta2complementary2} are preserved under relation $\sim$ which holds between $\omega_G$ and $\omega_{M^G}$ (recall \eqref{assoweightequiv}).
\qed\enddemo

\subsection{The $\Delta^2$-condition}\label{Deltasquaresection}
According to \cite[Chapter I, \S 6, 5]{orliczbook} and \cite[Def. 5, p. 40]{orliczbook1} we say that an N-function $F$ satisfies the $\Delta^2$-condition if
\begin{equation*}\label{deltasquare}
\exists\;k>1\;\exists\;t_0>0\;\forall\;t\ge t_0:\;\;\;F(t)^2\le F(kt).
\end{equation*}
As stated on \cite[p. 41]{orliczbook}, the $\Delta^2$-condition is preserved under relation \hyperlink{sim}{$\sim_{\mathfrak{c}}$} and $\Delta^2$ holds if and only if $F^{\alpha}\hyperlink{sim}{\sim_{\mathfrak{c}}}F$ for some/any $\alpha>1$. In fact, this equivalence holds for any non-decreasing $F:[0,+\infty)\rightarrow[0,+\infty)$ with $\lim_{t\rightarrow+\infty}F(t)=+\infty$. By \cite[Thm. 6.8]{orliczbook} and \cite[Sect. 2.5, Thm. 8]{orliczbook1} we know that $F$ has $\Delta^2$ if and only if the complementary function $F^c$ has
\begin{equation*}\label{deltasquareforcompl}
\exists\;k>1\;\exists\;t_0>0\;\forall\;t\ge t_0:\;\;\;F^c(t^2)\le ktF^c(t),
\end{equation*}
also known as the $\nabla^2$-condition, see e.g. \cite[Def. 7, p. 41]{orliczbook1}.

In view of Corollary \ref{princor} the associated N-function $F_M$ satisfies the $\Delta^2$-condition if and only if
\begin{equation*}\label{deltasquareforFM}
\exists\;k>1\;\exists\;t_0>0\;\forall\;t\ge t_0:\;\;\;(\omega_M(e^t))^2=\varphi_{\omega_M}(t)^2\le\varphi_{\omega_M}(kt)=\omega_M(e^{kt}),
\end{equation*}
i.e.
\begin{equation}\label{deltasquareforFM1}
\exists\;k>1\;\exists\;s_0>1\;\forall\;s\ge s_0:\;\;\;\omega_M(s)^2\le\omega_M(s^k).
\end{equation}
Note that \eqref{deltasquareforFM1} is not well-related to \eqref{Prop32Komatsu} and in general it seems to be difficult to obtain a characterization for $\Delta^2$ in terms of $M$ by using this formula.

However, we give the following characterization of $\Delta^2$ in the weight sequence setting. First we have to prove a technical result:

\begin{lemma}\label{Delta2lemma}
Let $M\in\hyperlink{LCset}{\mathcal{LC}}$ be given. Then the following are equivalent:
\begin{itemize}
\item[$(i)$] The counting function $\Sigma_M$ satisfies
\begin{equation}\label{Delta2lemmaequ}
\exists\;K>0\;\exists\;t_0>0\;\forall\;t\ge t_0:\;\;\;\Sigma_M(e^t)^2\le\Sigma_M(e^{tK}),
\end{equation}
i.e. \cite[$(6.9)$]{orliczbook} for $p=\Sigma_M\circ\exp$ (this is the $\Delta^2$-condition for $\Sigma_M\circ\exp$).

\item[$(ii)$] The sequence of quotients $\mu$ satisfies
\begin{equation}\label{Delta2lemmaequ1}
\exists\;A>0\;\exists\;j_0\in\NN\;\forall\;j\ge j_0:\;\;\;\mu_{j^2}\le\mu_j^A.
\end{equation}
\end{itemize}
\end{lemma}

The proof shows the correspondence $A=K$.

\demo{Proof}
$(i)\Rightarrow(ii)$ Write $s:=e^t$ and so $\Sigma_M(s)^2\le\Sigma_M(s^K)$ is satisfied for all $s\ge s_0:=e^{t_0}$. Let $j_0\in\NN_{>0}$ be minimal such that $\mu_{j_0}\ge s_0$. Let $j\ge j_0$ be such that $\mu_j<\mu_{j+1}$ and take $s$ with $\mu_j\le s<\mu_{j+1}$. Then $j^2=\Sigma_M(s)^2\le\Sigma_M(s^K)$ follows which implies $\mu_{j^2}\le s^K$. In particular, when taking $s:=\mu_j$ we have shown \eqref{Delta2lemmaequ1} with $A:=K$ and all $j\ge j_0$ such that $\mu_j<\mu_{j+1}$. If $j\ge j_0$ with $\mu_j=\dots=\mu_{j+\ell}<\mu_{j+\ell+1}$ for some $\ell\in\NN_{>0}$, then following the previous step we get $\mu_{(j+\ell)^2}\le\mu_i^K$ for all $j\le i\le j+\ell$. Since $(j+\ell)^2\ge i^2$ for all such indices $i$ and since by log-convexity $j\mapsto\mu_j$ is non-decreasing we are done for all $j\ge j_0$.\vspace{6pt}

$(ii)\Rightarrow(i)$ Let $s\ge\mu_{j_0}$ and so $\mu_j\le s<\mu_{j+1}$ for some $j\ge j_0$. Then $\Sigma_M(s)=j$ and $s^A\ge\mu_j^A\ge\mu_{j^2}$ which implies $\Sigma_M(s^A)\ge j^2=\Sigma_M(s)^2$. Thus \eqref{Delta2lemmaequ} is shown with $K:=A$ and $t_0:=\log(\mu_{j_0})$.
\qed\enddemo

Using this we get the following characterization.

\begin{theorem}\label{assoNfctDeltasquare}
Let $M\in\hyperlink{LCset}{\mathcal{LC}}$ be given. Then the following are equivalent:
\begin{itemize}
\item[$(i)$] The counting function $\Sigma_M$ satisfies \eqref{Delta2lemmaequ}, i.e. $\Sigma_M\circ\exp$ satisfies the $\Delta^2$-condition.

\item[$(ii)$] The sequence of quotients $\mu$ satisfies \eqref{Delta2lemmaequ1}.

\item[$(iii)$] The associated N-function $F_M$ (or equivalently $\varphi_{\omega_M}$) satisfies the $\Delta^2$-condition.
\end{itemize}
\end{theorem}

\demo{Proof}
$(i)\Leftrightarrow(ii)$ is verified in Lemma \ref{Delta2lemma}.

$(i)\Rightarrow(iii)$ In view of \eqref{princorequ} the function $f_M$ appearing in the representation \eqref{Ndef} of $F_M$ enjoys the $\Delta^2$-condition (with the same $K$ and for all $t\ge t_1$, $t_1$ possibly strictly larger than $t_0$ appearing in \eqref{Delta2lemmaequ}). Thus we can apply \cite[Lemma 6.1, Thm. 6.4]{orliczbook} in order to conclude. There the estimate in \eqref{Delta2lemmaequ} is assumed to be strict; however the proof only requires $\le$.\vspace{6pt}

$(iii)\Rightarrow(i)$ We repeat the arguments from \cite[Sect. 2.5, Thm. 8. $(i)\Rightarrow(ii)$]{orliczbook1} (for the representation involving the left-continuous density for $F_M$). First, since $F_M$ satisfies the $\Delta^2$-condition it also has the $\Delta_3$-condition, see Section \ref{Delta3section}. Consequently, by taking into account Lemma \ref{Delta3equiv} and Theorem \ref{Delta3equivforM} we get that $\Sigma_M\circ\exp$ resp. equivalently $f_M$ satisfies the $\Delta_3$-condition, too. Summarizing, we have
\begin{equation}\label{assoNfctDeltasquareequ}
\exists\;k>1\;\exists\;t_0>0\;\forall\;t\ge t_0:\;\;\;tf_M(t)\le f_M(kt),\hspace{30pt}F_M(t)^2\le F_M(kt).
\end{equation}
Using the second part of this and since $F_M(2x)=\int_0^{2x}f_M(s)ds\ge\int_x^{2x}f_M(s)ds\ge f_M(x)x$ for all $x\ge 0$ we estimate as follows for all $x>0$ with $x\ge t_0$:
$$f_M(x)^2\le\frac{1}{x^2}F_M(2x)^2\le\frac{1}{x^2}F_M(2kx)=\frac{1}{x^2}\int_0^{2kx}f_M(s)ds\le\frac{1}{x^2}f_M(2kx)2kx.$$
If $x\ge\max\{1,t_0\}$, then this estimate and the first part in \eqref{assoNfctDeltasquareequ} applied to $t:=2kx(>t_0)$ give
$$f_M(x)^2\le 2kxf_M(2kx)\le f_M(2k^2x),$$
i.e. the $\Delta^2$-condition for $f_M$ is verified with $t_1:=\max\{1,t_0\}$ and $k_1:=2k^2$. By \eqref{princorequ} this is equivalent to the fact that $\Sigma_M\circ\exp$ satisfies the $\Delta^2$-condition; i.e. \eqref{Delta2lemmaequ} for $\Sigma_M$.
\qed\enddemo

\begin{corollary}\label{assoNfctDeltasquarecor}
Let $G$ be an N-function and let $M^G$ be the associated weight sequence (see \eqref{assosequ}). Then the following are equivalent:
\begin{itemize}
\item[$(i)$] The sequence of quotients $\mu^G$ satisfies \eqref{Delta2lemmaequ1}.

\item[$(ii)$] $G$ satisfies the $\Delta^2$-condition.

\item[$(iii)$] The associated N-function $F_{M^G}$ satisfies the $\Delta^2$-condition.

\item[$(iv)$] The function $\varphi^G$ (see \eqref{Gvequation}) satisfies the $\Delta^2$-condition.
\end{itemize}
\end{corollary}

Recall that in order to verify \eqref{Delta2lemmaequ1} for $\mu^G$ for abstractly given N-functions $G$ the formula \eqref{rmucomparison} can be used.

\demo{Proof}
First, Theorem \ref{assoNfctDeltasquare} applied to $M^G$ yields $(i)\Leftrightarrow(iii)$. By Theorem \ref{assoseququotientcomparioncor} we have that $F_{M^G}$, $G$ and $\varphi^G$ are equivalent and since $\Delta^2$ is preserved under equivalence we are done.
\qed\enddemo

\begin{remark}\label{Deltasquarermark}
We gather some observations:
\begin{itemize}
\item[$(i)$] \eqref{Delta2lemmaequ1} means that the sequence of quotients has to increase ``relatively slowly'' and w.l.o.g. we can assume $A\in\NN_{\ge 2}$ in this condition (since $\mu_j\ge 1$ for all $j$).

\item[$(ii)$] \eqref{Delta2lemmaequ1} is preserved under relation \hyperlink{cong}{$\cong$}: Let $M,L\in\hyperlink{LCset}{\mathcal{LC}}$ such that $M\cong L$ and assume that $M$ has \eqref{Delta2lemmaequ1}. Then
$$\exists\;A>0\;\exists\;B\ge 1\;\exists\;j_0\in\NN\;\forall\;j\ge j_0:\;\;\;\frac{1}{B}\lambda_{j^2}\le\mu_{j^2}\le\mu_j^A\le B\lambda_j^A$$
follows and since $\lim_{j\rightarrow+\infty}\lambda_j=+\infty$ we get $B^2\lambda_j^A\le\lambda_j^{A'}$ for any $A'>A$ and for all $j\ge j_{A',B}$ sufficiently large. Thus $L$ satisfies \eqref{Delta2lemmaequ1} when choosing $A'>A$ and restricting to $j\ge\max\{j_0,j_{A',B}\}$.

\item[$(iii)$] In \cite[$(6.10)$, p. 43]{orliczbook} another sufficiency criterion is given. Suppose that an N-function $F$ has
$$\exists\;\alpha>0\;\exists\;t_0>0:\;\;\;t\mapsto\frac{\log(F(t))}{t^{\alpha}}\;\;\;\text{is not decreasing on}\;[t_0,+\infty),$$
then $F$ satisfies the $\Delta^2$-condition. One verifies that this condition yields $F^2\hyperlink{sim}{\sim_{\mathfrak{c}}}F$ and hence this implication holds for any non-decreasing $F:[0,+\infty)\rightarrow[0,+\infty)$ with $\lim_{t\rightarrow+\infty}F(t)=+\infty$.

In the weight sequence setting this expression amounts to the study of $\frac{\log(\varphi_{\omega_M}(t))}{t^{\alpha}}=\frac{\log(\omega_M(e^t))}{t^{\alpha}}=\frac{\log(\omega_M(s))}{\log(s)^{\alpha}}$ for all $s\ge s_0=e^{t_0}$. If there exists $\alpha>0$ such that $s\mapsto\frac{\log(\omega_M(s))}{\log(s)^{\alpha}}$ is non-decreasing for all large $s$, then $\varphi_{\omega_M}$ satisfies the $\Delta^2$-condition and so $F_M$ as well.
\end{itemize}
\end{remark}

\begin{example}\label{Deltasquareexample}
We comment on some examples.
\begin{itemize}
\item[$(*)$] All Gevrey-sequences $G^s$, $s>0$, satisfy \eqref{Delta2lemmaequ1}: This condition amounts to $(j^2)^s\le(j^s)^A$ and so the choices $A:=2$ and $j_0:=0$ are sufficient (for any $s>0$).

\item[$(*)$] If $M,L\in\hyperlink{LCset}{\mathcal{LC}}$ both have \eqref{Delta2lemmaequ1}, then also the product sequence $M\cdot L$: The corresponding sequence of quotients is given by the product $\mu\cdot\lambda$ and so \eqref{Delta2lemmaequ1} follows immediately. The same implication is not clear for the convolution product $M\star L$.

\item[$(*)$] The sequence $M^{q,2}$ does not satisfy this requirement because in this case the corresponding sequence of quotients is given by $(q^{2j-1})_j$. Then \eqref{Delta2lemmaequ1} transfers into $q^{2j^2-1}\le q^{(2j-1)A}$ but which is impossible for any choice $A\ge 1$ if $j\rightarrow+\infty$.

    Alternatively, one can also show that in this case \eqref{deltasquareforFM1} is violated: We have, see the proof of Corollary \ref{nabla2notDelta3cor} for more details and citations, that $\omega_{M^{q,2}}\sim\omega_2$ for all $q>1$, i.e. $\omega_{M^{q,2}}(t)=O(\omega_2(t))$ and $\omega_2(t)=O(\omega_{M^{q,2}}(t))$ as $t\rightarrow+\infty$ for each $q>1$ with $\omega_2(t):=\max\{0,\log(t)^2\}$. Fix $q>1$ and so \eqref{deltasquareforFM1} gives for some $C\ge 1$ and all $s\ge 1$
    $$-1+C^{-1}\log(s)^4\le\omega_{M^{q,2}}(s)^2\le\omega_{M^{q,2}}(s^k)\le C\log(s^k)^2+C=Ck^2\log(s)^2+C,$$
    yielding a contradiction as $s\rightarrow+\infty$.
\end{itemize}
\end{example}

Summarizing, by Examples \ref{Delta2rem}, \ref{mgstrongerrem} and \ref{Deltasquareexample} we get the following consequences which should be compared with the first diagram on \cite[p. 43]{orliczbook1}:
\begin{equation}\label{counterexamplesequ}
\Delta_2\wedge \nabla_2,\Delta_2,\nabla_2\nRightarrow\Delta^2,\hspace{25pt}\Delta^2\nRightarrow\Delta_2,\hspace{25pt}\nabla_2\nRightarrow\Delta_2.
\end{equation}

Let us study how condition \eqref{Delta2lemmaequ1} is related to moderate growth.

\begin{lemma}\label{sufficienDeltasquarelemma}
Let $M\in\hyperlink{LCset}{\mathcal{LC}}$ be given such that \hyperlink{mg}{$(\on{mg})$} holds and
\begin{equation}\label{sufficienDeltasquarelemmaequ}
\liminf_{n\rightarrow+\infty}(\mu_{2^n})^{1/(n+2)}>1.
\end{equation}
Then $M$ satisfies \eqref{Delta2lemmaequ1} and hence the associated N-function $F_M$ (or equivalently $\varphi_{\omega_M}$) satisfies the $\Delta^2$-condition.
\end{lemma}

\demo{Proof}
First, see e.g. \cite[Lemma 2.2]{whitneyextensionweightmatrix} and the citations there, \hyperlink{mg}{$(\on{mg})$} is equivalent to $\sup_{j\in\NN}\frac{\mu_{2j}}{\mu_j}<+\infty$. Consequently, we find some $A>1$ such that $\mu_{2^kj}\le A^k\mu_j$ for each $k,j\in\NN$.

Take now $j\in\NN_{>0}$ (the case $j=0$ is trivial), then $2^n\le j<2^{n+1}$ for some $n\in\NN$ and $2^{2n}\le j^2<2^{2n+2}$, so
$$\mu_{j^2}\le\mu_{2^{2n+2}}=\mu_{2^{n+2}2^n}\le A^{n+2}\mu_{2^n}.$$
By \eqref{sufficienDeltasquarelemmaequ} we find $n_0\in\NN$ and $\delta>1$ such that $(\mu_{2^n})^{1/(n+2)}\ge\delta$ for all $n\ge n_0$. Set $B:=\frac{\log(A)}{\log(\delta)}+1>1$, so $\delta=A^{1/(B-1)}$ and hence $A^{n+2}\le\mu_{2^n}^{B-1}$ for all $n\ge n_0$. This implies $\mu_{j^2}\le A^{n+2}\mu_{2^n}\le\mu_{2^n}^B\le\mu_j^B$ and so \eqref{Delta2lemmaequ1} is verified (for $j_0:=2^{n_0}$ and choosing $B$ as before).
\qed\enddemo

We finish by commenting on requirement \eqref{sufficienDeltasquarelemmaequ}:

\begin{itemize}
\item[$(*)$] \eqref{sufficienDeltasquarelemmaequ} is a mild extra growth assumption: It follows when $\liminf_{j\rightarrow\infty}\frac{\mu_j}{j}>0$ because then $\mu_j\ge j\epsilon$ for some $\epsilon>0$ and all $j\in\NN$. Thus $(\mu_{2^n})^{1/(n+2)}\ge 2^{n/(n+2)}\epsilon^{1/(n+2)}$ for all $n\in\NN$ and so $(\mu_{2^n})^{1/(n+2)}>1$ for all $n$ sufficiently large (depending on $\epsilon$).

\item[$(*)$] On the other hand note that $\lim_{n\rightarrow+\infty}(2^{ns})^{1/(n+2)}=2^s>1$ for any $s>0$ and hence each $G^s$ satisfies \eqref{sufficienDeltasquarelemmaequ}. However, $\lim_{j\rightarrow+\infty}\frac{j^s}{j}=0$ holds for all $0<s<1$.
\end{itemize}

\subsection{The $\Delta_3$-condition}\label{Delta3section}
An N-function $F$ satisfies the $\Delta_3$-condition, see \cite[Chapter I, \S 6, $(6.1)$]{orliczbook} and \cite[Sect. 2.5, Def. 1, p. 37]{orliczbook1}, if
\begin{equation}\label{Delta3}
\exists\;k>1\;\exists\;t_0>0\;\forall\;t\ge t_0:\;\;\;tF(t)\le F(kt).
\end{equation}
Since $F(t)\ge t$ for all large $t$ (recall the second part in \eqref{equ115116}) we immediately have that $\Delta^2$ implies $\Delta_3$; however the converse is not true in general, see \cite[p. 41]{orliczbook} and \cite[p. 40]{orliczbook1}. It is also known that $\Delta_3$ for $F$ implies $\Delta_2$ for $F^c$, i.e. $F$ has $\nabla_2$, see \cite[Thm. 6.5]{orliczbook}. Moreover, by \cite[Sect. 2.5, Thm. 3]{orliczbook1} we have that $F$ has $\Delta_3$ if and only if $F^c$ satisfies the $\nabla_3$-condition, see  again \cite[Sect. 2.5, Def. 1, p. 37]{orliczbook1}.

$\Delta_3$ is preserved under equivalence and since $tF(t)\ge F(t)$ for all $t\ge 1$ we get that

\begin{itemize}
\item[$(*)$] $F$ has $\Delta_3$ if and only if

\item[$(*)$] $F$ and $t\mapsto tF(t)$ are equivalent,
\end{itemize}

see \cite[Chapter I, \S 6, 1, p. 35]{orliczbook}. Moreover, we have the following reformulation for $\Delta_3$:

\begin{lemma}\label{Delta3equiv}
Let $F$ be an N-function. Then the following are equivalent:
\begin{itemize}
\item[$(i)$] $F$ satisfies the $\Delta_3$-condition.

\item[$(ii)$] We have $\widetilde{F}\hyperlink{sim}{\sim_{\mathfrak{c}}}F$ with
\begin{equation}\label{Delta3equivequ}
\widetilde{F}(t):=\int_0^{|t|}F(s)ds,\;\;\;t\in\RR.
\end{equation}

\item[$(iii)$] The function $f$ from the representation \eqref{Ndef} satisfies the $\Delta_3$-condition.

\item[$(iv)$] We have that $F\hyperlink{sim}{\sim_{\mathfrak{c}}}f$.
\end{itemize}
Consequently, if any of these equivalent conditions holds, then $\widetilde{F}$ has $\Delta_3$, too.
\end{lemma}

\demo{Proof}
$(i)\Rightarrow(ii)$ This is contained in the proof of \cite[Thm. 6.1]{orliczbook}. First, for any N-function $F$ we get for all $t\ge 1$ that
$$\widetilde{F}(2t)=\int_0^{2t}F(s)ds\ge\int_t^{2t}F(s)ds\ge F(t)t\ge F(t).$$
In fact this holds for any non-decreasing and non-negative $F$. On the other hand, by using $\Delta_3$ for some $k>1$ and all $t$ sufficiently large one has
$$\widetilde{F}(t)=\int_0^tF(s)ds\le F(t)t\le F(kt).$$

$(ii)\Rightarrow(i)$ By the equivalence we get $\widetilde{F}(t)\le F(kt)$ for some $k>1$ and all $t$ sufficiently large. Thus for all such large $t$ we estimate by
$$F(k2t)\ge\widetilde{F}(2t)=\int_0^{2t}F(s)ds=\int_0^tF(s)ds+\int_t^{2t}F(s)ds\ge\widetilde{F}(t)+F(t)t\ge F(t)t,$$
i.e. $\Delta_3$ with $k':=2k$.\vspace{6pt}

$(i)\Rightarrow(iii)$ This is shown in the proof of \cite[Sect. 2.5, Thm. 3. $(i)\Rightarrow(ii)$]{orliczbook1} (for the representation involving the left-continuous density for $F$); we repeat the details: First, we estimate by
$$xf(x)\le\int_x^{2x}f(s)ds\le\int_0^{2x}f(s)ds=F(2x)\le\frac{1}{2x}F(2kx)=\frac{1}{2x}\int_0^{2kx}f(s)ds\le\frac{1}{2x}2kxf(2kx),$$
thus $xf(x)\le 2kf(2kx)$ for all $x\ge t_0/2$, with $k>1$, $t_0>0$ from \eqref{Delta3}. We use this estimate and apply it also to $y:=2kx(>t_0)$ in order to get
$$xf(x)\le2kf(2kx)\le xf(2kx)\le f(4k^2x)$$
for all $x\ge t_1:=\max\{t_0,2k\}$. Thus $\Delta_3$ for $f$ is verified with $t_1$ and $k_1:=4k^2$.\vspace{6pt}

$(iii)\Rightarrow(i),(iv)$ We replace in $(i)\Rightarrow(ii)$ the function $F$ by $f$ and $\widetilde{F}$ by $F$ (i.e. \eqref{Delta3equivequ} turns into \eqref{Ndef}).

$(iv)\Rightarrow(iii)$ This follows by replacing in $(ii)\Rightarrow(i)$ the function $F$ by $f$ and $\widetilde{F}$ by $F$.
\qed\enddemo

We apply this characterization to the weight sequence setting.

\begin{theorem}\label{Delta3equivforM}
Let $M\in\hyperlink{LCset}{\mathcal{LC}}$ be given. Then the following are equivalent:
\begin{itemize}
\item[$(i)$] The associated N-function $F_M$ (equivalently $\varphi_{\omega_M}$) satisfies the $\Delta_3$-condition.

\item[$(ii)$] The function $\Sigma_M\circ\exp$ (equivalently $f_M$) satisfies the $\Delta_3$-condition.

\item[$(iii)$] $\widetilde{F}_M\hyperlink{sim}{\sim_{\mathfrak{c}}}F_M$ is valid (with $\widetilde{F}_M$ given by \eqref{Delta3equivequ}).

\item[$(iv)$] $\widetilde{\varphi}_{\omega_M}\hyperlink{sim}{\sim_{\mathfrak{c}}}\varphi_{\omega_M}$ is valid with
\begin{equation*}\label{Delta3equivforMequ}
\widetilde{\varphi}_{\omega_M}(t):=\int_0^{|t|}\varphi_{\omega_M}(s)ds,\;\;\;t\in\RR.
\end{equation*}

\item[$(v)$] We have that
\begin{equation}\label{Delta3equivforMequ1}
\exists\;k>1\;\exists\;s_0>1\;\forall\;s\ge s_0:\;\;\;\omega_M(s)\le\widetilde{\omega}_M(s^2)\le\omega_M(s^k),
\end{equation}
with
$$\widetilde{\omega}_M(t):=\widetilde{\varphi}_{\omega_M}(\log(t)),\;\;\;t\ge 1,\hspace{15pt}\widetilde{\omega}_M(t):=0,\;\;\;0\le t<1.$$
\end{itemize}
\end{theorem}

The function $\widetilde{\omega}_M$ admits the representation
\begin{equation}\label{Delta3equivforMequ2}
\widetilde{\omega}_M(s)=\int_0^{|s|}\frac{\omega_M(u)}{u}du=\int_1^{|s|}\frac{\omega_M(u)}{u}du,\;\;\;s\in\RR.
\end{equation}

\demo{Proof}
$(i)\Leftrightarrow(ii)\Leftrightarrow(iii)$ follows by Lemma \ref{Delta3equiv} applied to $F_M$, by \eqref{princorequ} and the fact that $\Delta_3$ is preserved under equivalence, see Corollary \ref{princor}.\vspace{6pt}

$(iii)\Rightarrow(iv)$ The estimate $\widetilde{\varphi}_{\omega_M}(2t)\ge\varphi_{\omega_M}(t)$ for $t\ge 1$ holds as in $(i)\Rightarrow(ii)$ in Lemma \ref{Delta3equiv} and for the converse we estimate as follows for $t$ sufficiently large:
\begin{align*}
\widetilde{\varphi}_{\omega_M}(t)&=\int_0^{t}\varphi_{\omega_M}(s)ds\le\int_0^t(F_M(s)+C)ds=\widetilde{F}_M(t)+Ct\le F_M(kt)+Ct
\\&
\le\varphi_{\omega_M}(kt)+Ct+D\le 2\varphi_{\omega_M}(kt)\le\varphi_{\omega_M}(2kt).
\end{align*}
The first estimate follows for some $C\ge 1$ (and all $s\ge 0$) by \eqref{equivwithCforM}, the second one since $\widetilde{F}_M\hyperlink{sim}{\sim_{\mathfrak{c}}}F_M$ by assumption, the third one again by \eqref{equivwithCforM}, the fourth since $\lim_{t\rightarrow+\infty}\frac{\varphi_{\omega_M}(t)}{t}=+\infty$, and finally the last one by convexity and normalization for $\varphi_{\omega_M}$ (see \eqref{equ114}).\vspace{6pt}

$(iv)\Rightarrow(iii)$ $\widetilde{F}_M(2t)\ge F_M(t)$ for all $t\ge 1$ is shown in $(i)\Rightarrow(ii)$ in Lemma \ref{Delta3equiv}. Conversely, by assumption $\widetilde{\varphi}_{\omega_M}(t)\le\varphi_{\omega_M}(kt)$ for some $k>1$ and all $t(\ge 1)$ large. Thus for all $t$ sufficiently large:
\begin{align*}
\widetilde{F}_M(t)&=\int_0^{t}F_M(s)ds\le\int_0^t(\varphi_{\omega_M}(s)+D)ds=\widetilde{\varphi}_{\omega_M}(t)+Dt\le\varphi_{\omega_M}(kt)+Dt
\\&
\le F_M(kt)+Dt+C\le 2F_M(kt)\le F_M(2kt).
\end{align*}
The first estimate follows by \eqref{equivwithCforM} (for all $s\ge 0$), the second one by assumption, the third one again by \eqref{equivwithCforM}, for the fourth estimate we have used the second part in \eqref{equ115116} and finally the last one holds by \eqref{equ114}.\vspace{6pt}

$(iv)\Leftrightarrow(v)$ First, for all $t\ge 0$ we have
$$\widetilde{\omega}_M(e^t)=\widetilde{\varphi}_{\omega_M}(t)=\int_0^t\omega_M(e^s)ds=\int_1^{e^t}\frac{\omega_M(u)}{u}du=\int_0^{e^t}\frac{\omega_M(u)}{u}du,$$
because $\omega_M(t)=0$ for $0\le t\le 1$. Thus $\widetilde{\omega}_M(t)=\int_0^t\frac{\omega_M(u)}{u}du$ for all $t\ge 1$ and in fact even for all $t\ge 0$ (since $\widetilde{\omega}_M(t):=0$ for $0\le t\le 1$). Thus \eqref{Delta3equivforMequ2} is verified.

Moreover, $\widetilde{\varphi}_{\omega_M}\hyperlink{sim}{\sim_{\mathfrak{c}}}\varphi_{\omega_M}$ holds if and only if
\begin{equation}\label{Delta3equivforMequ1proof}
\exists\;k>1\;\exists\;t_0>0\;\forall\;t\ge t_0:\;\;\;\varphi_{\omega_M}(t)\le\widetilde{\varphi}_{\omega_M}(2t)\le\varphi_{\omega_M}(kt),
\end{equation}
since the first estimate holds for all $t\ge 1$ as mentioned in $(ii)\Rightarrow(iii)$. \eqref{Delta3equivforMequ1proof} is obviously equivalent to \eqref{Delta3equivforMequ1}.
\qed\enddemo

Using this characterization we give two applications.

\begin{corollary}\label{nabla2notDelta3cor0}
Let $M,L\in\hyperlink{LCset}{\mathcal{LC}}$ be given. If both $F_M$ and $F_L$ have the $\Delta_3$-condition, then $F_{M\star L}$, too. \end{corollary}

\demo{Proof}
Recall that $M\star L\in\hyperlink{LCset}{\mathcal{LC}}$ and $\omega_{M\star L}=\omega_M+\omega_L$. Then \eqref{Delta3equivforMequ2} implies $\widetilde{\omega}_{M\star L}=\widetilde{\omega}_M+\widetilde{\omega}_L$ and by assumption we have \eqref{Delta3equivforMequ1} for both $\omega_M$ and $\omega_L$ and so for $\omega_{M\star L}$ as well. Theorem \ref{Delta3equivforM} yields that $F_{M\star L}$ satisfies the $\Delta_3$-condition, too.
\qed\enddemo

\begin{corollary}\label{nabla2notDelta3cor}
There exist N-functions having $\Delta_2$ and $\nabla_2$ but not $\Delta_3$.
\end{corollary}

This statement should be compared with the first diagram on \cite[p. 43]{orliczbook1}.

\demo{Proof}
Let us consider the sequence(s) $M^{q,n}$ with $q,n>1$. As seen in Examples \ref{Delta2rem}, \ref{mgstrongerrem} each sequence yields an associated N-function $F_{M^{q,n}}$ having both $\Delta_2$ and $\nabla_2$.

We prove now that \eqref{Delta3equivforMequ1} is violated. For this first recall results from \cite[Sect. 5.5]{whitneyextensionweightmatrix} and \cite[Sect. 3.10]{dissertation}: Let $n>1$ be arbitrary but from now on fixed, then each $M^{q,n}$ is an element of the weight matrix (i.e. the one-parameter family of weight sequences) associated with the weight $\omega_s(t):=\max\{0,\log(t)^s\}$, $s>1$, such that $\frac{1}{s}+\frac{1}{n}=1$. Therefore, $\omega_{M^{q,n}}\sim\omega_s$ for all $q>1$, i.e. $\omega_{M^{q,n}}(t)=O(\omega_s(t))$ and $\omega_s(t)=O(\omega_{M^{q,n}}(t))$ as $t\rightarrow+\infty$ for each $q>1$, see \cite[Thm. 4.0.3]{dissertation} and \cite[Lemma 5.7]{compositionpaper}.

For all $x\ge 1$ we have
$$\int_0^x\frac{\omega_s(t)}{t}dt=\int_1^x\frac{\log(t)^s}{t}dt=\left[\frac{\log(t)^{s+1}}{s+1}\right]_{t=1}^{t=x}=\frac{\log(x)^{s+1}}{s+1},$$
and so by the above (recall also \eqref{Delta3equivforMequ2})
$$\forall\;q>1\;\exists\;D\ge 1\;\forall\;x\ge 1:\;\;\; D^{-1}\frac{\log(x)^{s+1}}{s+1}-\log(x)\le\widetilde{\omega}_{M^{q,n}}(x)\le D\frac{\log(x)^{s+1}}{s+1}+D\log(x).$$
Fix now $q>1$ and then, when \eqref{Delta3equivforMequ1} is valid, we obtain
\begin{align*}
&D^{-1}\frac{\log(x^2)^{s+1}}{s+1}-\log(x^2)\le\widetilde{\omega}_{M^{q,n}}(x^2)\le\omega_{M^{q,n}}(x^k)\le D\log(x^k)^s+D
\\&
\Longrightarrow 2^{s+1}\log(x)^{s+1}\le D^2(s+1)k^s\log(x)^s+D^2(s+1)+2(s+1)D\log(x).
\end{align*}
But this is impossible as $x\rightarrow+\infty$ for any choices of $D$ and $k$.
\qed\enddemo

By applying Theorem \ref{Delta3equivforM} to the associated sequence $M^G$, Proposition \ref{assoseququotientcomparion}, Theorem \ref{assoseququotientcomparioncor} and recalling that $\Delta_3$ is preserved under equivalence we obtain:

\begin{corollary}\label{nabla2notDelta3cor1}
Let $G$ be an N-function and let $M^G$ be the associated sequence (see \eqref{assosequ}). Then the following are equivalent:
\begin{itemize}
\item[$(i)$] $G$ satisfies the $\Delta_3$-condition.

\item[$(ii)$] The associated N-function $F_{M^G}$ satisfies the $\Delta_3$-condition.

\item[$(iii)$] $\varphi^G$ (see \eqref{Gvequation}) satisfies the $\Delta_3$-condition.

\item[$(iv)$] The function $\Sigma_{M^G}\circ\exp$ (equivalently $f_{M^G}$) satisfies the $\Delta_3$-condition.

\item[$(v)$] The function $\Sigma^G\circ\exp$ (see \eqref{Fcounting}) satisfies the $\Delta_3$-condition.

\item[$(vi)$] Assertions $(iv)$ resp. $(v)$ listed in Theorem \ref{Delta3equivforM} are valid for the associated sequence $M^G$.
\end{itemize}
\end{corollary}

We close with the following observation concerning assertion $(ii)$ in Theorem \ref{Delta3equivforM}.

\begin{lemma}
Let $M\in\hyperlink{LCset}{\mathcal{LC}}$ be given. Consider the following conditions:
\begin{itemize}
\item[$(i)$] The sequence of quotients $\mu$ satisfies
$$\exists\;k>1\;\exists\;j_0\in\NN\;\forall\;j\ge j_0:\;\;\;\mu_j^k\ge\mu_{\lceil j\log(\mu_{j+1})\rceil}.$$
\item[$(ii)$] The function $\Sigma_M\circ\exp$ (equivalently $f_M$) satisfies the $\Delta_3$-condition.

\item[$(iii)$] The sequence of quotients $\mu$ satisfies
$$\exists\;k>1\;\exists\;j_0\in\NN\;\forall\;j\ge j_0:\;\;\;\mu_j^k\ge\mu_{\lceil j\log(\mu_{j})\rceil}.$$
\end{itemize}
Then $(i)\Rightarrow(ii)\Rightarrow(iii)$ holds.
\end{lemma}

\demo{Proof}
$(i)\Rightarrow(ii)$ Let $t\ge 0$ be such that $\mu_j\le t<\mu_{j+1}$ for some $j\ge j_0$. Then $\Sigma_M(t)=j$ and $t^k\ge\mu_j^k\ge\mu_{\lceil j\log(\mu_{j+1})\rceil}$, hence $\Sigma_M(t^k)\ge\lceil j\log(\mu_{j+1})\rceil\ge j\log(\mu_{j+1})=\Sigma_M(t)\log(\mu_{j+1})\ge\Sigma_M(t)\log(t)$. Hence the $\Delta_3$-condition is verified with the same $k$ and $t_0:=\mu_{j_0}$.\vspace{6pt}

$(ii)\Rightarrow(iii)$ Let $j\in\NN$ such that $\mu_{j+1}>\mu_j\ge t_0$ with $t_0$ the value appearing in the $\Delta_3$-condition. Then this property evaluated at $t=\mu_j$ yields $\log(\mu_j)j=\log(\mu_j)\Sigma_M(\mu_j)\le\Sigma_M(\mu_j^k)$ which gives the desired estimate with the same $k$ since $\Sigma_M(\mu_j^k)\in\NN$. If $\mu_j\ge t_0$ such that $\mu_j=\dots=\mu_{j+d}<\mu_{j+d+1}$, then $t=\mu_j$ yields $\log(\mu_{j+i})(j+i)\le\log(\mu_{j+i})(j+d)=\log(t)\Sigma_M(t)\le\Sigma_M(t^k)$ for all $0\le i\le d$ and so, since $\mu_j^k=\dots=\mu_{j+d}^k$ and $\mu$ is non-decreasing, we get $\mu_{j+i}^k\ge\mu_{\lceil\log(\mu_{j+i})(j+i)\rceil}$ for all $0\le i\le d$. Thus we are done by taking the same $k$ and $j_0\in\NN$ minimal such that $\mu_{j_0}\ge t_0$.
\qed\enddemo

\subsection{The $\Delta'$-condition}\label{Deltaprimesection}
According to \cite[Chapter I, \S 5]{orliczbook} and \cite[Sect. 2.3, Def. 7, p. 28]{orliczbook1} we say that an N-function $F$ satisfies the $\Delta'$-condition if
\begin{equation*}\label{Deltaprime}
\exists\;k>0\;\exists\;u_0>0\;\forall\;t,s\ge u_0:\;\;\;F(ts)\le kF(t)F(s).
\end{equation*}
In the weight sequence setting in view of Corollary \ref{princor} and since $\Delta'$ is preserved under equivalence, see \cite[Chapter I, \S 5, p. 30]{orliczbook}, this condition means that
\begin{equation}\label{Deltaprimeweighsequ}
\exists\;k>0\;\exists\;u_0>0\;\forall\;t,s\ge u_0:\;\;\;\omega_M(t^{\log(s)})\le k\omega_M(s)\omega_M(t).
\end{equation}
By \cite[Lemma 5.1]{orliczbook} and \cite[Sect. 2.3, Lemma 8]{orliczbook1} we know that $\Delta'$ implies $\Delta_2$ and on \cite[p. 30-31]{orliczbook} (see also \cite[p. 29]{orliczbook1}) it is shown that in general this implication is strict. Moreover, by \cite[Thm. 6.6]{orliczbook} it follows that if $F$ satisfies $\Delta^2$, then $F^c$ has $\Delta'$ and by \cite[Sect. 2.3, Thm. 11]{orliczbook1} we know that $F$ has $\Delta'$ if and only if $F^c$ satisfies the so-called $\nabla'$-condition; see \cite[Sect. 2.3, Def. 7, p. 28 $(11)$]{orliczbook1}. \cite[Sect. 2.3, Prop. 12]{orliczbook1} tells us that an N-function $F\in\Delta'\cap\nabla'$ if and only if $F$ is equivalent to $t\mapsto|t|^s$ for some $s>1$; i.e. this corresponds to the weight $\omega_s(t):=\max\{0,\log(t)^s\}$ and hence to the sequences $M^{q,n}$ such that $\frac{1}{s}+\frac{1}{n}=1$, see the proof of Corollary \ref{nabla2notDelta3cor}.\vspace{6pt}

A direct check of \eqref{Deltaprimeweighsequ} seems to be quite technical resp. hardly possible since due to the multiplicative nature this estimate is not well-related w.r.t. formula \eqref{Prop32Komatsu}. The same comment applies to $\nabla'$ and to the counting function $\Sigma_M$.

In \cite[Thm. 5.1]{orliczbook} a sufficiency criterion for $\Delta'$ is shown:

\begin{theorem}\label{Thm51orliczbook}
Let $F(t)=\int_0^{|t|}f(s)ds$ be a given N-function (recall the representation \eqref{Ndef}). Then $F$ satisfies the $\Delta'$-condition provided that $f$ has the following growth property which we abbreviate by $(\Delta'_f)$ from now on:

There exists some $t_0>1$ such that for every fixed $t\ge t_0$ the function $h_f$ given by $h_f(s):=\frac{f(st)}{f(s)}$ is not increasing on $[t_0,+\infty)$.
\end{theorem}

In the weight sequence setting this result takes the following form:

\begin{corollary}\label{Thm51orliczbookcor}
Let $M\in\hyperlink{LCset}{\mathcal{LC}}$ be given. If $\Sigma_M\circ\exp$ satisfies $(\Delta'_{\Sigma_M\circ\exp})$ then the associated N-function $F_M$ (resp. equivalently $\varphi_{\omega_M}$) satisfies the $\Delta'$-condition.
\end{corollary}

\demo{Proof}
In view of \eqref{princorequ}, we get $(\Delta'_{\Sigma_M\circ\exp})$ if and only if $(\Delta'_{f_M})$ when enlarging $t_0$ sufficiently if necessary. Then Theorem \ref{Thm51orliczbook} applied to $f=f_M$ and $F=F_M$ yields the conclusion.
\qed\enddemo

However, we show now that in general $(\Delta'_{\Sigma_M\circ\exp})$ fails in the weight sequence setting.

\begin{proposition}\label{Deltaprimefails}
Let $M\in\hyperlink{LCset}{\mathcal{LC}}$ be given such that $1\le\mu_1<\mu_2<\dots$, i.e. the sequence of quotients is strictly increasing (see Remark \ref{strictincrem}). Then $(\Delta'_{\Sigma_M\circ\exp})$ (resp. equivalently $(\Delta'_{f_M})$) is violated.
\end{proposition}

\demo{Proof}
First, with $u:=e^s$ we get $\frac{\Sigma_M(e^{ts})}{\Sigma_M(e^s)}=\frac{\Sigma_M(u^t)}{\Sigma_M(u)}$ and hence $(\Delta'_{\Sigma_M\circ\exp})$ precisely means:
$$\exists\;t_0>1\;\forall\;t\ge t_0:\;\;\;u\mapsto\frac{\Sigma_M(u^t)}{\Sigma_M(u)}\;\text{is not increasing on}\;\;\;[e^{t_0},+\infty)=:[u_0,+\infty).$$
Let now $t\ge t_0>1$ be arbitrary but fixed. Then $\mu_{j_0}\le u_0<\mu_{j_0+1}$ and $\mu_k\le u_0^t<\mu_{k+1}$ for some (large) $j_0,k\in\NN_{>0}$. Note that $k$ is depending on $t$ and $\frac{\Sigma_M(u_0^t)}{\Sigma_M(u_0)}=\frac{k}{j_0}$. Since $u_0^t\ge u_0$ we clearly have $k\ge j_0$.

We show that assuming $(\Delta'_{\Sigma_M\circ\exp})$ yields a contradiction. Let $u\in[u_0,+\infty)$ increase and we split the argument in several steps:

\begin{itemize}
\item[$(*)$] If $u\ge u_0$ is given with $\mu_{j_0}<u<\mu_{j_0+1}$ and $\mu_k<u^t<\mu_{k+1}$, then for all $u'\in[\epsilon-u,u+\epsilon]$ with $\epsilon>0$ sufficiently small the quotient appearing in $(\Delta'_{\Sigma_M\circ\exp})$ remains constant, i.e. in this case the crucial expression is locally constant.

\item[$(*)$] Even $k>j_0$ is valid: If $k=j_0$, then in order to have $(\Delta'_{\Sigma_M\circ\exp})$ we need $\mu_{j_0}\le u<u^t<\mu_{j_0+1}$ for all $u\ge u_0$ with $\mu_{j_0}<u<\mu_{j_0+1}$, a contradiction as $u\rightarrow\mu_{j_0+1}$.

\item[$(*)$] If $(\Delta'_{\Sigma_M\circ\exp})$ holds, then for all $u$ with $\mu_{j_0}\le u_0\le u<\mu_{j_0+1}$ it follows that $u^t<\mu_{k+1}$: Otherwise, if $u^t\ge\mu_{k+1}$, then $\frac{\Sigma_M(u^t)}{\Sigma_M(u)}\ge\frac{k+1}{j_0}>\frac{k}{j_0}=\frac{\Sigma_M(u_0^t)}{\Sigma_M(u_0)}$, a contradiction.

\item[$(*)$] On the other hand, for all $u\ge u_0$ such that $\mu_k\le u^t<\mu_{k+1}$ it is allowed that $u\ge\mu_{j_0+d}$ for some $d\in\NN_{>0}$: In this case, if $\mu_{j_0+d}\le u<\mu_{j_0+d+1}$ and still $\mu_k\le u^t<\mu_{k+1}$, then $\frac{\Sigma_M(u^t)}{\Sigma_M(u)}=\frac{k}{j_0+d}<\frac{k}{j_0}$.

\item[$(*)$] Take $u=(\mu_{k+1})^{1/t}$ and so $u>u_0$. We get that $\mu_{j_0+d}\le u<\mu_{j_0+d+1}$ for some $d\in\NN$. Thus $\Sigma_M(u^t)=k+1$ holds (since $\mu$ is strictly increasing!) and $\Sigma_M(u)=j_0+d$.

    We distinguish: If $\mu_{j_0+d}<u<\mu_{j_0+d+1}$, then we can find $\epsilon>0$ sufficiently small to ensure $u-\epsilon>\mu_{j_0+d}$ and $\mu_k\le(u-\epsilon)^t<\mu_{k+1}$. In this case $\frac{\Sigma_M(u^t)}{\Sigma_M(u)}=\frac{k+1}{j_0+d}>\frac{k}{j_0+d}=\frac{\Sigma_M((u-\epsilon)^t)}{\Sigma_M(u-\epsilon)}$, a contradiction to $(\Delta'_{\Sigma_M\circ\exp})$.

    In particular this case happens if $d=0$ because then $\mu_{j_0}\le u_0<u$ by assumption.

\item[$(*)$] If now $\mu_{j_0+d}=u<\mu_{j_0+d+1}$ for some $d\ge 1$ and $u^t=\mu_{k+1}$, then take $\epsilon>0$ sufficiently small to ensure $u-\epsilon>\mu_{j_0+d-1}$ and $\mu_k\le(u-\epsilon)^t<\mu_{k+1}$. Both estimates are possible since the sequence $\mu$ is assumed to be strictly increasing.

    Thus $\frac{\Sigma_M(u^t)}{\Sigma_M(u)}=\frac{k+1}{j_0+d}\le\frac{k}{j_0+d-1}=\frac{\Sigma_M((u-\epsilon)^t)}{\Sigma_M(u-\epsilon)}$ and which verifies $(\Delta'_{\Sigma_M\circ\exp})$ in this case because this estimate is equivalent to $j_0+d-1\le k$ and this is clear since $\mu_{j_0+d}=u<u^t=\mu_{k+1}$.

\item[$(*)$] Summarizing all the information, a necessary condition to ensure $(\Delta'_{\Sigma_M\circ\exp})$ is that
\begin{equation}\label{starviolating}
\exists\;j_0\in\NN_{>0}\;\exists\;t_0>1\;\forall\;t\ge t_0\;\exists\;k\in\NN_{>0},\;k>j_0,\;\exists\;d\in\NN_{>0}:\;\;\;\mu_{k+1}=(\mu_{j_0+d})^t.
\end{equation}

\item[$(*)$] The equality precisely means $t=\frac{\log(\mu_{k+1})}{\log(\mu_{j_0+d})}$. However, the expression on the right-hand side can only take countable many values whereas $t$ is required to belong to an uncountable set and therefore \eqref{starviolating} is impossible.
\end{itemize}
\qed\enddemo

We finish with the following consequence showing that in general \cite[Thm. 5.1]{orliczbook} does not provide a characterization. This should be compared with \cite[Sect. 2.3, Thm. 11]{orliczbook1} where it has been shown that an N-function $F$ satisfies $\Delta'$ if and only if the left-derivative of $F$ has $\Delta'$.

\begin{corollary}
There exist N-functions $F$ such that $F$ satisfies the $\Delta'$-condition but $(\Delta'_{f})$ fails.
\end{corollary}

\demo{Proof}
Consider the function(s) $\omega_s$, $s>1$. Then $\varphi_{\omega_s}(t)=t^s$ for all $t\ge 0$ and so the $\Delta'$-condition is satisfied. Since $\omega_{M^{q,n}}\sim\omega_s$ for all $q>1$ and $n>1$ such that $\frac{1}{s}+\frac{1}{n}=1$ we get that $\varphi_{\omega_{M^{q,n}}}\sim\varphi_{\omega_s}$ as well (see the proof of $(ii)\Leftrightarrow(iii)$ of Theorem \ref{comparisonprop0}). It is immediate that the $\Delta'$-condition is also preserved under $\sim$. Hence $\varphi_{\omega_{M^{q,n}}}$ and finally $F_{M^{q,n}}$ satisfy the $\Delta'$-condition. However, for any $q>1$ the corresponding sequence of quotients is clearly strictly increasing (recall that $\mu^{q,n}_j=q^{j^n-(j-1)^n}$, $j\ge 1$) and so by Proposition \ref{Deltaprimefails} property $(\Delta'_{f_{M^{q,n}}})$ fails.
\qed\enddemo

\bibliographystyle{plain}
\bibliography{Bibliography}

\end{document}